\documentstyle[12pt]{article}
\begin{document}
\author{S.V. Ludkovsky.}
\title{Line antiderivations over local fields and their applications.}
\date{25 November 2003}
\maketitle
\begin{abstract}
A non-Archimedean antiderivational line analog of the Cauchy-type
line integration is defined
and investigated over local fields. Classes of non-Archimedean
holomorphic functions are defined and studied. Residues of functions
are studied, Lorent series representations are described.
Moreover, non-Archimedean antiderivational
analogs of integral representations
of functions and differential forms such as the Cauchy-Green,
Martinelli-Bochner, Leray, Koppelman and Koppelman-Leray formulas
are investigated. Applications to manifold and operator theories
are studied.
\end{abstract}
\section{Introduction}
Line (Cauchy) integration is the cornerstone in the complex analysis
and integral formulas of functions and differential forms such as
the Cauchy-Green, Martinelli-Bochner, Leray, Koppelman and
Koppelman-Leray formulas play very important role in it and in
analysis on complex manifolds and theory of Stein and K\"ahler
manifolds and theory of holomorphic functions (see, for example,
\cite{henlei,shabat}).
In the non-Archimedean case there is not so developed analog
of complex analysis. Though there are few works devoted to
non-Archimedean holomorphic functions over the complex non-Archimedean
field $\bf C_p$ and the Levi-Civit\'a fields, which are not locally compact
(see \cite{koblbr,berz} and references therein). In that works
M.M. Vishik and M. Berz have obtained analogs of residues and the
Cauchy formula, but integrals that they have used were of
combinatorial-algebraic nature and they have operated with power series
mainly for their analogs of holomorphic functions. On the other hand,
there is not any
measure equaivalent with the Haar measure on such nonlocally compact fields
because of the A. Weil \cite{weil40} theorem stating that the existence
of such nontrivial measure on a topological group implies its local
compactness. This article is devoted to others non-Archimedean analogs
of integral representation theorems, that were not yet considered by others
authors. Moreover, this article operates with locally compact non-Archimedean
fields of characteristic zero (local fields) and the corresponding
analogs of complex planes. Apart from the classical case in the
non-Archimedean case there is not any indefinite integral. Instead of it
antiderivation operators by Schikhof \cite{sch1} are used.
\par It is necessary to note that in this article are considered
not only manifolds treated by the rigid geometry, but
much wider classes continuing the previous work \cite{lustpnam}.
For them the existence of an exponential
mapping is proved. A rigid non-Archimedean geometry
serves mainly for needs of the cohomology theory on such manifolds,
but it is too restrictive and operates with narrow classes of analytic
functions \cite{freput}. It was introduced at the beginning
of sixties of the 20-th century. Few years later wider classes
of functions were investigated by Schikhof \cite{sch1}.
In this paper classes of functions and antiderivation operators
by Schikhof and their generalizations from works \cite{lutmf99,luambp00}
are used.
\par Section 2 is devoted to the definition and investigations of the
non-Archimedean analogs of the line integration over local fields.
Classes of non-Archimedean holomorphic functions are defined and studied.
For this specific non-Archimedean geometry's definitions
and theorems are given (see also definitions and notations in
\cite{luambp00,lutmf99,lustpnam,lujmsq2,lu8}).
It is necessary to note that definitions, formulations of theorems,
propositions, etc. and their proofs differ substantially from the classical
case (over $\bf C$). Residues of functions
are studied, Lorent series representations are described.
In Section 3 non-Archimedean antiderivational
analogs of integral representations
of functions and differential forms such as the Cauchy-Green,
Martinelli-Bochner, Leray, Koppelman and Koppelman-Leray formulas
are investigated. These studied are accomplished on domains
in finite dimensional Banach spaces over local fileds and also
on manifolds over local fields. All results of this paper are
obtained for the first time. Finally applications
of the obtained results to the theory of non-Archimedean manifolds
and linear operators in non-Archimedean Banach spaces are outlined.
In works of Vishik (see \cite{koblbr} and references therein)
the theory of non-Archimedean (Krasner) analytic operators
with compact spectra in $\bf C_p$ was developed. In this article
operators may have noncompact spectra in a field $\bf L$
such that ${\bf Q_p}\subset \bf L$ (may be also ${\bf L}\supset
\bf C_p$ and ${\bf L}\ne \bf C_p$) continuing the investigation
of \cite{luddia}.

\section{Line antiderivation over local fields}
\par To avoid misunderstandings we first present our
specific definitions.
\par {\bf 2.1. Notation and Remarks.} Let $\bf K$ denotes a local field,
that is, a finite algebraic extension of the field $\bf Q_p$
of $p$-adic numbers with a norm extending that of $\bf Q_p$
\cite{wei}. Denote by $\bf C_p$ the field of complex numbers
with the norm extending that of $\bf Q_p$ \cite{kobl}.
If $i\in \bf K$ take $\alpha \in {\bf C_p}\setminus \bf K$
such that there exists ${\tilde m}\in \bf N$ with $\alpha ^{\tilde m}
\in \bf K$, where $\tilde m$ is such a minimal natural number,
${\tilde m}={\tilde m}(\alpha )$, $i:=(-1)^{1/2}$.
If $i\notin \bf K$ take $\alpha =i$. Denote by ${\bf K}(\alpha )$
a local field which is the extension of $\bf K$ with the help
of $\alpha $.
\par Suppose $U$ is a clopen compact perfect (that is, dense in itself)
subset in $\bf K$ and $\mbox{ }_U\sigma :=\sigma $ is its
approximation of the identity:
there is a sequence of maps $\sigma _l: U \to U$, where $0\le l\in \bf Z$,
such that
\par $(i)$ $\sigma _0$ is constant;
\par $(ii)$ $\sigma _l\circ \sigma _n=\sigma _n\circ \sigma _l=
\sigma _n$ for each $l\ge n$;
\par $(iii)$ there exists a constant $0< \rho <1$ such that
for each $x, y \in U$ the inequality $|x-y| < \rho ^n$ implies
$\sigma _n(x)=\sigma _n(y)$;
\par $(iv)$ $|\sigma _n(x)-x| <\rho ^n$ for each integer $n\ge 0$.
Consider spaces $C^n(U,{\bf L})$ of all $n$-times continuously differentiable
in the sence of difference quotients functions $f: U\to \bf L$,
where $\bf L$ is a field containing $\bf K$ with the multiplicative norm
$|*|_{\bf L}$ which is the extension of the multiplicative norm
$|*|_{\bf K}$ in $\bf K$.
Then there exists an antifderivation:
\par $(1)$ $\mbox{ }_UP^n: C^{n-1}(U,{\bf L})\to C^n(U,{\bf L})$ given by the formula:
\par $(2)$ $\mbox{ }_UP^nf(x):=\sum_{l=0}^{\infty }\sum_{j=0}^{n-1}
f^{(j)}(x_l)(x_{l+1}-x_l)^{j+1}/(j+1)!$, \\
where $x_l:=\sigma _l(x)$, $x\in U$, $n\ge 1$ (see \S 80 \cite{sch1}).
Formula $(2)$ shows, that if $\mbox{ }_UP^n$ is defined on
$C^{n-1}(U,{\bf K})$,
then it is defined on $C^{n-1}(U,Y)$ for each field $\bf L$
which is complete relative to its norm such that ${\bf K}\subset \bf L$
and a Banach space $Y$ over $\bf L$.
\par Since $P^n$ is the $\bf L$-linear operator, then there exists
the $\bf L$-linear space $\mbox{ }_PC^n_0(U,Y):=
P^n(C^{n-1}(U,Y)) $, put $\mbox{ }_PC^n(U,Y):=
\mbox{ }_PC^n_0(U,Y) \oplus Y$, where $n\ge 1$,
$Y$ is a Banach space over $\bf L$.
For a clopen subset $\Omega $ in $({\bf K}\oplus \alpha {\bf K})^m$
such that $\Omega \subset U^m\times U^m$
consider the antiderivation $\mbox{ }_{\Omega }P^nf(z)$
as the restriction of $\mbox{ }_{U^m\times U^m}P^nf(z)$ on $\Omega $,
\par $(3)\quad \mbox{ }_{\Omega }P^nf(z):=
\mbox{ }_{U^m\times U^m}P^n|_{\Omega }f(z)=
\mbox{ }_{U^m\times U^m}P^nf(z)\chi _{\Omega }(z)$, where
\par $(4)\quad \mbox{ }_{U^m\times U^m}P^nf(z):=\mbox{ }_UP^n_{x_1}...
\mbox{ }_UP^n_{x_m} \mbox{ }_UP^n_{y_1}...\mbox{ }_UP^n_{y_m}f(z)$, \\
$\chi _{\Omega }(z)$ denotes the characteristic function of $\Omega $,
$\chi _{\Omega }(z)=1$ for each $z\in \Omega $, $\chi _{\Omega }(z)=0$
for each $z\in {\bf K}^{2m}\setminus \Omega $,
$z=(x,y)$, $x, y\in U^m\subset {\bf K}^m$,
$x=(x_1,...,x_m)$, $x_1,...,x_m\in \bf K$, $\mbox{ }_UP^n_{x_l}$ means the
antiderivation by the variable $x_l$. This is correct, since each
$f\in C^{(0,n-1)}(\Omega ,{\bf L}):=C((0,n-1),\Omega \to {\bf L})$
(see \S I.2.4 \cite{luambp00}) has a $C^{(0,n-1)}$-extension on
$U^m\times U^m$, for example, $f|_{U^m\times U^m\setminus \Omega }=0$.
This means, that
$\mbox{ }_{U^m\times U^m}P^nf(z)$ is the antiderivation defined with
the help approximation of the unity on $U^m\times U^m$ such that
$\mbox{ }_{U^m\times U^m}\sigma =(\mbox{ }_U\sigma ,...,\mbox{ }_U \sigma )$.
\par The condition of compactness of $\Omega $ is not very restrictive,
since each locally compact subset in $({\bf K}\oplus \alpha {\bf K})^m$
has a one-point (Alexandroff) compactification which is totally disconnected
and hence homeomorphic to a clopen subset in
$({\bf K}\oplus \alpha {\bf K})^m$ (see \S 3.5 and Theorem 6.2.16
about universality of the Cantor cube in \cite{eng}).
If $\rho (z_1,z_2):=|z_1-z_2|$ is the metric in
$({\bf K}\oplus \alpha {\bf K})^m$, then the metric
$\rho '(z_1,z_2):=\rho (z_1,z_2)/[1+\rho (z_1,z_2)]$ has the extension
on the one-point compactification $A({\bf K}\oplus \alpha {\bf K})^m:=
({\bf K}\oplus \alpha {\bf K})^m\cup \{ A \} $, where $A$ is a singleton.
If $Y$ is a mertic space with a metric $\rho $, then
$B(Y,y,r):=\{ z\in Y: \rho (z,y)\le r \} $ denotes the ball
of radius $r>0$ and containing a point $y\in Y$.
\par {\bf 2.2. Notes and Definitions. 1.}
For a local field $\bf K$ there exists a prime $p$ such that
$\bf K$ is a finite algebraic extension of $\bf Q_p$.
In view of Theorems 1.1, 4.6 and Proposition 4.4 \cite{wei} there exists
a prime element $\pi \in \bf K$ such that $P=\pi R=R\pi $,
$R/P$ is a finite field 
$\bf F_{p^n}$ consisting of $p^n$ elements for some $n\in \bf N$
\cite{wei}, $mod_{\bf K}(\pi ):=q^{-1}$ and $\Gamma _{\bf K}:=mod ({\bf K})$,
where $mod_{\bf K}$ is the modular function of $\bf K$ associated
with the nonnegative Haar measure $\mu $ on $\bf K$ such that
$\mu (xS)=mod_{\bf K}(x)\mu (S)$ for each $0\ne x\in \bf K$,
$mod_{\bf K}(0):=0$ and
each Borel subset $S$ in $\bf K$ with $\mu (S)<\infty $,
$P:=\{ x\in {\bf K}: |x|<1 \} $, $R:=B({\bf K},0,1)$.
Then each $x\in \bf K$ can be written in the form
$x=\sum_lx_l\pi ^l$, where $x_l\in \{ 0,\theta _1,...,\theta _{p^n-1} \} $,
$\min_{x_l\ne 0}l=:-ord_{\bf K}(x)>-\infty $,
$\theta _0+P$, $\theta _1+P$,...,$\theta _{p^n-1}+P$ is the disjoint covering
of $R$, $\theta _0:=0$.
Consider in $\bf K$ the linear ordering $a\triangle b$ if
$a_k=b_k$,..., $a_s=b_s$, $a_{s+1}<b_{s+1}$, where
$a$, $b\in \bf K$, by our definition $\theta _s<\theta _v$ for each
$s<v$, $k:= \min (ord_{\bf K}(a),ord_{\bf K}(b))$.
In $B({\bf K},0,1)$ the largest element relative to such linear
ordering is $\beta :=\sum_{l=0}^{\infty }\theta _{(p^n-1)}\pi ^l=
\theta _{(p^n-1)}/(1-\pi )$.
\par Though this linear ordering is preserved
neither by additive nor by multiplicative structures of $\bf K$
it is useful (see, for example, \cite{put68} and \S 62 \cite{sch1}).
\par {\bf 2.2.2.} Let $v_0,...,v_k\in {\bf K}(\alpha )^m$
such that vectors $v_1-v_0$,...,$v_k-v_0$ are $\bf K$-linearly independent,
then the subset $s:=[v_0,...,v_k]:=\{ z\in {\bf K}(\alpha )^m:
z=a_0v_0+...+a_kv_k; a_0+...+a_k=1; a_0,...,a_k\in
B({\bf K},0,1) \} $ is called the simplex of dimension $k$
over $\bf K$, $k=dim_{\bf K}s$. A polyhedron $P$
is by our definition the union of a locally finite
family $\Psi _P$ of simplexes. For compact $P$ a family
$\Psi _P$ can be chosen finite. An oriented $k$-dimensional
simplex is a simplex together with a class of linear orderings
of its vertices $v_0,...,v_k$. Two linear orderings are equivalent
if they differ on an even transposition of vertices.
For a simplicial complex $S$ let $C_q(S)$ be an Abelian group
generated by simplices $s^q$ of dimension $q$ over $\bf K$
and relations $s_1^q+s_2^q=0$, if $s_1^q$ and $s_2^q$ are differently
oriented simplices (see the real case in Chapter 4 \cite{span}).
Then there exists the homomorphism
$\partial_q: C_q(S)\to C_{q-1}(S)$ such that
$\partial_q [v_0,...,v_q]:=\sum_{l=0}^q (-1)^l[v_0,...,v_{l-1},
v_{l+1},...,v_q]$ and $\partial_q [v_0,...,v_q]$ is called the oriented
$\bf K$-boundary of $s^q$.
\par {\bf 2.2.3.} A clopen compact subset
$\Omega $ in $({\bf K}\oplus \alpha {\bf K})^m$ is totally
disconnected and its topological boundary is empty.
Nevertheless, using the following affine construction
it is possible to introduce convention about certain curves
and boundaries which will serve for the antiderivation operators.
\par Let $\Omega $ be a locally $\bf K$-convex subset
in ${\bf K}(\alpha )^m$ for which there exists a sequence
$\Omega _n$ of polyhedra with $\Omega _n\subset \Omega _{n+1}$ for
each $n\in \bf N$, $\Omega =cl (\bigcup_n \Omega _n)$,
where $cl (S)$ denotes the closure of a subset
$S$ in ${\bf K}(\alpha )^m$. Suppose each $\Omega _n$ is the union
of simplices $s_{j,n}$ with vertices $v^j_{0,n},...,v^j_{k,n}$,
$j=1,...,b(n)\in \bf N$, moreover, $dim_{\bf K}(s_{j,n}\cap s_{j',n})< k$
for each $j\ne j'$ and each $n$, where $k>0$ is fixed.
Then define the oriented $\bf K$-border
$\partial \Omega _n:=\sum_{j,l}(-1)^l[v^j_{0,n},...,
v^j_{l-1,n},v^j_{l+1,n},...,v^j_{k,n}]$.
Consider $\Omega _n$ for each $n$ such that if
$dim_{\bf K}(s_{j,n}\cap s_{j',n})=k-1$ for some $j\ne j'$,
then $s_{j,n}\cap s_{j',n}=[v^j_{0,n},...,v^j_{l-1,n},v^j_{l+1,n},...,
v^j_{k,n}]=[v^{j'}_{0,n},...,v^{j'}_{l'-1,n},v^{j'}_{l'+1,n},...,
v^{j'}_{k,n}]$ and $(l-l')$ is odd.
For each $n$ choose a set of vertices
generating $\Omega _n$ of minimal cardinality and
such that the sequence $\{ \partial \Omega _n: n \} $ converges relative to
the distance function $d(S,B):=\max (\sup_{x\in S} \rho (x,B),
\sup_{b\in B} \rho (b,S))$, where $\rho (x,B):=\inf_{b\in B} \rho (x,b)$
and $\rho (x,b):=|x-b|$. Then by our definition
$\partial \Omega :=\lim_{n\to \infty }\partial \Omega _n$.
\par Evidently, each clopen compact subset $\Omega $ has such decomposition
into simplices and the described $\partial \Omega $, since
$\Omega $ is the finite union of balls, but for two balls
$B_1$ and $B_2$ in $({\bf K}\oplus \alpha {\bf K})^m$ either
$B_1\subset B_2$ or $B_2\subset B_1$ or $B_1\cap B_2=\emptyset $
due to the ultrametric inequality and each ball $B$ has such decomposition
into simplices as described above.
\par {\bf 2.2.4.}
We say that a subset $\Omega $ in $A({\bf K}\oplus \alpha {\bf K})^m$
encompasses a point $z$ if $z\in \Omega $.
\par For the unit ball relative to the metric $\rho (z_1,z_2):=|z_1-z_2|$
let its non-Archimedean canonical
oriented $\bf K$-border $\partial_cB_{\rho }({\bf K}\oplus
\alpha {\bf K},0,1)$
be given by the set $[(-\beta ,-\beta ),(\beta ,-\beta )] \cup
[(\beta ,-\beta ),(\beta ,\beta )] \cup [(\beta ,\beta ),
(-\beta ,\beta )] \cup [(-\beta ,\beta ),(-\beta ,-\beta )]$, where
$[a,b]:= \{ z\in {\bf K}\oplus \alpha {\bf K}: z=(1-t/\beta )a+(t/\beta )b,
t\in B({\bf K},0,1) \} $ for each $a, b\in {\bf K}\oplus \alpha {\bf K}$.
Then $\partial_c B(({\bf K}\oplus \alpha {\bf K})^m,0,1):=
\bigcup_{l=1}^m B({\bf K}\oplus \alpha {\bf K},0,1)^{l-1}\times
\partial_cB({\bf K}\oplus \alpha {\bf K},0,1) \times
B({\bf K}\oplus \alpha {\bf K},0,1)^{m-l}$,
$\partial_c B(({\bf K}\oplus \alpha {\bf K})^m,z,q^k):=
z+\pi ^{-k}\partial_cB(({\bf K}\oplus \alpha {\bf K})^m,0,1)$.
This is the particular case of \S 2.2.
\par A continuous mapping $\gamma : B({\bf K},0,1)\to A({\bf K}(\alpha ))^m$
is called a path. We say that $\gamma $ encompass a point
$z\in A({\bf K}(\alpha ))^m$, if
\par $(i)$ $z\in \Omega $, where $\partial \Omega = \gamma $,
$dim_{\bf K}\Omega =2$,
\par $(ii)$ $z\notin \gamma (B({\bf K},0,1))$,
\par $(iii)$ $|z|< \sup_{\theta \in B({\bf K},0,1)} |\gamma (\theta )|$
for $z\ne A$, $\sup_{\theta \in B({\bf K},0,1)} |\gamma (\theta )|
<\infty $ for $z=A$.
\par A path $\gamma $ we call a locally affine, if there exists
a finite partition $\cal Z$ of $\gamma (B({\bf K},0,1))$ such that
$\gamma =\bigcup_{l=1}^n\tau _l$, where ${\cal Z}:= \{ z_0,z_1,...,z_n \} $,
$\tau _l:=[z_{l-1},z_l]$ for each $l=1,...,n$.
We consider the family ${\cal F}_q$ of all paths $\gamma $ for which
there exists a sequence $\{ \gamma _n: n \} \subset {\cal F}_a$
converging relative to the distance function $d'(S,B):=
\max (\sup_{x\in S}\rho '(x,B), \sup_{b\in B}\rho '(b,S))$
to $\gamma $ in $(A({\bf K}(\alpha ))^m, \rho ')$
and such that there there exists a homeomorphism $\nu $
of $\gamma (B({\bf K},0,1))$ with $B({\bf K},0,1)$ and $\nu $ is a
piecewise $\mbox{ }_PC^{q+1}$-diffeomorphism
with it, where ${\cal F}_a$ denotes the family of all locally affine paths,
$q\in \bf N$.
In addition we take $\Omega $ and $\gamma $
such that $\gamma =\partial \Omega $ in accordance with \S 2.2.
\par Since $A{\bf K}(\alpha )^m$ and $A({\bf K}\oplus \alpha
{\bf K})^m$ are compact, then a clopen compact set $\Omega $
in $A{\bf K}(\alpha )^m$ or in $A({\bf K}\oplus \alpha {\bf K})^m$
is homeomorphic with a clopen compact subset $\kappa (\Omega )$
in ${\bf K}(\alpha )^m$ or $({\bf K}\oplus \alpha {\bf K})^m$
respectively (see Theorem 6.2.16 and Corollary 6.2.17 about
universality of the Cantor cube for zero dimensional spaces
\cite{eng}), where $\kappa : \Omega \to
\kappa (\Omega )$ is the homeomorphism. Therefore, we can consider
$\mbox{ }_{\Omega }P^n$, $\partial \Omega $ and $\mbox{ }_{\partial
\Omega }P^n$ induced by $\kappa $ of such sets $\Omega $ also.
\par {\bf 2.2.5.} Let $M$ be a $C^{\xi + (1,0)}$-manifold
of dimension $k$ over $\bf K$ such that $\xi =(q,n-1)$,
where spaces $C^{\xi }({\bf K}^a,{\bf K}^b):=
C(\xi ,{\bf K}^a\to {\bf K}^b)$ and $C^{\xi }$-manifolds
and uniform spaces $C^{\xi }(M,N)$ of all $C^{\xi }$-mappings
$f: M\to N$  were defined
in \S I.2.4 \cite{luambp00}, $0\le q \in \bf Z$, $0< n\in \bf Z$,
$\mbox{ }_PC^{\xi +(0,1)}_0(\Omega ,{\bf L}^b):=P^n(C^{\xi }
(\Omega ,{\bf L}^b))$ and
$\mbox{ }_PC^{\xi +(0,1)}(\Omega ,{\bf L}^b):=
\mbox{ }_PC^{\xi +(0,1)}_0 (\Omega ,{\bf L}^b)\oplus {\bf L}^b$
was described in Lemma 3.4 \cite{lustpnam}.
Suppose that charts $(V_j,\phi _j)$ of the atlas $At(M)$ of $M$
are such that $V_j$ are clopen in $M$, $\bigcup_jV_j=M$,
$\phi _j: V_j\to \phi _j(V_j)\subset U^k$ are homeomorphisms
on clopen subsets in $U^k$, where $\phi _{i,j}:=\phi _i\circ \phi _j^{-1}
\in \mbox{ }_{P,x_l}C^{\xi +(1,0)}\cap C^{\xi }(W_{i,j},{\bf K}^k)$
for each $i\ne j$ with $U_i\cap U_j\ne \emptyset $
and each coordinate $x_l$ induced from $\bf K$, $l=1,...,k$,
$W_{i,j}:=dom (\phi _{i,j})$, $\mbox{ }_{P,x_l}C^{\xi +(1,0)}
(\Omega ,Y):=\mbox{ }_UP^n_{x_l}(C^{\xi }(\Omega ,Y))\oplus Y$
for a Banach space $Y$ over $\bf L$, ${\bf K}\subset \bf L$.
Then we call such $M$ the $\mbox{ }_SC^{\xi +(1,0)}$-manifold,
such mappings $\phi _{i,j}$ we call the $\mbox{ }_SC^{\xi +(1,0)}$-mappings.
Then $\mbox{ }_SC^{\xi +(1,0)}(\Omega ,Y):=
\{ f \in C^{\xi +(1,0)}(\Omega ,Y):$
$f(x_1,...,x_k)\in \mbox{ }_{P,x_l}C^{\xi +(1,0)}(\Omega ,Y)$
for each $l=1,...,k \} $, where $Y$ is a Banach space over $\bf L$.
In particular, $\mbox{ }_SC^{\xi +(1,0)}(U,Y)=
\mbox{ }_PC^{\xi +(0,1)}(U,Y) =\mbox{ }_PC^{q+n}(U,Y)$,
but for $dim_{\bf K}\Omega >1$ these spaces are different
$\mbox{ }_SC^{\xi +(1,0)}(\Omega ,Y)\ne \mbox{ }_PC^{\xi +(0,1)}(\Omega ,Y)$.
Tensor fields over $M$ were defined in \S \S 3.1 and 3.5 \cite{lustpnam}.
Then the bundle of $r$-differential forms is the antisymmetrized
bundle $\psi _r: \Lambda ^rM\to M$ of the bundle $\tau _r:
T_r M\to M$ of $r$-fold covariant tensors.
\par Consider the $\mbox{ }_SC^{\xi +(1,0)}$-diffeomorphism
$\phi : \tau \to \phi (\tau )$,
that is, $\phi $ is surjective and bijective with $\phi $ and
$\phi ^{-1}\in \mbox{ }_SC^{\xi +(1,0)}$, where
\par $(i)\quad \tau =[v_0,v_1]\times [v_1,v_2]\times ... [v_{k-1},v_k]$ \\
is the parallelepiped in ${\bf K}^k$,
vectors $v_1-v_0,...,v_k-v_0$ are $\bf K$-linearly independent.
Then for a $k$-differential $C^{(0,n-1)}$-form $w$ on $\phi (\tau )$ define
\par $(1)\quad \mbox{ }_{\phi (\tau )}P^nw:=\mbox{ }_{\tau }P^n\phi ^*w$, \\
where $\phi ^*w$ is the pull back of $w$ such that 
\par $(2)\quad \mbox{ }_{\phi (\tau )}P^nw=0$ for $dim_{\bf K}\tau \ne k$, \\
since $w=0$ for $k>dim_{\bf K}M$. Without loss of generality
take $0\in U$ and $\sigma _0(0)=0$, then $\sigma _l(0)=0$ for each
$l\in \bf N$, consequently, $\mbox{ }_UP^n|_{ \{ 0 \} } =0$.
Therefore, $\mbox{ }_{U^m}P^n|_{(U^m\cap {\bf K}^k\times \{ 0 \} ^{m-k}
)} w=0$ for $k< dim_{\bf K} \Omega =m$.
Each such parallelepiped is the finite union of simplices
satisfying conditions of \S 2.2.3.
The orientation of $\partial \tau $ is induced by the orientations
of constituting its simplices which are consistent.
Consider such parallelepipeds
$\tau _{j,q,l}$ with $l=1,...,b(q)\in \bf N$ and
\par $(ii)\quad \dim_{\bf K}(\tau _{j,q,l}\cap \tau _{j,q,l'})<k$ 
for each $l\ne l'$ and 
\par $(iii)\quad cl (\bigcup_q\kappa _{j,q})=\phi _j(V_j)$, where 
\par $(iv)\quad \bigcup_{l=1}^{b(q)}\tau _{j,q,l}=:\kappa _{j,q}$, 
\par $(v)\quad \lim_{q\to \infty }\max_l diam (\tau _{j,q,l})=0$.  \\
Since $\mbox{ }_{\tau _{j,q,l}}P^nv+\mbox{ }_{\tau _{j,q,l'}}P^nv=
\mbox{ }_{\tau _{j,q,l}\cup \tau _{j,q,l'}}P^nv$
for each differential $C^{(0,n-1)}-k$-form $v$ with support in
$U^k$ and each $l\ne l'$ and $\mbox{ }_{U^k}P^n$
is the continuous operator from
$C^{(0,n-1)}(U^k,{\bf L})$ to $C^{(0,n)}(U^k,{\bf L})$
then there exists
\par $(vi)\quad \lim_{q\to \infty }
\sum_{l=1}^{b(q)} \mbox{ }_{\tau _{j,q,l}}P^nv=:
\mbox{ }_{\phi _j(V_j)}P^nv.$  \\
Using transition mappings $\phi _{i,j}$ and considering
clopen disjoint covering
\par $(vii)\quad W_j:=V_j\setminus \bigcup_{l=1}^{j-1}V_j$ of $M$ we get
\par $(viii)\quad \mbox{ }_MP^nw=\sum_j\mbox{ }_{W_j}P^nw$ \\
independent on the choise of local coordinates in $M$.
Mention, that since $|\beta |=1$, then $B({\bf K}^l,z,r)$
can be represented as the parallelepiped with the desribed above
$\bf K$-boundary $\partial_c B({\bf K}^l,z,r)$ due the ultrametric inequality.
Due to $(vi-viii)$ there is defined $\mbox{ }_{\gamma }P^nv$
for locally affine path $\gamma $, which is the
$\mbox{ }_PC^n$-manifold, that will be supposed henceforth.
\par Each compact manifold $M$ has a finite dimension over $\bf K$
and using $W_j$ we get an embedding into ${\bf K}^b$ for some
$b\in \bf K$. Let $\phi : \Omega \to M$ be
such that $\phi $ is surjective and bijective,
$\phi $ and $\phi ^{-1}\in \mbox{ }_SC^{\xi +(1,0)}$,
which means that $\phi _j\circ \phi \in \mbox{ }_SC^{\xi +(1,0)}
(\phi ^{-1}(V_j),{\bf K}^k)$
and $\phi ^{-1}\circ \phi _j^{-1}\in \mbox{ }_SC^{\xi +(1,0)}
(\phi _j(V_j),{\bf K}^k)$ for each $j$, where
$\phi ^{-1}(M)=\Omega \subset U^k$ satisfies conditions of \S 2.2.3.
Such $\phi $ we call the $\mbox{ }_SC^{\xi +(1,0)}$-diffeomorphism.
Then $M$ is oriented together with $\Omega $.
Then $\partial M:=\phi (\partial \Omega )$
is the oriented boundary. We also can consider the analytic manifold
$M$ and the analytic diffeomorphism $\phi $.
Each compact $C^{\xi }$-manifold $M$ can be supplied
with the analytic manifold structure using a disjoint covering
refined into $At(M)$.
\par {\bf 2.2.6. Theorem. } {\it Let $M$ be a compact
$\mbox{ }_SC^{\xi }$ or $\mbox{ }_PC^{\xi }$-manifold
over the local field $\bf K$ with dimension $dim_{\bf K}M=k$
and an atlas $At (M)=\{ (V_j, \phi _j): j=1,..., n \} $,
where $\xi =(q,n)$, $1\le q\in \bf N$, $0\le n\in \bf Z$,
then there exists a $\mbox{ }_SC^{\xi }$ or $\mbox{ }_PC^{\xi }$-embedding
of $M$ into ${\bf K}^{nk}$ respectively.}
\par {\bf Proof.} Let $(V_j,\phi _j)$ be the chart of the atlas
$At(M)$, where $V_j$ is clopen in $M$, hence $M\setminus V_j$
is clopen in $M$. Therefore, there exists a $\mbox{ }_SC^{\xi }$
or $\mbox{ }_PC^{\xi }$-mapping
$\psi _j$ of $M$ into ${\bf K}^k$ such that $\psi _j(M\setminus V_j)
= \{ x_j \} $ is the singleton and $\psi _j: V_j\to \psi _j(V_j)$
is the $\mbox{ }_SC^{\xi }$ or $\mbox{ }_PC^{\xi }$-diffeomorphism
onto the clopen subset $\psi _j(V_j)$ in ${\bf K}^k$ correspondingly,
$x_j\in {\bf K}^k\setminus \psi _j(V_j)$,
since the operator $\mbox{ }_MP^n$ is $\bf K$-linear,
$\mbox{ }_MP^n0=0$ and the covering
$\{ V_j: j \} $ of $M$ has a disjoint finite refinement
$\{ W_k: k \} $ such that $P^n_{x_l}[f]=
P^n_{x_l}[\sum_kf\chi _{W_k}]=\sum_kP^n_{x_l}[f\chi _{W_k}]$
for each $f\in C^{(q,n-1)}(M,{\bf K})$
and each coordinate $x_l$ (see \S 2.1 and \S 2.2.5).
Then the mapping $\psi (z):=(\psi _1(z),...,\psi _n(z))$
is the embedding into ${\bf K}^{nk}$, since the rank
$rank [d_z\psi (z)]=k$ at each point $z\in M$, because
$rank [d_z\psi _j(z)]=k$ for each $z\in V_j$ and $dim_{\bf K}\psi (V_j)\le
dim_{\bf K}M=k$. Moreover, $\psi (z)\ne \psi (y)$ for each
$z\ne y\in V_j$, since $\psi _j(z)\ne \psi _j(y)$.
If $z\in V_j$ and $y\in M\setminus V_j$, then there exists
$l\ne j$ such that $y\in V_l\setminus V_j$,
$\psi _j(z)\ne \psi _j(y)=x_j$.
\par {\bf 2.3.1. Theorem.} {\it Let $M$ be a compact oriented manifold
over $\bf K$ of dimension $\dim_{\bf K}M=k>0$ with an oriented
boundary $\partial M$
and let $w$ be a differential $(k-1)$-form as in \S 2.2.5
such that its pull back $\phi ^*w$ is a differential $(k-1)$
$\mbox{ }_SC^{(1,n-1)}$-form, then \\
$(1)\quad \mbox{ }_MP^ndw=\mbox{ }_{\partial M}P^nw$.}
\par {\bf Proof.} Since $M$ is the manifold of $dim_{\bf K}M=k>0$,
then $M$ is dense in itself and compact, hence $\Omega $ is dense
in itself and compact (see Chapter 1 and Theorems 3.1.2, 3.1.10 \cite{eng})
and the approximation of the identity can be applied
to $\Omega $. In view of Formulas $2.1.(1-4)$ and $2.5.(1,2)$
on the space of $C^{\xi }$
differential forms operators $\mbox{ }_UP^n_{x_q}$
and $\mbox{ }_UP^n_{x_s}$ commute for each $1\le q, s \le k$. Then
\par $(i)\quad \mbox{ }_UP^nf|^b_a=-\mbox{ }_UP^nf|^a_b$, where
$\mbox{ }_UP^nf|^b_a:=\mbox{ }_UP^nf(b)-\mbox{ }_UP^nf(a)$.
Using conditions imposed on the manifold $M$, partitions
of $\Omega _n$ into unions of parallelepipeds, which are finite
unions of simplices as in \S 2, Formula $(i)$ and $2.5.(1)$,
also using the limit $2.5.(vi)$ and Formula $2.5.(viii)$,
it is sufficient to verify Formula $(1)$ for a parallelepiped
and an arbitrary term $\psi :=f(z)dz_1\wedge ... \wedge dz_{q-1}\wedge
dz_{q+1}\wedge ... \wedge dz_k$ corresponding
to the differential $(k-1)$ $\mbox{ }_SC^{(1,n-1)}$-form $\phi ^*w$.
Consider in ${\bf K}^k$ the standard orthonormal base
$e_1,...,e_k$, where $e_l:=(0,...,0,1,0,...,0)$ is the vector
with $1$ in $l$-th place. Without loss of generality using limits
we can take the parallelepipeds $\tau =[v_0,v_1]\times ...\times
[v_{k-1},v_k]$ with $v_l-v_{l-1}=\lambda _le_l$ for each $l=1,..,k$, where
$0\ne \lambda _l\in \bf K$.
Therefore, $df(z)=(-1)^{q-1}(\partial f(z)/\partial z_q)
dz_1\wedge ... \wedge dz_k$. Since $f\in \mbox{ }_SC^{(1,n-1)}$,
then $\mbox{ }_UP^n_{z_q}(\partial f(z)/\partial z_q)|^b_a=
f(z_1,...,z_{q-1},b,z_{q+1},...,z_k)-f(z_1,...,z_{q-1},a,z_{q+1},...,z_k)$
for each $l=1,...,k$. Consequently, \\
$(ii)\quad \mbox{ }_{\tau }P^nd\psi =
(-1)^{q-1}\mbox{ }_{[v_0,v_1]\times ...[v_{q-2},v_{q-1}]
\times [v_{q+1},v_{q+2}]\times ... [v_{k-1},v_k]} P^ndz_1\wedge ...
\wedge dz_{q-1} \wedge dz_{q+1} \wedge ... \wedge dz_k
\mbox{ }_{[v_{q-1},v_q]}P^n(\partial f(z)/\partial z_q)dz_q$ \\
$=(-1)^{q-1} \mbox{ }_{[v_0,v_1]\times ...[v_{q-2},v_{q-1}]
\times [v_{q+1},v_{q+2}]\times ... [v_{k-1},v_k]}P^n \{
f(z_1,...,z_{q-1},v_q,dz_{q+1},...,dz_k)-
f(z_1,...,z_{q-1},v_{q-1},dz_{q+1},...,dz_k) \}
dz_1\wedge ... \wedge dz_{q-1} \wedge dz_{q+1} \wedge ... \wedge dz_k$ \\
for each $q=1,...,k$. In view of $2.5.(2)$ antiderivations
of $\psi $ by others pieces $(-1)^{s-1}[v_0,v_1]\times ...
[v_{s-2},v_{s-1}] \times ( \{ v_s \} - \{ v_{s-1} \} ) \times
[v_s,v_{s+1}] \times ...
[v_{k-1},v_k]$ corresponding to $s\ne q$ of the $\bf K$-border are zero.
\par {\bf 2.3.2. Corollary.} {\it Let $M$ be a compact oriented manifold
over $\bf K$ of dimension $\dim_{\bf K}M=k>0$
with an oriented boundary $\partial M$
and let $w$ be a differential $(k-1)$ $C^{(1,n-1)}$-form as in \S 2.2.5
such that its pull back $\phi ^*w=\sum_{j_1<...<j_{k-1}}
f_{j_1,...,j_{k-1}}dz_{j_1}\wedge ... \wedge dz_{j_{k-1}}$
has $f_{j_1,....,j_{k-1}}$ in $\mbox{ }_{P,z_j}C^n(U,{\bf L})$
by the variable $z_j$ for each $j$ such that $j\in \{ 1,...,k \}
\setminus \{ j_1,...,j_{k-1} \} $, then \\
$(1)\quad \mbox{ }_MP^ndw=\mbox{ }_{\partial M}P^nw$.}
\par {\bf Proof.} Repeating the proof of Theorem 3.1 for
each term $f_{j_1,...,j_{k-1}}dz_{j_1}\wedge ... \wedge dz_{j_{k-1}}$
of $w$ and applying Formulas $3.1.(i,ii)$
we get the statement of this corollary.
\par {\bf 2.4.1. Remarks and Notations.} Let $f\in C^1({\bf K}(\alpha ),
Y)$, where $Y$ is a Banach space over $\bf L$, $\bf L$ is a field
containing ${\bf K}(\alpha )$ such that $\bf L$ is complete relative to its
uniformity, the multiplicative norm in $\bf L$ is the extension of the
multiplicative norm in ${\bf K}(\alpha )$.
As the Banach space ${\bf K}(\alpha )$
over $\bf K$ is isomorphic with ${\bf K}^r$, where $2\le r\in \bf N$.
Consider such structure over $\bf K$.
Then each $\zeta \in {\bf K}(\alpha )$ we write in the form
$\zeta =x+\alpha y$, where $x\in \bf K$, $y\in {\bf K}^{r-1}$.
Denote by ${\bar {\zeta }}:=x-\alpha y$ the so called conjugate element
to $\zeta $. Then $x=(\zeta + {\bar {\zeta }})/2$ and $y=
(\zeta - {\bar {\zeta }})/(2\alpha )$. Therefore,
\par $(i) \quad \partial f(\zeta ,{\bar {\zeta }})/\partial x=
\partial f(\zeta ,{\bar {\zeta }})/\partial \zeta +
\partial f(\zeta ,{\bar {\zeta }})/\partial {\bar {\zeta }}$  and
\par $(ii) \quad \partial f(\zeta ,{\bar {\zeta }})/\partial y=
\alpha \partial f(\zeta ,{\bar {\zeta }})/\partial \zeta -
\alpha \partial f(\zeta ,{\bar {\zeta }})/\partial {\bar {\zeta }},$
consequently,
\par $(iii) \quad \partial f(\zeta ,{\bar {\zeta }})/\partial \zeta =
[\partial f(\zeta ,{\bar {\zeta }})/\partial x +
\alpha ^{-1} \partial f(\zeta ,{\bar {\zeta }})/\partial y]/2$  and
\par $(iv) \quad \partial f(\zeta ,{\bar {\zeta }})/\partial {\bar {\zeta }}=
[\partial f(\zeta ,{\bar {\zeta }})/\partial x -
\alpha ^{-1}\partial f(\zeta ,{\bar {\zeta }})/\partial y]/2$.  \\
In particular, the external differentiation of differential
$C^1$-forms $w$ on a clopen subset $\Omega $ in
$({\bf K}\oplus \alpha {\bf K})^m$ has the form
\par $(v) \quad dw=\partial w+{\bar {\partial }}w$, where
\par $(vi) \quad w=\sum_{I,J}w_{I,J}(\zeta ,{\bar {\zeta }})
d{\zeta }^{\wedge I}\wedge d{\bar {\zeta }}^{\wedge J}$,
\par $(vii) \quad \partial w=\sum_{I,J,l}(\partial w_{I,J}/\partial
\zeta _l)d{\zeta }_l\wedge d{\zeta }^I\wedge d{\bar {\zeta }}^{\wedge J}$,
\par $(viii) \quad {\bar {\partial }} w=
(-1)^{|I|}\sum_{I,J,l}(\partial w_{I,J}/\partial
{\bar {\zeta }}_l)dz^I\wedge d{\bar {\zeta }}^l\wedge
d{\bar {\zeta }}^{\wedge J}$, where
$d\zeta ^{\wedge I}:=d\zeta _{I_1}\wedge ... \wedge d\zeta _{I_b}$,
$d{\bar {\zeta }}^{\wedge J}:=d{\bar {\zeta }}_{J_1}\wedge ... \wedge
d{\bar {\zeta }}_{J_c}$,
$1\le I_1<...<I_b\le m$, $1\le J_1<...<J_c\le m$,
such that $w$ is the $(b,c)$-form
with coefficients $w_{I,J}\in C^1(\Omega ,Y)$, $|I|:=b$. \\
If $r>2$, then the differential $s$-form $w$ can be written as
\par $(ix)\quad w=\sum_{J,|J|=s} w_Jdz^{\wedge J}$,
$z=(z_1,...,z_{rm})$, $z_l\in \bf K$ for each $l=1,...,rm$,
$dz^{\wedge J}:=dz_{J_1}\wedge ... \wedge dz_{J_s}$,
$1\le J_1<...<J_s\le rm$.
Let $\Lambda ({\bf K}(\alpha )^m)$ denote the Grassmann algebra
(exterior algebra) of ${\bf K}(\alpha )^m$, where ${\bf K}(\alpha )$
is considered as a ${\bf K}$-linear space, $\Lambda ({\bf K}(\alpha )^m)=
\bigoplus_{l=0}^{rm}\Lambda ^l({\bf K}(\alpha )^m)$.
Then $w\in C^{\xi }(\Omega ,L(\Lambda ({\bf K}(\alpha )^m),Y))$
is the differential form, since the space $({\bf K}(\alpha )^m)^*$
of ${\bf K}$-linear functionals on ${\bf K}(\alpha )^m$
is the space isomorphic with ${\bf K}(\alpha )^m$ due to discretness
of $\Gamma _{\bf K}$, where $L(\Lambda ({\bf K}(\alpha )^m),Y)$
is the Banach space of ${\bf K}$-linear operators from
$\Lambda ({\bf K}(\alpha )^m)$ into $Y$.
\par Consider $\omega $ such that $\omega \subset E$,
where $E:=\{ z\in {\bf K}(\alpha ): |z|<p^{1/(1-p)} \} $,
since $\exp $ is the bijective analytic function on $E$, therefore, we put
\par $(x)\quad \exp (\omega )=\Omega $, that is, $\omega =log (\Omega )$
for $\Omega \subset 1+E$ (see \S \S 25 and 44
\cite{sch1}). Henceforth, if on a manifold $M$ there will be considered
functions $f$ having the property ${\bar \partial }f=0$,
then it will be supposed that ${\bar \partial }\phi _{i,j}=0$ for
each transition mapping $\phi _{i,j}$, if another will not be specified.
\par Consider an extension of $log$.
Denote by ${\bf C_p}^+:= \{ z\in {\bf C_p}: |z-1|<1 \} $
and ${\bf K}(\alpha )^+:={\bf K}(\alpha )\cap {\bf C_p}^+$.
Then ${\bf K}(\alpha )^+$ is the Abelian subgroup in the additive
group ${\bf C_p}^+$ and ${\bf C_p}^{\times }:= {\bf C_p}\setminus \{ 0 \} $
is the Abelian multiplicative group. The group ${\bf C_p}^{\times }$
is divisible,
that is, for each $y\in {\bf C_p}^{\times }$ and each $n\in \bf N$
there exists $x\in {\bf C_p}^{\times }$ such that $x^n=y$.
Let $X$ be a proper divisible subgroup in ${\bf C_p}^{\times }$ such that
${\bf C_p}^+\subset X$. Let $G$ be a subgroup generated by
$X$ and $y\in {\bf C_p}^{\times }\setminus X$.
Suppose $y^n\notin X$ for each
$n\in \bf N$, then for each $g\in G$ there exist unique
$n\in \bf Z$ and $x\in X$ such that $g=y^nx$. Choose
$z\in \bf C_p$, then put $Log (g):=nz+Log (x)$.
The second possibility is: $y^n\in X$ for some $n\in \bf N$, $n>1$.
For each $g\in G$ there exist unique $n\in \{ 0,1,...,m-1 \} $
and $x\in X$ such that $g=y^nx$, where $m:=\min_{y^n\in X; n\in
{\bf N}}n$. Since $\bf C_p$ is divisible, there exists
$z\in \bf C_p$ such that $z_m=Log (y_m)$, therefore, define
$Log (g):=nz+Log (x)$. Using the Zorn's Lemma we can extend $Log $
from ${\bf C_p}^+$ on ${\bf C_p}^{\times }$. In particular we can consider
values of $Log (i)$ and $Log (\alpha )$ using identities
$Log (1)=0$, $i^4=1$, $\alpha ^m \in \bf K$, $\alpha ^n=1$ for some minimal
$n\in \bf N$. In view of Theorem 45.9 \cite{sch1} we can choose
an infinite family of branches of $Log $ indexed by $\bf Z$.
For the convenience put $Log (0):=A$.
\par From the consideration above it follows, that
the extension $Exp$ of $\exp $ on $\bf C_p$
and the extension $Log$ of $log$ on ${\bf C_p}\setminus \{ 0 \} $
can be chosen such that directed going (defined by going from
$0$ to $\beta $ in linearly ordered $B({\bf K},0,1)$,
see \S 2.2) by the oriented loop
$\partial_cB({\bf K}(\alpha ),0,p^{-2})$
changes a branch $\mbox{ }_n Log$ of $Log$ on $1$ in the following
manner: $\mbox{ }_{n+1} Log (x)-\mbox{ }_n Log (x)=:\delta \ne 0$
for each $n\in \bf N$, where $Exp (\delta )=1$, $\delta $ is independent
of $n$. This is possible, since
algebraically $\bf C_p$ and $\bf C$ are isomorphic fields \cite{kobl},
also points $p^2(-1,-1)$, $p^2(1,-1)$, $p^2(1,1)$ and $p^2(1,-1)$ belong
to $\partial_cB({\bf K}(\alpha ),0,p^{-2})$.
\par {\bf 2.4.2. Theorem.} {\it Let $M$ be a compact
$\mbox{ }_SC^{(q,n)}$-manifold over $\bf K$ satisfying conditions of
\S 2.2.5 and \S 2.4.1 
for which $\phi ^{-1}(M)=\Omega \subset {\bf K}(\alpha )$ with
a $\bf K$-boundary $\gamma :=\partial M$,
$dim_{\bf K}M=2$, $2\le r\in \bf N$, $0\le q\in \bf Z$, $1\le n\in \bf N$,
then there exists a constant
$0\ne C:=C_n(\alpha )\in {\bf K}(\alpha )$, such that
\par $(1)\quad f(z)=C^{-1}\mbox{ }_{\partial M}P^n \{ f(\zeta)
(\zeta -z)^{-1} d\zeta \}
- C^{-1}\mbox{ }_MP^n \{ ({\overline {\partial }}f\wedge d\zeta )/
(\zeta -z) \} $   \\
for each
\par $(2)$ $f_1(z+Exp (\eta ))=:\psi (\eta )
\in \mbox{ }_SC^{(1,n-1)}(\omega _{\epsilon },Y)$
and each marked $z\in M$ encompassed by $\gamma $, \\
each $\epsilon =\epsilon _j$,
$0<\epsilon _j$ for  each $j\in \bf N$, $\{ \epsilon _j : j \} $
is a sequence in $\Gamma _{\bf K}$ with $\lim_{j\to \infty }
\epsilon _j=0$, $f_1:=f\circ \phi $,
where $\omega := \omega (z) :=
\{ \eta \in {\bf K}(\alpha ):$ $z+Exp (\eta )\in \Omega \} $,
$\omega _{\epsilon }:=\omega \setminus
Log (B({\bf K}(\alpha ),z,\epsilon ))$, $z\in \Omega $.
Moreover, $C_n(\alpha )=C_1(\alpha )=\delta $ for each $n\in \bf N$. }
\par {\bf Proof.} Using the $\mbox{ }_SC^{(q,n)}$-diffeomorphism
$\phi $ reduce the proof to the case of $f$ on $\Omega $.
Consider the differential form $w:=f(\zeta )(\zeta -z)^{-1}d\zeta $
on $\Omega \setminus \{ z \} $, then $dw= -(\zeta -z)^{-1}
(\partial f/\partial {\bar {\zeta }}) d\zeta \wedge d{\bar {\zeta }}$.
Let $s\in \bf Z$ be such that
$\inf_{\zeta \in \partial \Omega }|\zeta -z|=|\pi |^s$.
Take the change of variables $\zeta =z+Exp (\eta )$,
hence $(\zeta -z)^{-1}d\zeta =d\eta $; also
take $l>s$, then from Corollary 2.3.2 and \S \S 2.2.5, 2.4.1 it follows
\par $(i)\quad \mbox{ }_{\Omega \setminus B({\bf K}(\alpha ),z,|\pi |^l)}
P^ndw=\mbox{ }_{\partial \Omega }P^nw- \mbox{ }_{\partial_cB({\bf K}(\alpha ),
z,|\pi |^l)}P^nw$, \\
since $f_1(z+Exp (\eta ))=:\psi (\eta )
\in \mbox{ }_SC^{(1,n-1)}(\omega _{\epsilon },Y)$
and from $\psi \in \mbox{ }_SC^{(1,n-1)}(\omega _{\epsilon },Y)$
it follows $\psi (x,y)\in \mbox{ }_{P,x}C^{(0,n)}(\omega _{\epsilon ,x},Y)$
and $\psi (x,y)\in \mbox{ }_{P,y}C^{(0,n)}(\omega _{\epsilon ,y},Y)$
for each $\epsilon =\epsilon _j$ and for each $x, y$, where
$z=(x,y)$, $\omega _{\epsilon ,x}=\pi _x(\omega _{\epsilon })$,
$\omega _{\epsilon ,y}=\pi _y(\omega _{\epsilon })$,
$\pi _x: {\bf K}\oplus \alpha {\bf K}\to \bf K$ and
$\pi _y: {\bf K}\oplus \alpha {\bf K}\to \alpha \bf K$
are projections, $Y$ is a Banach space over $\bf L$ such that
${\bf K}(\alpha )\subset \bf L$.
The differential form $w$ can be written as $w=f(\zeta )d Log (\zeta -z)$.
From $Log (xz)=Log (x)+Log (z)$ for each $x, z\in {\bf C_p}^{\times }$
it follows, that directed going by the oriented loop
$\partial_cB({\bf K}(\alpha ),0,|\pi |^l)$
changes a branch $\mbox{ }_n Log$ of $Log$ on $1$ in the following
manner: $\mbox{ }_{n+1} Log (x)-\mbox{ }_n Log (x)=:\delta \ne 0$
for each $s\in \bf Z$.
In view of \S 2.4.1 there exists
\par $\lim_{l\to \infty }
\mbox{ }_{\partial_cB({\bf K}(\alpha ), z,|\pi |^l)}P^nw=:
C_n(\alpha )f(z)$.
Finally $\mbox{ }_{\Omega }P^n((\zeta -z)^{-1} {\bar {\partial }}
f(\zeta )\wedge d\zeta )=-
\mbox{ }_{\Omega }P^n((\zeta -z)^{-1} ({\partial }f(\zeta )/\partial
{\bar {\zeta }})d\zeta \wedge d{\bar {\zeta }} )$, where for short
we write $f=f(\zeta )=f(\zeta ,{\bar {\zeta }})$.
\par In view of Formulas $2.1.(2,3)$ and the non-Archimedean Taylor
formula for $C^n$-functions (see Theorem 29.4 \cite{sch1})
\par  $\mbox{ }_{\partial_cB({\bf K}(\alpha ),z,|\pi |^l)}
P^n[(\zeta -z)^{-1}d\zeta ]=
\mbox{ }_{\partial_cB({\bf K}(\alpha ),z,|\pi |^l)}
P^1[d Log (\zeta -z)] +\epsilon (\pi ^l)$  \\
such that there exists a constant $0<b<\infty $ for which
$| \epsilon (\pi ^l)|\le b | \pi |^l$ for each $l\in \bf N$.
On the other hand, due to the Taylor formula for $C^1$-functions
and Formulas $2.1.(2,3)$:
\par  $\mbox{ }_{\partial_cB({\bf K}(\alpha ),z,|\pi |^l)}
P^1 [d Log (\zeta -z) ]= \delta +\eta (\pi ^l)$, \\
where $\lim_{l\to \infty } \eta (\pi ^l)=0.$ Therefore,
$C_n(\alpha )=C_1(\alpha )=\delta \ne 0$.
\par {\bf 2.4.3. Corollary.} {\it Let suppositions of Theorem 2.4.2
be satisfied for each $z\in M$ encompassed by $\partial M$,
then $\partial f(z)/\partial {\bar z}=0$
for each $z\in M$ encompassed by $\partial M$ if and only if
\par $(1)\quad f(z)=C^{-1}\mbox{ }_{\partial M}P^n \{ f(\zeta)(\zeta -z)^{-1}
d\zeta \} $ \\
for each $z\in M$ encompassed by $\partial M$.}
\par {\bf Proof.}  If $\partial f/\partial {\bar {\zeta }}=0$ on $M$,
then the second term in Formula $2.4.2.(1)$ is equal to zero,
that gives Formula $2.4.3.(1)$.
Vice versa, let Formula $2.4.3.(1)$ be satisfied for each $z\in M$
encompassed by $\partial M$. Since $(\partial z /\partial
{\bar z})=0$, $\partial {\bar z}/\partial z =0$,
then $\partial ((\zeta -z)^{-1})/\partial {\bar z}=0$,
consequently, $\partial f(z)/\partial {\bar z}=0$.
\par {\bf 2.4.4. Corollary.} {\it Let suppositions of Theorem 2.4.2
be satisfied for each $z\in M$ encompassed by $\partial M$
and $\partial f(z)/\partial {\bar z}=0$
for each $z\in M$ encompassed by $\partial M$, then
$f$ is locally $z$-analytic in a neighbourhood of each point
$\zeta $ in $M$ encompassed by $\partial M$.}
\par {\bf Proof.} Using the mapping $\phi $ we can consider $\Omega $
instead of $M$. Let $z\in \Omega $ and $B({\bf K}\oplus \alpha {\bf K},z,
R)\subset \Omega $ such that  $0<R<\inf \{ |z-y|: y\in \partial \Omega \} .$
Consider $x\in B({\bf K}\oplus \alpha {\bf K},z,R/p)$, then \\
$(\zeta -x)^{-1}=(\zeta -z +z-x)^{-1}=(\zeta -z)^{-1}\sum_{l=0}^{\infty }
(x-z)^l/(\zeta -z)^l \in {\bf K}(\alpha )$, \\
where $\zeta \in B({\bf K}\oplus \alpha {\bf K},z,R)$.
Applying Formula $2.4.3.(1)$ we get:
\par $(i)\quad f(x)=C^{-1}\mbox{ }_{\partial_cB}
P^n((\zeta -x)^{-1}f(\zeta )d\zeta )=C^{-1}\sum_{l=0}^{\infty }(x-z)^l
\mbox{ }_{\partial_cB}P^n[(\zeta -z)^{-l-1}f(\zeta )d\zeta ]$ \\
for each $x\in B({\bf K}\oplus \alpha {\bf K},z,R/p)$, since
\par $(ii)\quad |\mbox{ }_{\partial_cB}P^n[(\zeta -z)^{-l-1}f(\zeta )
d\zeta ]| \le \| f \|_{C^{n-1}(\partial _cB,{\bf K}(\alpha ))}
\max_{j,s =0,...,n-1}(R^{j-l-s}/|(j+1)!|) $ \\
and the series is uniformly converging on
$B({\bf K}\oplus \alpha {\bf K},z,R/p)$,
where $B=B({\bf K}\oplus \alpha {\bf K},z,R),$
${\bf K}\oplus \alpha {\bf K}\subset {\bf K}(\alpha )$,
hence $f(x)$ is locally $x$-analytic.
\par {\bf 2.4.5. Definition.} Let $\Omega $ be as in \S 2.2.3.
Two paths $\gamma _0: B({\bf K},0,1)\to \Omega $ and
$\gamma _1: B({\bf K},0,1)\to \Omega $ with common ends
$\gamma _0(0)=\gamma _1(0)=a$, $\gamma _0(\beta )=\gamma _1(\beta )=b$
are called affine homotopic in $\Omega $, if there exists
a continuous mapping $\gamma (x,y): B({\bf K},0,1)^2\to \Omega $
such that
\par $(i)$ $\gamma (0,y)=\gamma _0(y)$, $\gamma (\beta ,y)=\gamma _1(y)$
for each $y\in B({\bf K},0,1)$,
\par $(ii)$ $\gamma (x,0)=a$, $\gamma (x,\beta )=b$ for each
$x\in B({\bf K},0,1)$,
\par $(iii)$ there exists a sequence $\{ \gamma _n(x,y): n\in {\bf N} \} $
of continuous mappings, $\gamma _n: B({\bf K},0,1)^2\to \Omega $
such that each $\gamma _n$ is locally affine and $\{ \gamma _n : n \} $
converges uniformly to $\gamma $ on $B({\bf K},0,1)^2$, where
$\gamma _n(x,y)=(1-x/\beta )\gamma _n(0,y)+x\gamma _n(\beta ,y)/\beta $
for each $x\in B({\bf K},0,1),$ $\gamma _n(0,y)$ and $\gamma _n(\beta ,y)$
are locally affine (see \S 2.2.4). In particular, for $a=b$ this
produces the definition of affine homotopic loops. We call $\Omega $
(or $M$) affine homotopic to a point, if $\partial \Omega $
(or $\partial M$ respectively) is affine homotopic to a point
$z$ in $\Omega $
(or $z$ in $M$ correspondingly, see \S 2.2.5).
\par {\bf 2.4.6. Theorem.} { \it Let conditions of Theorem 2.4.2
be satisfied for each $z\in M$ and let $M$ be affine homotopic
to a point, where $\partial f(z,{\bar z})/\partial {\bar z}=0$
for each $z\in M$ encompassed by $\partial M$. Then
\par $(1)$ $\mbox{ }_{\gamma _0}P^n[fd\zeta ]=
\mbox{ }_{\gamma _1}P^n[fd\zeta ]$ \\
for each two paths $\gamma _0$ and $\gamma _1$ which are affine
homotopic in $M$.}
\par {\bf Proof.} Using the diffeomorphism $\phi $ we can consider
$\Omega $ instead of $M$. For each $\epsilon >0$ there exists
a finite partition of a suitable subset $\Omega _{\epsilon }$
into finite union of parallelepipeds of diameter less, than
$\epsilon $ in the proof of Theorem 2.3.1, where $\Omega _{\epsilon }
\subset \{ z\in \Omega : d(z,\partial \Omega )< \epsilon \} $,
$cl (\bigcup_{\epsilon >0} \Omega _{\epsilon })=\Omega $.
In view of Corollary 2.4.3 $0=f(z)(z-\zeta )|_{z=\zeta }=
C(\alpha )^{-1}\mbox{ }_{\partial \xi }P^n[f(\zeta )d\zeta ]$
for each such parallelepiped $\xi $. Therefore, there exists
a sequence $\{ \gamma _l: l \} $ of affine homotopy such that
$\gamma _l(0,y)$ and $\gamma _l(\beta ,y)$ are contained in the
union $\bigcup_{\xi \subset \Omega }\partial \xi $ for each
$\epsilon _l=|\pi |^l$, $l\in \bf N$. Since
$\mbox{ }_{\gamma _l(0,*)}P^n[fd\zeta ]=
\mbox{ }_{\gamma _l(\beta ,*)}P^n[fd\zeta ]$
for each $l$ and taking $l$ tending to the infinity
we get $(1)$ due to continuity of the operator $P^n$.
\par {\bf 2.4.7. Corollary.} {\it Let $f$ satisfies conditions
of Corollary 2.4.4 with $\Omega =B({\bf K}\oplus \alpha {\bf K},z,R)$.
Then
\par $(i)$ $|f(x)|\le |C|^{-1} \max_{j, s=0,...,n-1}
(\| f^{(j-s)} \|_{C^0(\partial _cB,Y)} s!{j\choose s}
R^{j-s-1}/(j+1)!$\\
$\le |C|^{-1} \| f \|_{C^{(n-1)}(\partial _cB,Y)}
\max_{j, s =0,...,n-1} (R^{j-s-1}/(j+1)!)$.}
\par {\bf Proof.} From \\
$\partial ^j_{\zeta }f(\zeta )(\zeta -x)^{-1}
=\sum_{s=0}^j s!{j\choose s}(-1)^sf^{(j-s)}(\zeta )(\zeta -x)^{-1-s}$ \\
and $|\zeta _{l+1}-\zeta _l|\le  R$ on $\partial _cB$ and \S 2.1
it follows Inequality $(i)$.
\par {\bf 2.4.8. Remark.} The field $\bf K$ is locally compact,
then ${\bf T}_q$ is not contained in $\bf K$, where ${\bf T}_q$
is a group of all $q^n$-roots $b$ of the unity: $b^l=1$,
$l=q^n$, $n\in \bf N$, $q$ is the prime number, since
$dim_{\bf Q_p}{\bf Q_pT}_q=\infty $ for ${\bf Q_p}\subset
\bf K$ and $\bf K$ would be nonlocally compact whenever
${\bf T}_q\subset \bf K$, which is impossible by the supposition
on $\bf K$. Therefore, there exists $\min \{ s\in {\bf N}:
b^q\in {\bf K}, b\notin {\bf K},$ where $b\ne 1$ is the
$q^{s+1}$-root of the unity $ \} $. Hence there exists
$\zeta \in \bf K$ such that $\zeta ^{1/q}\notin \bf K$.
In particular, it is true for $q=2$. Therefore, each local
field $\bf K$ has a quadratic extension ${\bf K}(\alpha )$ such that
$\alpha \notin \bf K$. In the particular case ${\bf K}=\bf Q_p$
there exists the finite field ${\bf F_p}:=R/P$ (see \S 2.2.1).
Then ${\bf F_p}\setminus \{ 0 \} $ is the multiplicative
group consisting of $p-1$ elements. If $p=4n+1$, where $1\le n\in \bf N$,
then $\bf Q_p$ contains $i=(-1)^{1/2}$.
\par {\bf 2.5.1. Lemma.} {\it If $f$ is locally $z$-analytic on
$M$, where $M$ is a locally compact $C^{(0,n)}$-manifold satisfying
conditions of \S \S 2.2.5 and 2.4.1, $\phi ^{-1}(M)=\Omega \subset
{\bf K}(\alpha )$, $dim_{\bf K}M=2$, $2\le r\in \bf N$,
then $\partial f(z,{\bar z})/\partial {\bar z}=0$ on $M$.}
\par {\bf Proof.} Using the diffeomorphism $\phi $ we can
consider $\Omega $ instead of $M$. Since for each $z\in \Omega $
there exists $0<R<\infty $ such that $B:=B({\bf K}\oplus \alpha
{\bf K},z,R)\subset \Omega $ and $f(z,{\bar z})=
\sum_{k=0}^{\infty } (\zeta -z)^kf_k$ on $B$, where $f_k\in Y$,
then there exist $\partial f/\partial z$ and $\partial f/\partial
{\bar z}=0$ on $B$. Since $z\in \Omega $ is arbitrary and such balls
form the covering of $\Omega $, then
$\partial f/\partial {\bar z}=0$ on $\Omega $.
\par {\bf 2.5.2. Remark.} Let $n\ge 1$, then $(d/dz)\mbox{ }_{\Omega }
P^n=I: C^{n-1}(\Omega ,{\bf L})\to C^{n-1}(\Omega ,{\bf L})$.
But $P^n d/dz\ne I$ on $C^n(\Omega ,{\bf L})$, where
$P^n d/dz: C^n(\Omega ,{\bf L})\to C^n(\Omega ,{\bf L}).$
If $\mbox{ }_PC^n(\Omega ,{\bf L})$ would be dense in
$C^n(\Omega ,{\bf L})$, then $P^nd/dz$ would have the continuous
extension $I$ on $C^n(\Omega ,{\bf L})$, since $P^nd/dz$ is the continuous
operator from $C^n$ into $C^n$ and $P^n(d/dz)|_{\mbox{ }_PC^n_0}=I$.
Therefore, $\mbox{ }_PC^n(\Omega ,{\bf L})$ is not dense in
$C^n(\Omega ,{\bf L})$. On the other hand,
$C^1(\Omega ,{\bf L})=
\mbox{ }_PC^1_0(\Omega ,{\bf L})\oplus N^1$,
where $N^1:= \{ f\in C^1: f'=0 \} $ is the closed $\bf L$-linear
subspace in $C^1$ (see Theorem 5.1 and Corollary 5.5 \cite{sch2}).
\par {\bf 2.5.3. Theorem.} {\it Let $f$ be a function on $M$
over $\bf K$ satisfying
Conditons $2.2.5$ and $2.4.1$ and $\gamma $ be a loop in $M$ satisfying
Conditions $2.2.4$ and let $\gamma $ be affine homotopic to a point
in $M$, $dim_{\bf K}M=2$, $\phi ^{-1}(M)=\Omega \subset {\bf K}(\alpha )$,
$2\le r\in \bf N$, $f$ satisfies Condition $2.4.2.(2)$ for each $z\in M$
and $\partial f(z,{\bar z})/\partial {\bar z}=0$ on $M$, then
$\mbox{ }_{\gamma }P^nf=0$.}
\par {\bf Proof.} Let $V$ be a submanifold in $M$ such that
$\partial V=\gamma $. In view of Theorem $2.4.6$
$\mbox{ }_{\gamma }P^n[fd\zeta ]=\mbox{ }_{\gamma _{\epsilon }}P^n
[fd\zeta ]$, where $\gamma $ and $\gamma _{\epsilon }$ are affine
homotopic and $0< diam (\gamma _{\epsilon })<\epsilon $. In view
of continuity of the operator $P^n$ there exists $\lim_{\epsilon \to 0}
\mbox{ }_{\gamma _{\epsilon }}P^n[fd\zeta ]=0$.
\par {\bf 2.5.4. Theorem.} {\it If $f$ satisfies Condition $2.4.2.(2)$,
a manifold $M$ over $\bf K$ satisfies Conditions $2.2.5$ and $2.4.1$
and $M$ is affine homotopic to a point, $dim_{\bf K}M=2$,
$2\le r\in \bf N$ and $\mbox{ }_{\gamma }P^nf=0$
for each loop $\gamma $ in $M$ satisfying Conditions $2.2.4$, then
$\partial f(z,{\bar z})/\partial {\bar z}=0$ for each $z\in M$
encompassed by $\partial M$.}
\par {\bf Proof.} Using the diffeomorphism $\phi $ we can consider
$\Omega $ instead of $M$. Choose a marked point $z_0$ in $M$.
Let $\eta $ be a path joining points $z_0$ and $z$ and satisfying
Conditions $2.2.4$. From $\mbox{ }_{\gamma }P^nf=0$ it follows,
that $\mbox{ }_{\eta }P^nf$ does not depend on $\eta $ besides points
$z_0=\eta (0)$ and $z=\eta (\beta )$,
since each two points in $\Omega $ can be joined by an affine path,
hence it is possible to put
$F(z):=\mbox{ }_{\eta ; \eta (0)=z_0; \eta (\beta )=z}P^nf$
such that $F$ is a function on $\Omega $. In view of Formulas $2.4.1.(i-iv)$, 
\par $(i)\quad \partial F(z)/\partial z=f(z)$. \\
In view of theorem $2.4.2$
\par $(ii)\quad 0=\mbox{ }_{\gamma }P^n(f(\zeta )d\zeta )=-C^{-1}
\mbox{ }_UP^n((\partial f(\zeta ,{\bar {\zeta }})/\partial {\bar {\zeta }})
d\zeta \wedge d{\bar {\zeta }} ) $ \\
for each submanifold $V$ in $M$ with the loop $\gamma =\partial V$,
$dim_{\bf K}V=2$.
Since $V$ is arbitrary, then ${\bar {\partial }}f(z,{\bar z})=0$
at each point $z\in M$ encompassed by $\partial M$.
\par {\bf 2.5.5. Corollary.} {\it Let conditions of Theorem 2.5.4
be satisfied, then $f$ has an antiderivative $F$ such that
$F'=f$ on $M$.}
\par {\bf 2.6.1. Lemma.} {\it Let $\Omega $ be a clopen compact subset
in ${\bf K}^m$, then for each $y \in \Omega $ there exists
a ball $B$ such that $y\in B\subset \Omega $ and
$\mbox{ }_PC^{\xi }(\Omega ,Y)|_B=\mbox{ }_PC^{\xi }(B,Y)$
and $\mbox{ }_SC^{\xi }(\Omega ,Y)|_B=\mbox{ }_SC^{\xi }(B,Y)$
for each $\xi $, where $\mbox{ }_PC^{\xi }(\Omega ,Y)|_B:=
\{ g|_B: g\in \mbox{ }_PC^{\xi }(\Omega ,Y) \} $
and $\mbox{ }_SC^{\xi }(\Omega ,Y)|_B:=
\{ g|_B: g\in \mbox{ }_SC^{\xi }(\Omega ,Y) \} $,
$Y$ is a Banach space over $\bf L$,
$\xi =(t,n)$, $0\le n\in \bf Z$, $0\le t\in \bf Z$,
$1\le n$ for $\mbox{ }_PC^{\xi }$, $1\le t$
for $\mbox{ }_SC^{\xi }$.}
\par {\bf Proof.} Let $\sigma $ be an approximation of the unity in
$U$. In view of \S 2.1 it is sufficient to consider the case
$m=1$. Choose $R=\rho ^{s+1/2}$ for sufficiently large $s\in \bf N$
such that $0\in B:=B({\bf K},0,R)$. If $x\in B({\bf K},y,R)$, then
$\sigma _s(x)=\sigma _s(y)$ due to $2.1.(iii)$. From Formula $2.1.(ii)$
it follows, that $\sigma _l(x)=\sigma _l(y)$ for each $l<s$.
Moreover, $\sigma _l(x)=:x_l\in B$ for each $l\ge s$, since
$\rho ^{s+1}< R < \rho ^s$ and the valuation group $\Gamma _{\bf K}
:= \{ |q|_{\bf K}: 0\ne q \in {\bf K} \} $ of $\bf K$ is discrete,
since $\bf K$ is locally compact. Therefore,
\par $(i)\quad [\mbox{ }_UP^nf(x) ] - [\mbox{ }_BP^nf(x)]=
\sum_{j=0}^{n-1}\sum_{k=0}^{s-1} f^{(j)}(x_k)(x_{k+1}-x_k)^{j+1}/[(j+1)!]$ \\
for each $f\in C^{\xi }(U,Y)$,
where $x_k=y_k$ is fixed, and the term on the right-hand side
of $(i)$ is independent of $x\in B$, that is, constant on $B$.
Hence $g\in \mbox{ }_PC^{\xi }(U,Y)$
if and only if $g|_B\in \mbox{ }_PC^{\xi }(B,Y)$.
From $2.1.(3)$ and $\chi _{\Omega } \chi _B= \chi _B= \chi _{\Omega }|_B$
the statement of this lemma follows.
\par {\bf 2.6.2. Definition.} Let a manifold $M$ be satifying
Conditions $2.4.2$, $f\in C^{(q,n-1)} (M,Y)$, $0\le q\in \bf Z$,
$1\le n\in \bf N$, $Y$ ia a Banach space over $\bf L$, ${\bf K}(\alpha )
\subset \bf L$.
Then put in the sence of distributions:
\par $(i)\quad \mbox{ }_MP^n(gf'):=-\mbox{ }_MP^n(g'f)$ \\
for each $g\in \mbox{ }_SC^{(1,n-1)}(M,Y^*)$
with $supp (g)\subset {\tilde M}:=\{ z\in M:$ $z\mbox{ is encompassed by }
\partial M \} $,
where $M\hookrightarrow {\bf K}(\alpha )^N$ (see Theorem 2.2.6),
$Y^*$ is the topologically dual space of all $\bf L$-linear
continuous functionals $\theta : Y\to \bf L$, the valuation
group $\Gamma _{\bf L}$ of $\bf L$ is discrete.
\par {\bf 2.6.3. Theorem.} {\it Let a manifold $M$ satisfy Conditions
$2.4.2$ and let $f$ satisfy $2.4.2.(2)$ for each $z\in M$, then the function
\par $(1) \quad u(z) :=
{\bar z}C_n(\alpha )^{-1}\mbox{ }_{\partial M}P^n[f(\zeta )
(\zeta -z)^{-1}d\zeta ] - C_n(\alpha )^{-1} \mbox{ }_MP^n[f(\zeta )
(\zeta -z)^{-1} d{\bar {\zeta }} \wedge d\zeta ]$  \\
is a solution of the equation
\par $(2)\quad \partial u (z)/\partial {\bar z}=f(z)$ \\
in the sence of distributions 
for each $z\in M$ encompassed by $\partial M$.}
\par {\bf Proof.} The space $\mbox{ }_PC^{(q,n)}(M,Y)$
is dense in $C^{(q,n-1)}(M,Y)$. Indeed, for each $\delta >0$ and for each
continuous function $f\circ \phi $ on $\Omega $
or a continuous partial difference quotient
$w_q:={\bar \Phi }^qf\circ \phi (x;h_1^{\otimes s_1},...,h_m^{\otimes s_m};
\zeta _1,...,\zeta _q)$ on a domain contained
$\Omega ^{q+1}\times B({\bf K},0,1)^t$
with $0\le t\le (n-1)m$, $0\le s_j\le n$ for each
$j=1,...,m$, $t=s_1+...+s_m$, $x, x+\zeta _jh_j\in \Omega $, $h_j\in V$,
$\zeta _j\in B({\bf K},0,1)$, $V$ is a neighbourhood of $0$ in
$\bf K^m$, $\Omega +V\subset \Omega $ (see  \cite{luambp00}),
$m:=dim_{\bf K}M$, there exists a finite partition
of $\Omega ^{q+1}$ into disjoint union of balls $B_j$ such that on each
$B_j$ the variation $var (w_q):= \sup _{x, y\in B_j}
|w_q(x)-w_q(y)| < \delta $, since $M$ is compact and for each covering
of $M$ by such balls there exists a finite subcovering.
Therefore, in $C^{(q,n-1)}(\Omega ,Y)$ the subspace
$\Sigma ^{(q,n-1)}(\Omega ,Y)$ of all
$C^{(q,n-1)}(\Omega ,Y)$-functions $f$ such that
$w_{(n-1)m}$ corresponding to $f$ is locally constant on the diagonal
$\Delta \Omega ^{(n-1)m+1}:=
\{ (y_1,...,y_{(n-1)m+1})\in \Omega ^{(n-1)m+1}: y_1=...=y_{(n-1)m+1} \} $
is dense. Since the operator $\mbox{ }_{\Omega }P^n$ is continuous, then
$\mbox{ }_{\Omega }P^n(\Sigma ^{(q,n-1)}(\Omega ,Y))$
is dense in $C^{(q,n-1)}_0(\Omega ,Y)$
and $\mbox{ }_S\Sigma ^{(q+1,n-1)}(M,Y):= \{ f\in \mbox{ }_SC^{(q+1,n-1)}
(M,Y):$ $f=y_l+\mbox{ }_MP^n_{x_l}g_l$ $\forall l=1,...m, y_l\in Y,$
$g_l\in \Sigma ^{(q,n-1)}(M,Y) \} $ is dense in $\mbox{ }_SC^{(q+1,n-1)}
(M,Y).$ From $\mbox{ }_S\Sigma ^{(q+1,n-1)}(M,Y)\subset
\mbox{ }_SC^{(q+1,n-1)}(M,Y)\subset \mbox{ }_PC^{(q,n)}(M,Y)
\subset C^{(q,n-1)}(M,Y)$ it follows,
that $\mbox{ }_SC^{(q+1,n-1)}(M,Y)$ is dense in $C^{(q,n-1)}(M,Y)$.
In particular, take $\bf L$ such that ${\bf K}(\alpha )\subset \bf L$.
Since $\mbox{ }_PC^{(0,n)}(M,Y^*)\subset
\{ g': g\in \mbox{ }_SC^{(1,n-1)}(M,Y^*) \} \subset 
C^{(0,n-1)}(M,Y^*)$, then the family of functionals $\{ \mbox{ }_M
P^n(g'f): g\in \mbox{ }_SC^{(1,n-1)}(M,Y^*) \} $ separates points
of $C^{(q,n-1)}(M,Y)$, since $Y^*$ separates points of $Y$
for discrete $\Gamma _{\bf L}$ (see Theorem 4.15 in \cite{roo}).
In view of Formula $2.6.2(i)$ it is sufficient to prove
this theorem for $f\circ \phi (z+Exp (\eta )) =: \psi (\eta )
\in \mbox{ }_PC^{(0,n)}(\omega _{\epsilon ,x},Y)\cap
\mbox{ }_PC^{(0,n)}(\omega _{\epsilon ,y},Y)$ for each
$\epsilon =\epsilon _j$, where $\omega (z):=
\omega := \{ \eta \in {\bf K}(\alpha ):
z+Exp (\eta )\in \Omega \}$, $z\in \Omega $, $\omega _{\epsilon }
=\omega \setminus Log (B({\bf K}(\alpha ),z,\epsilon ))$,
$\epsilon =\epsilon _j$, $\omega _{\epsilon ,x}=
\pi _x(\omega _{\epsilon })$, $\omega _{\epsilon ,y}:
=\pi _y(\omega _{\epsilon })$. 
\par Using the diffeomorphism $\phi $ we can consider
$\Omega $ instead of $M$. Choose a clopen ball $B:=B({\bf K}\oplus
{\bf K}(\alpha ),z_0,R)\subset \Omega $ containing a point
$z_0 \in \Omega $ and its characteristic function $\chi :=\chi _B$. Then
$(f\chi )_1\in \mbox{ }_PC^{(0,n)}(b_x,Y)\cap
\mbox{ }_PC^{(0,n)}(b_y,Y)$ for suitable
$0< R < \infty $, where $b:= \{ \eta \in {\bf K}(\alpha ): z_0+Exp(\eta )
\in B \} $, $b_x:=\pi _x(b)$, $b_y:=\pi _y(b)$ (see Lemma $2.6.1$).
Using the affine mapping $z\mapsto (z-z_0)$
we can consider $0$ instead of $z_0$.
Then $B$ is the additive group.
We can take $R>0$ sufficiently small such that each point of $B$
is encompassed by $\partial \Omega $. Therefore,
\par $(3)\quad u=u_1+u_2$, where
\par $(4) \quad u_1(z) :=C^{-1}{\bar z}
\mbox{ }_{\partial \Omega }P^n [f(\zeta )(\zeta -z)^{-1}d\zeta ]
-C^{-1} \mbox{ }_{\Omega } P^n[\chi (\zeta )f(\zeta )
(\zeta -z)^{-1} d{\bar {\zeta }} \wedge d\zeta ] $,
\par $(5) \quad u_2(z) := -C^{-1} \mbox{ }_{\Omega }
P^n[ (1 - \chi (\zeta )) f(\zeta )
(\zeta -z)^{-1} d{\bar {\zeta }} \wedge d\zeta ]$. \\
From $(5)$ it follows, that $\partial u_2(z) / \partial {\bar z} =0$
on $B$. From $(4)$ it follows \\
$u_1(z)=C^{-1}{\bar z}
\mbox{ }_{\partial \Omega }P^n [f(\zeta )(\zeta -z)^{-1}d\zeta ] $ \\
$-C^{-1} \mbox{ }_{\Omega } P^n[\chi (\zeta +z)f(\zeta +z)
\zeta ^{-1} d({\bar {\zeta }}+{\bar z}) \wedge d(\zeta +z)]$  \\
for each $z\in B$, since $B+B=B\subset \Omega $.
Since $\partial _cB$ encompasses $z$ and \\
$\mbox{ }_{\partial _cB} P^n[f(\zeta )(\zeta -z)^{-1}d\zeta ]=
\sum_{l=0}^{\infty }(z-y)^l
\mbox{ }_{\partial _cB} P^n[f(\zeta )(\zeta -y)^{-l-1}d\zeta ]$ \\
for each $z\in B$ with $|z-y|<R$
due to Formula $2.4.4.(ii)$ for this antiderivative, then
\par $(6)$ $\partial \{
\mbox{ }_{\partial _cB} P^n[f(\zeta )(\zeta -z)^{-1}d\zeta ] \}
/\partial {\bar z}=0$.
In view of Formulas $2.1.(1-4)$ and \S 2.2.5  \\
$\partial u_1(z)/ \partial {\bar z}=
C^{-1}\mbox{ }_{\partial \Omega }P^n [f(\zeta )(\zeta -z)^{-1}d\zeta ] $ \\
$-C^{-1}\mbox{ }_{\Omega } P^n \{ [{\bar {\partial }}_{\zeta +z}
[\chi (\zeta +z)f(\zeta +z)]
\zeta ^{-1} \wedge d(\zeta +z) \} $, consequently,  \\
$\partial u_1(z)/ \partial {\bar z}=
C^{-1}\mbox{ }_{\partial \Omega }P^n [f(\zeta )(\zeta -z)^{-1}d\zeta ] $ \\
$-C^{-1}\mbox{ }_{\Omega } P^n \{ [{\bar {\partial }}_{\zeta }
(\chi (\zeta )f(\zeta ))\wedge d\zeta ](\zeta -z)^{-1} \} $.  \\
In view of Theorem $2.4.2$ we get the statement of this theorem,
since $B$ and $y$ are arbitrary
forming covering of each point $z\in \Omega $ encompassed by
$\partial \Omega $.
\par {\bf 2.7.1. Definition.} Let $M$ be a manifold
over $\bf K$ satisfying Conditions $2.4.2$. If
$f\in C^{(q,n-1)}(M,Y)$ and for each
loop $\gamma $ in $M$ $\mbox{ }_{\gamma }P^nf=0$,
then we call $f$ $(q,n)$-antiderivationally holomorphic
on $M$, where $Y$ is a Banach space over $\bf L$, ${\bf K}(\alpha )
\subset \bf L$. If $f\in C^{(q,n)}(M,Y)$
and ${\bar {\partial }}f(z)=0$ for each $z\in M$, then
$f$ we call derivationally $(q,n)$-holomorphic.
\par {\bf 2.7.2. Theorem.} {\it  Let $\Omega $ be a clopen compact subset
in $({\bf K}\oplus \alpha {\bf K})^m$.
Consider the following conditions.
\par $(i)$. $f$ satisfies $2.4.2.(2)$ and
${\bar {\partial }}f(z)=0$ for each $z\in \Omega $ with $z_j$
encompassed by $\partial \Omega _j$ for each $j=1,...,m$, where
$\Omega _j=\pi _j(\Omega )$, $\pi _j(\zeta )=\zeta _j$ for each $\zeta =(
\zeta _1,...,\zeta _m)$, $\zeta _j\in {\bf K}\oplus \alpha {\bf K}$.
\par $(ii)$. $f$ is $(0,n)$-antiderivationally holomorphic on  \\
${\tilde {\Omega }}:= \{ z\in \Omega :$ $z_j\mbox{ is encompassed by }
\partial \Omega _j $ $\forall j=1,...,m \} $.
\par $(iii)$. $f\in C^{(0,n-1)}(\Omega ,Y)$ and for each
polydisc $B=B_1\times ... \times B_m \subset \Omega $,
$B_j=B({\bf K}(\alpha ),z_{0,j},R_j)$ for each $j=1,...,m$,
$f(z)$ is given by the antiderivative
\par $(1)$ $f(z)=C(\alpha )^{-m}\mbox{ }_{\partial B_1}P^n
... \mbox{ }_{\partial B_m}P^n[f(\zeta )(\zeta _1-z_1)^{-1}...
(\zeta _m-z_m)^{-1}d\zeta _1\wedge ... \wedge d\zeta _m]$
for each $z\in B$ with $z_j$ encompassed by $\partial B_j$ for each $j$.
\par $(iv)$ $f$ is locally $z$-analytic, that is,
\par $(2)$ $f(z)=\sum_k a_k (z-\zeta )^k$
in some neighbourhood of $\zeta \in {\tilde {\Omega }}$, $a_k\in Y$,
$k=(k_1,...,k_m)$, $0\le k_j\in \bf Z$, $z^k:=z_1^{k_1}...z_m^{k_m}$,
$z=(z_1,...,z_m)$, $z_j\in {\bf K}(\alpha )$.
\par $(v)$ $f\in C^{\infty }(\Omega ,Y)$.
\par $(vi)$ $f\in C^{(0,n-1)}(\Omega ,Y)$ and
for every polydisc $B$ as in $(iii)$  and each multiorder $k$
as in $(iv)$ derivatives are given by
\par $(3)$ $\partial _z^kf(z)=k!C(\alpha )^{-m}\mbox{ }_{\partial B_1}P^n
... \mbox{ }_{\partial B_m}P^n[f(\zeta )(\zeta _1-z_1)^{-k_1-1}...
(\zeta _m-z_m)^{-k_m-1}d\zeta _1\wedge ... \wedge d\zeta _m]$.
\par $(vii)$ the coefficients in Formula $(2)$ are determined by
the equation:
\par $(4)$ $a_k=\partial _z^kf(z)/k!$.
\par $(viii)$ The power series $(2)$ converges uniformly
in each polydics $B\subset {\tilde {\Omega }}$ with sufficiently small
$b:= \max (R_1,...,R_m,1)$.
\par Then from $(i)$ Properties $(ii-viii)$ follow.
Properties $(iii)$ and $(vi)$ are equivalent.
From $(iii)$ Properties $(iv,v,vii,viii)$ follow.
In the subspace \\
$\{ f\in C^{(0,n-1)}(\Omega ,Y):
f(z_1,...,z_{l-1},z_l+Exp (\eta ),z_{l+1},...,z_m)=:\psi _l(\eta ) $ \\
$\in \mbox{ }_SC^{(1,n-1)}(\omega _{l,\epsilon },Y) $ $\mbox{for each}$
$l=1,...,m$ $\mbox{and each } \epsilon =\epsilon _j \} $, \\
where $\omega _l:=\omega _l(z):= \{ \eta \in {\bf K}(\alpha ):$
$(z_1,...,z_{l-1},z_l+Exp (\eta ),z_{l+1},...,z_m)\in \Omega  \} $,
$\omega _{l,\epsilon }:=\omega _l\setminus Log (B({\bf K}(\alpha ),
z_l,\epsilon ))$, $z\in \Omega $,
$Y$ is a Banach space over $\bf L$ such that ${\bf K}(\alpha )
\subset \bf L$, Properties $(i-iv)$ are equivalent.}
\par  {\bf Proof.} From $(i)$ it follows $(iv)$ due to
repeated application of Corollary $2.4.4$. From $(i)$ it follows
$(iii)$ due to repeated application of Corollary $2.4.3$.
Others statements follow from Theorems $2.5.3$, $2.5.4$,
Lemma $2.5.1$ and Formulas $(i,ii)$ in \S 2.4.4,
since from Formula $(3)$ it follows
\par $(5)$ $|\partial_z^kf(z)/k!|\le |C(\alpha )|^{-m}\sup_{\zeta \in \Omega ,
l} |\partial_z^l f(\zeta )| \max_l[b^{|l|-|k|}/(l+{\bar e})!]
<\infty $, \\
where $l=(l_1,...,l_m)$, $0\le l_j\in \bf Z$, $l_j<n$ for each $j=1,...,m$,
$|l|=l_1+...+l_m$, ${\bar e}:=(1,...,1)\in \bf Z^m$,
$b:=\max (R_1,...,R_m)$. The series $(2)$ with
$a_k$ given by $(4)$ converges uniformly in $B$, when
${\overline {\lim }}_k |a_k|^{1/|k|}b<1$. 
\par {\bf 2.7.3. Corollary.} {\it Spaces $C^{la}(\Omega ,{\bf K}(\alpha ))$
of locally analytic functions $f: \Omega \to {\bf K}(\alpha )$ and the space
$C^{(q,n),dh}(\Omega ,{\bf K}(\alpha ))$
of all derivationally $(q,n)$-holomorphic functions are
rings. If $f$ is derivationally $(q,n)$-holomorphic and $f\ne 0$ on
$\Omega $, then $1/f$ is derivationally $(q,n)$-holomorphic on $\Omega $.}
\par {\bf 2.7.4. Corollary.} {\it If $f$ satisfies Condition $2.4.2.(2)$
and there exists $\zeta \in \Omega $ encompassed by $\partial \Omega $
such that $\partial _z^kf(\zeta )=0$ for each $k$,
then there exists a polydisc $B\subset \tilde {\Omega }$
(see \S 2.7.2) such that $f=0$ on $B$.}
\par {\bf Proof.} In view of Theorem $2.7.2$ there exists
a polydics $B$ such that on it Formulas $2.7.2.(2,4)$ are accomplished.
\par {\bf 2.7.5. Remark.} $\mbox{ }_UP^nz^k|_a^b\ne
(b^{k+1}-a^{k+1})/(k+1)$ for each $a\ne b\in U$, where $k>0$.
In view of Corollary 54.2 and Theorem 54.4
\cite{sch1} and \cite{ami} the spaces $C^{la}(\Omega ,{\bf K}(\alpha ))
\cap C^{((q,n),dh)}(\Omega ,{\bf K}(\alpha ))$,
$C^{la}(\Omega ,{\bf K}(\alpha ))
\cap \mbox{ }_PC^{(q,n)}(\Omega ,{\bf K}(\alpha ))$,
$\{ f\in C^{la}(\Omega ,{\bf K}(\alpha )):
$f$ $\mbox{ is }$(q,n)-$$\mbox{antiderivationally}$ \\
$\mbox{holomorphic}$ $ \} $
are infinite dimensional over ${\bf K}(\alpha )$, since the condition
of the local analyticity means that the expansion coefficients
$a(m,f)$ of the function $f$ in the Amice polynomial basis
${\bar Q}_m$ are such that $\lim_{|m|\to \infty } a(m,f)/
P_m({\tilde u}(m))=0$, where $P_m$ are definite polynomials
(see Formulas $2.6.(i-iii)$ \cite{luseam3}).
\par {\bf 2.7.6. Theorem.} {\it Let $\Omega $ and $f$ be as in
$2.7.2.(i)$. If $\zeta $ is zero of $f$ such that $f$ does not
coincide with $0$ on each neighbourhood of $\zeta $, then there
exists $n\in \bf N$ such that
\par $(1)$ $f(z)=(z-\zeta )^ng(z)$,  \\
where $g$ is analytic and $\phi \ne 0$ on some neighbourhood
of $z$.}
\par {\bf Proof.} In view of Theorem $2.7.2$ there exists a neighborhood
$V$ of $\zeta $ such that $f$ has a decomposition into converging
series $2.7.2.(2)$. If $a_k=0$ for each $k$, then $f|_V=0$.
Therefore, there exists a minimal $k$ denoted by $l$ such that \\
$f(z)=\sum_{k=l}^{\infty }a_k(z-\zeta )^k$, \\
put $g(z)=\sum_{k=0}^{\infty }a_{k+l}(z-\zeta )^k$. \\
Since $a_l\ne 0$, then there exists a neighborhood $\zeta \in W\subset V$
such that $g|_W\ne 0$.
\par {\bf 2.7.7. Theorem.} {\it Let $\Omega $ and two functions
$f_1$ and $f_2$ be satisfying $2.7.2.(i)$ such that $f_1(z)=f_2(z)$
for each $z\in E$, where $E\subset \Omega $ and $E$ contains
a limit point $\zeta \in E'$. Then there exists a clopen subset
$W$ in $\Omega $ such that $\zeta \in W$ and $f_1|_W=f_2|_W$.}
\par {\bf Proof.} Put $f:=f_1-f_2$, then $f$ satisfies Condition
$2.7.2.(i)$ and $f(\zeta )=0.$ In view of Theorem $2.7.6$ $f|_W=0$
for some clopen $W$ in $\Omega $, where $\zeta \in W$.
\par {\bf 2.7.8. Theorem.} {\it Let $f$ satisfy $2.4.2.(2)$
and $f$ be derivationally $(0,n)$-holomorphic on $\Omega :=
\{ z\in ({\bf K}\oplus \alpha {\bf K})^m:$ $R_1\le |z-\xi |\le R_2 \} $,
where $0<R_1<R_2<\infty $, $R_1$ and $R_2\in \Gamma _{\bf K}$. Then
\par $(1)$ $f(z)=\sum_ka_k(z-\xi )^k$ \\
for each $z\in \Omega $ with $R_2>|z|>R_1$, where \\
$(2)$ $a_k=C(\alpha )^{-m}
\mbox{ }_{\partial _cB_{R,1}}P^n ... \mbox{ }_{\partial _cB_{R,m}}P^n
[(\zeta _1-\xi _1)^{-k_1-1}...(\zeta _m-\xi _m)^{-k_m-1}f(\zeta )d\zeta _1
\wedge ... \wedge d\zeta _m]$   \\
for each $k\in \bf Z^m$, $R_1<R<R_2$, $B_{R,l}:=
\{ z_l\in {\bf K}\oplus \alpha {\bf K}:$ $|z_l-\xi _l|\le R \} $,
$k=(k_1,...,k_m)$, $k_l\in \bf Z$, $l=1,...,m$.}
\par {\bf Proof.} Let $\pi _l(z)=z_l$ for each $z=(z_1,...,z_m)\in
({\bf K}\oplus \alpha {\bf K})^m$, where $z_l\in {\bf K}\oplus
\alpha {\bf K}$. Then $\pi _l(\Omega )= \{ z_l\in {\bf K}\oplus
\alpha {\bf K}:$ $R_1\le |z_l-\xi _l|\le R_2 \} $. To prove
the theorem consider $f$ by each variable $z_l$. That is, consider
$z=z_l$ and $m=1$. Let $R_3$ and $R_4$ be such that $R_1<R_3<R_4<R_2$
and $z\in W\subset \Omega $, where $W=\{ z\in {\bf K}\oplus \alpha
{\bf K}:$ $R_3\le |z-\xi |\le R_4 \} $. In view of Theorems $2.4.6$
and $2.7.2.(iv,vi)$
\par $(3)$ $f(z)=C(\alpha )^{-1}\mbox{ }_{\partial W}P^n[f(\zeta )
(\zeta -z)^{-1}d\zeta ]$ \\
$=C(\alpha )^{-1}\mbox{ }_{\partial _cB_{R_4}}P^n
[f(\zeta )(\zeta -z)^{-1}d\zeta ] -C(\alpha )^{-1}
\mbox{ }_{\partial _cB_{R_3}}P^n [f(\zeta )(\zeta -z)^{-1}d\zeta ]$,
since $W$ is the union $W_1\cup W_2$, $dim_{\bf K}(W_1\cap W_2)=1$,
where $W_1$ and $W_2$ satisfy $2.4.6$ and $2.7.2$. The part of the path
$\gamma _{1,2}$ in $W_1\cap W_2$ joining $\partial _cB_{R_4}$ with
$\partial _cB_{R_3}$ and forming two paths $\gamma _1$ and $\gamma _2$
affine homotopic to points in $W_1$ and $W_2$, $\gamma _1\subset W_1$,
$\gamma _2\subset W_2$, such that $\gamma _{1,2}$ is being gone twice
in one and opposite directions. This gives $(3)$.
For each $\zeta \in \partial_cB_{R_4}$ we have
$|(z-\xi ) (\zeta -\xi )^{-1}| <1$, hence $(\zeta -z)^{-1}=
\sum_{k=0}^{\infty }(z-\xi )^k(\zeta -\xi )^{-k-1}$ and inevitably \\
$\mbox{ }_{\partial _cB_{R_4}}P^n[f(\zeta )(\zeta -z)^{-1}d\zeta ]=
\sum_{k=0}^{\infty }a_k(z-\xi )^k$, \\
where $a_k=C(\alpha )^{-1}\mbox{ }_{\partial _cB_{R_4}}P^n
[f(\zeta )(\zeta -\xi )^{-k-1}d\zeta ]$ for each $0\le k\in \bf Z$.
\par If $\zeta \in \partial _cB_{R_3}$, then $|(\zeta -\xi )(z-\xi )^{-1}|
<1$ and $(\zeta -z)^{-1}=-\sum_{k=1}^{\infty }
(\zeta -\xi )^{k-1}(z-\xi )^{-k}$, hence due to continuity of $P^n$ \\
$-C(\alpha )^{-1}\mbox{ }_{\partial _cB_{R_3}}P^n[f(\zeta )(\zeta -z)^{-1}
d\zeta ]=\sum_{k=1}^{\infty }a_{-k}(z-\xi )^{-k}$, \\
where $a_{-k}=C(\alpha )^{-1}\mbox{ }_{\partial _cB_{R_3}}P^n[
f(\zeta )(\zeta -\xi )^{k-1}d\zeta ]$ for each $k\ge 1$. In view of
Theorem $2.4.6$ we get Formula $(2)$.
\par {\bf 2.7.9. Definitions.} A point $z\in A({\bf K}\oplus
\alpha {\bf K})$ is called an isolated critical point of a function
$f$, if there exists a set $B({\bf K}\oplus \alpha {\bf K},z,R)
\setminus \{ z \} $ for $z\ne A$, and $ \{ \zeta \in {\bf K}
\oplus \alpha {\bf K}:$ $R<|\zeta |<\infty \} $ for $z=A$, on which
$f$ is $(q,n)$-antiderivationally holomorphic. An isolated critical point
$z$ of the function $f$ is called removable, if there exists a limit
$\lim_{\zeta \to z}f(\zeta )=g\in Y$; it is called a pole if there
exists $\lim_{\zeta \to z}\| f(\zeta ) \| =\infty $; it is called
essentially critical point, if there exists neither finite nor
infinite limit point, when $\zeta $ tends to $z$.
\par {\bf 2.7.10. Theorem.} {\it Let $f$ satisfy $2.7.2.(i)$ on
$\Omega \setminus \{ z \} $. A point $z\in {\bf K}\oplus \alpha
{\bf K}$ is removable if and only if decomposition $2.7.8.(1)$
does not contain the main part:
\par $(i)$ $f(\zeta )=\sum_{k=0}^{\infty }a_k(\zeta -z)^k$.}
\par {\bf 2.7.11. Theorem.} {\it Let $f$ satisfy $2.7.2.(i)$ on
$\Omega \setminus \{ z \} $. An isolated critical point
$z\in {\bf K}\oplus \alpha {\bf K}$ is a pole if and only if
the main part of series $2.7.8.(1)$ contains only a finite
and positive number of nonzero terms:
\par $(i)$ $f(\zeta )=\sum_{k=-N}^{\infty }a_k(\zeta -z)^k$, $N>0$.}
\par {\bf 2.7.12. Theorem.} {\it Let $f$ satisfy $2.7.2.(i)$ for
$Y={\bf K}(\alpha )$ on $\Omega \setminus \{ z \} $. An isolated
critical point $z$ of $f$ is essentially critical if and only if
the main part of series $2.7.8.(1)$ in a neighborhood of $z$ contains
an infinite family $ \{ a_k\ne 0:$ $ k<0 \} $. If $z$ is an
essentially critical point of $f$, $r=2$, that is,
${\bf K}(\alpha )={\bf K}\oplus \alpha {\bf K}$, then for each
$\xi \in A{\bf K}(\alpha )$ there exists a sequence $\{ z_n: $
$n\in {\bf N} \}  $, $\lim_{n\to \infty }z_n=z$ such that
$\lim_{n\to \infty }f(z_n)=\xi $.}
\par The {\bf proof} of these latter three theorems is analogous
to the classical case (see, for example, \S II.7 \cite{shabat})
due to the given above Theorems $2.7.2$ and $2.7.8$.
\par {\bf 2.7.13. Definition.} Let $f\in C^{(q,n)}(\Omega ,Y)$
and $B:=B({\bf K}\oplus \alpha {\bf K},z,R)\subset \Omega $,
$0<R<\infty $, $f$ is $(q,n)$-antiderivationally holomorphic
on $B\setminus \{ z \} $, then
\par $(i)$ $res_zf:=C(\alpha )^{-1}\mbox{ }_{\partial _cB}P^n
[f(\zeta )d\zeta ]$ \\
is called the residue of $f$, where $Y$ is a Banach space
over $\bf L$ such that ${\bf K}(\alpha )\subset \bf L$.
\par {\bf 2.7.14. Theorem.} {\it Let $f$ satisfy $2.7.2.(i)$
on $\Omega \setminus \bigcup_{l=1}^{\nu } \{ z_l \} $ such that
$\partial \Omega $ does not contain critical points $z_l$
of $f$ and all of them are encompassed by $\partial \Omega $,
$\nu \in \bf N$. Then
\par $(i)$ $\mbox{ }_{\partial \Omega }P^n[f(\zeta )d\zeta ]=
C(\alpha )\sum_{z_l\in \Omega }res_{z_l}f$, \\
where  $res_{z_l}f$ is independent of $n$ and $R$ in $2.7.13$,
\par $(ii)$ $res_{z_l}f=a_{-1}$, $a_k$ is as in $2.7.8.(1)$.}
\par {\bf Proof.} In view  of Theorems $2.4.2$ and $2.4.6$
$C(\alpha )=C_n(\alpha )$ is independent of $n$ and
$res_{z_l}f$ is independent of $n$ and $R$. From $2.7.8.(2)$
it follows $(ii)$.
\par {\bf 2.7.15. Definition.} Let $f\in C^{(q,n)}(\Omega ,Y)$
and let $A\in \Omega \subset A({\bf K}\oplus \alpha {\bf K})$
be the isolated critical point of $f$, put
\par $res_Af:=-C(\alpha )^{-1}\mbox{ }_{\partial B}P^n[f(\zeta )d\zeta ]$.
\par {\bf 2.7.16. Theorem.} {\it Let $f$ satisfy $2.7.2.(i)$ on
$({\bf K}\oplus \alpha {\bf K})\setminus \bigcup_{l=1}^{\nu }
\{ z_l \} $, then
\par $(i)$ $res_Af+\sum_{l=1}^{\nu }res_{z_l}f=0$.}
\par {\bf Proof.} Take a ball $B_R:=B({\bf K}\oplus \alpha
{\bf K},0,R)$ of sufficiently large $0<R_0<R<\infty $
such that it contains all $\{ z_l:$ $l=1,...,\nu \} $,
$\Omega =A({\bf K}\oplus \alpha {\bf K})$ and $\kappa (\Omega )
\subset {\bf K}\oplus \alpha {\bf K}$ (see \S 2.2.4).
In view of Theorem $2.7.14$
\par $(ii)$ ${\partial _cB_R}P^n[f(\zeta )d\zeta ]=
\sum_{l=1}^{\nu }res_{z_l}f$ \\
and it is independent of $R$ for each $R>R_0$, $R<\infty $.
In accordance with Definition $2.7.15$ and Theorem $2.4.6$
\par $(iii)$ $\mbox{ }_{\partial _cB_R}P^n[f(\zeta )d\zeta ]
=-res_Af$. \\
Therefore, from $(ii,iii)$ it follows $(i)$.
\par {\bf 2.7.17. Definitions.} Let $f\in C^{(q,n)}(\Omega ,
{\bf K}(\alpha ))$ and let $f$ be $(q,n)$-antiderivationally holomorphic
on $B({\bf K}\oplus \alpha {\bf K},z,R)\setminus \{ z \} $,
where $\Omega \subset {\bf K}\oplus \alpha {\bf K}$, $f(z)\ne 0$. Then
\par $res_zf'(z)/f(z)$ is called the logarithmic residue of $f$ at the
point $z$. Let us count each zero and pole of $f$ a number of times
equal to its order.
\par A function $f$ is called $(q,n)$-antiderivationally meromorphic,
if it is $(q,n)$-antiderivationally holomorphic on $\Omega $ besides
a set of poles.
\par {\bf 2.7.18. Theorem.} {\it Let $f$ be $(q,n)$-antiderivationally
meromorphic on $\Omega $ and let $Log (f)$ satisfy $2.7.2.(i)$ on
$\Omega \setminus \bigcup_{l=1}^{\nu } \{ z_l \} $, where
$z_l$ is the pole of $f$ for each $l=1,...,\nu $, all zeros
and poles of $f$ are encompassed by $\partial \Omega $,
$Y={\bf K}(\alpha )$. Then
\par $(1)$ $N-P=C(\alpha )^{-1}\mbox{ }_{\partial \Omega }P^n
[d Log (f(\zeta )) ]$, \\
where $N$ and $P$ denote total numbers of zeros and poles in $\Omega $.}
\par {\bf Proof.} Since $\Omega $ is compact, then $N$ and $P$
are finite. In view of Theorem $2.7.8$
\par $(2)$ $C(\alpha )^{-1}\mbox{ }_{\partial \Omega }P^n[
d Log (f(\zeta ))]=\sum_{l=1}^{\nu }res_{z_l} Log (f)+
\sum_{l=1}^{\mu }res_{\xi _l} Log (f)$, \\
where $z_l$ is the pole of $f$ and $\xi _l$ is the zero
of $f$ for each $l$. On the other hand,
\par $f'(z)/f(z)=[k(z-\xi _l)^{k-1}\phi (z)+(z-\xi _l)^k\phi '(z)]
(z-\xi _l)^{-k}/ \phi (z)$
$=(z-\xi _l)^{-1}[k\phi (z)+(z-\xi _l)\phi '(z)]/\phi (z)$, \\
where $f(z)=(z-\xi _l)^k\phi (z)$, $k=k_l$ is the order of zero
$\xi _l$, $\phi (z)\ne 0$ in a neighborhood of $\xi _l$. Therefore,
\par $(3)$ $res_{\xi _l}Log (f)=k_l$ and $res_{z_l}Log (f)=-s_l$, \\
where $s_l$ is the order of pole $z_l$. Hence, from $(2,3)$ it follows
Formula $(1)$.
\par {\bf 2.8. Theorem.} {\it Let $\Omega $ be a clopen compact subset in,
$B(({\bf K}\oplus \alpha {\bf K})^m,y,R)$,
$0<R<p^{1/(1-p)}$. Suppose  \\
$(i)$ $f_j(z_1,...,z_{l-1},z_l+Exp (\eta ),z_{l+1},...,z_m):=
\psi _{j,l}(\eta ): \omega _{l,\epsilon }\to Y$ belongs to
$\mbox{ }_SC^{(q+1,n-1)}(\omega _{l,\epsilon },Y)$
for each $j, l=1,...,m$ and each $z=(z_1,...,z_m)\in \Omega $,
where $\omega _l=\{ \eta \in {\bf K}(\alpha ):$ $(z_1,...,z_{l-1},z_l
+Exp (\eta ), z_{l+1},...,z_m)\in \Omega \} $, $\omega _{l,\epsilon }
=\omega _l\setminus Log (B({\bf K}(\alpha ),z_l,\epsilon ))$,
$\epsilon =\epsilon _k$, $0<\epsilon _k$ for each $k\in \bf N$,
$\lim_{k\to \infty }\epsilon _k=0$, $0\le q\in \bf Z$, $1\le n\in \bf N$,
$Y$ is a Banach space over $\bf L$ such that ${\bf K}(\alpha )\subset \bf L$.
Assume:
\par $(1)\quad \partial f_j/\partial {\bar z}_l=
\partial f_l/\partial {\bar z}_j$ for each $j, l=1,...,m$. \\
Then there exists $u\in C^{(q,n-1)}(\Omega ,{\bf K}(\alpha ))$ such that
\par $(2)\quad \partial u(z)/\partial {\bar z}_j=f_j(z)$ for each
$j=1,...,m$ and each $z\in \tilde {\Omega }$ (see \S 2.7.2).}
\par {\bf Proof.} Define
\par $(3)\quad u(z):=
C(\alpha )^{-1}\sum_{j=1}^m{\bar z}_j\mbox{ }_{\partial \Omega _j}P^n
[f_j(z_1,...,z_{j-1},\zeta ,z_{j+1},...,z_m)(\zeta -z_j)^{-1}d\zeta ]$ \\
$-C(\alpha )^{-1}\mbox{ }_{\Omega }P^n
[f_1(\zeta ,z_2,...,z_m)(\zeta -z_1)^{-1}d{\bar {\zeta }}\wedge d\zeta ]$. \\
Hence $u(z)=
C(\alpha )^{-1}\sum_{j=1}^m{\bar z}_j\mbox{ }_{\partial \Omega _j}P^n[
f_j(z_1,...,z_{j-1},\zeta ,z_{j+1},...,z_m)(\zeta -z_j)^{-1}d\zeta ]$ \\
$+C(\alpha )^{-1}\mbox{ }_{U^2}P^n[(\chi _{\Omega _1}
f_1)(z_1-\eta ,z_2,...,z_m)\eta ^{-1}d({\overline {z_1-\eta }})
\wedge d(z_1-\eta )]$, where $\eta :=z_1-\zeta $ and we can take
$U=B({\bf K},0,R)$ such that $U^2+U^2=U^2$ and $U^2$
is the additive group,
$\Omega _j:= \{ \xi \in U^2:$ $(z_1,...,z_{j-1},\xi ,z_{j+1},...,z_m)
\in \Omega \} .$ Therefore, $u\in C^{(q,n-1)}(\Omega ,Y)$. Then \\
$\partial u/\partial {\bar z}_j=C(\alpha )^{-1}\mbox{ }_{\partial \Omega _j}
P^n[f_j(z_1,...,z_{j-1},\zeta ,z_{j+1},...,z_m)(\zeta -z_j)^{-1}d\zeta ]$ \\
$-C(\alpha )^{-1}\mbox{ }_{U^2}P^n[(\chi _{\Omega _1}
\partial f_1(\zeta ,z_2,...,z_m)/\partial {\bar z}_j)(\zeta -z_1)^{-1}
d{\bar {\zeta }}\wedge d\zeta ]$. In view of Condition $(1)$
and Formula $2.6.3.(6)$ and  \\
$f_j(z)=C(\alpha )^{-1}\mbox{ }_{\partial \Omega _j}
P^n[f_j(z_1,...,z_{j-1},\zeta ,z_{j+1},...,z_m)(\zeta -z_j)^{-1}d\zeta ]$ \\
$-C(\alpha )^{-1}\mbox{ }_{U^2}P^n[(\chi _{\Omega _j}
(\partial f_j(\zeta ,z_2,...,z_m)/\partial {\bar {\zeta }})
(\zeta -z_1)^{-1}d{\bar {\zeta }}\wedge d\zeta ]$
(see Theorem $2.4.2$) it follows $(2)$.
\section{Antiderivational representations of functions and
differential forms}
\par {\bf 3.1. Remark and Notation.} Let $\Omega $ be a clopen compact subset
in $({\bf K}\oplus \alpha {\bf K})^m$. Put
\par $(1)$ $w(z,\zeta ):=\sum_{j=1}^m(-1)^{j+1}
(\zeta _j-z_j)^{-1}d\zeta _j\wedge _{l\ne j}
[(\xi ({\bar {\zeta }}-{\bar z}))^{-1}
d_{\bar {\zeta }}\xi _l({\bar {\zeta }}-{\bar z})
\wedge (\xi (\zeta - z))^{-1} d_{\zeta }\xi _l(\zeta -z)]$ \\
for each $z\ne \zeta \in \Omega ^2$, 
where $(\alpha ')^j\ne \alpha ^t$ for each $j=1,...,q'$;
$t=1,...,q$, $m\le q'\le {\tilde m}(\alpha ')$,
$r\le q\le {\tilde m}(\alpha )$  (see \S \S 2.1 and 2.4.1); there are
constants $0<\epsilon _1<\epsilon _2<\infty $ such that
\par $(2)$ $\epsilon _1 |\pi |^{-s} |\zeta | \le 
| Log (\xi (\zeta )) | \le \epsilon _2 |\pi |^{-s} |\zeta |$
and $\xi (\zeta )\ne 0$ for each $\zeta \in \Omega -z$, $\xi (0)=1$,
where $z\in \Omega $, $\xi (\zeta )\in C^{(q,n)}(\Omega -z,{\bf C_p})$,
$Log (\xi (\zeta ))\in ({\bf K}\oplus \alpha {\bf K})^m$
for each $\zeta \in \Omega -z$,
here the embeddings are used: $({\bf K}\oplus \alpha {\bf K})^m
\hookrightarrow ({\bf K}(\alpha ))^m\hookrightarrow
{\bf K}(\alpha ,\alpha ')\hookrightarrow {\bf C_p}$;
\par $(3)$ $\xi $ is such that $d_{\zeta }w(z,\zeta )=0$ on
$\Omega \setminus \{ z \} $;
\par $(4)\quad s:=s(\zeta ):=
-ord_{{\bf K}(\alpha ,\alpha ')}(\zeta )$
for each $\zeta \in \Omega -z$, $z_j\ne \zeta _j$ for each $j$, \\
$\alpha '$ is the root of $1$ in $\bf C_p$ such that
${\bf K}(\alpha )^m$ is embedded into ${\bf K}(\alpha ,\alpha ')=
({\bf K}(\alpha ))(\alpha ')$, $|z|_{{\bf K}(\alpha ,\alpha ')}=
|\pi |^{-ord_{{\bf K}(\alpha ,\alpha ')}(z)}$, $\pi $ is the same as
in \S 2.1;
\par $(5)$ $\lim_{l\to \infty }\mbox{ }_{\partial _cB(({\bf K}
\oplus \alpha {\bf K})^m,z,|\pi |^l)} P^n[w(z,\zeta )]=:q_m\ne 0$.
If $f$ is a $1$-form of class $C^{(0,n-1)}$ we define:
\par $(6)$ $(B^n_{\Omega }f)(z):=q_m^{-1}
\mbox{ }_{\Omega }P^n[f(\zeta )\wedge w(z,\zeta )]$ for each
$z\in \Omega $ encompassed by $\partial \Omega $.  \\
If $f\in C^{(0,n-1)}(\Omega ,Y)$, we define
\par $(7)$ $(B^n_{\partial \Omega }f)(z):=
q_m^{-1}\mbox{ }_{\zeta \in \partial \Omega }P^n[f(\zeta ) w(z,\zeta )]$
for each $z\in \Omega $ encompassed by $\partial \Omega $, where
$Y$ is a Banach space over $\bf L$ such that ${\bf K}(\alpha )
\subset \bf L$. 
\par {\bf 3.2. Theorem.} {\it Let $\Omega $ be a clopen compact
subset in $B(({\bf K}\oplus \alpha {\bf K})^m,y,R)$, $0<R<p^{1/(1-p)}$,
$B^n_{\Omega }$, $B^n_{\partial \Omega }$ be given by \S 3.1,
$f\in C^{(q+1,n-1)}(\Omega ,Y)$. Then
\par $(1)$ $f(z)=(B^n_{\partial \Omega }f)(z)-
(B^n_{\Omega }{\bar {\partial }}f)(z)$  \\
for each $z=(z_1,...,z_m) \in \tilde {\Omega }$ (see \S 2.7.2)
such that  \\
$(2)$ $(fw)(z_1,...,z_{l-1},z_l+Exp (\eta ),z_{l+1},...,z_m)=:
{\tilde \psi }_l(\eta )$ \\
$\in \mbox{ }_SC^{(q+1,n-1)}(\omega _{l,\epsilon },
L(\Lambda {\bf K}(\alpha ),Y))$ for each $l=1,...,m$ and each $\epsilon
=\epsilon _j$, where
$\omega _l:= \{ \eta \in {\bf K}(\alpha ):$ $(z_1,...,z_{l-1},
z_l+Exp (\eta ),z_{l+1},...,z_m)\in \Omega \} $,
$\omega _{l,\epsilon }:=\omega _l\setminus Log
(B({\bf K}(\alpha ),z_l,\epsilon ))$,
$0<\epsilon _j$ for each $j\in \bf N$, $\lim_{j\to \infty }\epsilon _j=0$,
$0\le q\in \bf Z$, $1\le n\in \bf N$ (see \S 2.4.1).}
\par {\bf Proof.} Fix $z\in \tilde {\Omega }$. In the particular case
$\xi ({\bar \zeta }-{\bar z})=Exp ( \pi ^{-s}(
{\bar \zeta }-{\bar z}))$ properties $3.1.(2,3)$ are
satisfied and $q_m=C(\alpha )m(2\alpha )^{m-1}$ due to
Formulas $2.1.(2,3)$, \S 2.2.4 and \S 2.4.2, since
$d{\bar z}\wedge dz=2\alpha dx\wedge dy$ and $\mbox{ }_UP^n[dx]|^b_a=b-a$
for each $a, b\in U$, where $z=x+\alpha y$, $x, y\in U\subset \bf K$.
Therefore, the family of such $\xi $ and $w$ satisfying Conditions
$3.1.(1-5)$ is nonvoid.
In view of $3.1.(3)$ $d_{\zeta }w(z,\zeta )=0$ on $\Omega
\setminus \{ z \} $, hence $d(f(\zeta )w(z,\zeta ))={\bar \partial }
f(\zeta )\wedge w(z,\zeta )$ on $\Omega \setminus \{ z \} $,
since $\partial f(\zeta )\wedge w(z,\zeta )=0$ on
$\Omega \setminus \{ z \} $. From Corollary $2.3.2$ it follows, that
there exists $\delta >0$ such that for each $0<\epsilon <\delta $,
$\epsilon \in \Gamma _{\bf K}$, there is the inclusion
$B({\epsilon }):=
B(({\bf K}\oplus \alpha {\bf K})^m,z,\epsilon )\subset \Omega $
and the equality is satified: 
\par $\mbox{ }_{\zeta \in \partial_cB({\epsilon })}P^n[f(\zeta )w(z,\zeta )]
=\mbox{ }_{\partial \Omega }P^n[f(\zeta )w(z,\zeta )]-
\mbox{ }_{\Omega ({\epsilon })}P^n[{\bar {\partial }}
f(\zeta )\wedge w(z,\zeta )]$ \\
for each $z\in \tilde {\Omega }$ and satisfying Condition $2.3.(2)$, where
$\Omega ({\epsilon }):= \{ \zeta \in \Omega :$ $|\zeta -z|\ge \epsilon \} $.
In the particular case, $\xi ({\bar {\zeta }}-{\bar z})=
Exp (\pi ^{-s}({\bar {\zeta }}-{\bar z}))$
due to $2.3.2$: \\
$\mbox{ }_{\zeta \in \partial _cB({\epsilon })}P^n[f(\zeta )w(z,\zeta )]=
\pi ^{-2(m-1)s}\mbox{ }_{\zeta \in \partial_cB({\epsilon })}P^n[
f(\zeta )\sum_{j=1}^m(-1)^{j+1}(\zeta _j-z_j)^{-1}d\zeta _j
\wedge _{l\ne j} (d{\bar {\zeta }}_l\wedge d\zeta _l)]$ \\
$= \pi ^{-2(m-1)s}(2\alpha )^{m-1}(\sum_{j=1}^m(-1)^{j+1}
\mbox{ }_{B(({\bf K}\oplus \alpha {\bf K})^{m-1},z',\epsilon )}
P^n \{ \mbox{ }_{\zeta _j\in \partial_cB({\bf K}\oplus \alpha {\bf K},z_j,
\epsilon )}P^n[f(\zeta )(\zeta _j-z_j)^{-1}d\zeta _j]\wedge _{l\ne j}
(dx_{2l-1}\wedge dx_{2l}) \} )$, \\
where $z'=(z_1,...,z_{j-1},z_{j+1},...,z_m)$,
$x_{2l-1}, x_{2l}\in U$, $z_l=x_{2l-1}+\alpha x_{2l}$ for each $l=1,...,m$.
Therefore, there exists
\par $\lim_{\epsilon \to 0}\mbox{ }_{\zeta \in \partial_cB({\epsilon })}
P^n[f(\zeta )w(z,\zeta )]=f(z)q_m$ \\
due to $3.1.(2,5)$, since there exists $C=const >0$, $C<\infty $,
such that 
\par $| \mbox{ }_{B({\epsilon })}P^n[f(\zeta )-f(z))w(z,\zeta )] |
\le C\epsilon \| f(\zeta )-f(z) \|_{C^{(0,n-1)}(B({\epsilon }),Y)}$  \\
for each $0<\epsilon <\delta $.
\par {\bf 3.3. Corollary.} {\it Let $\Omega $ and $f$ be as
in Theorem $3.2$ and $f$ be derivationally $(q,n)$-holomorphic on
$\Omega $, then
\par $(1)$ $f(z)=(B^n_{\partial \Omega }f)(z)$ for each $z\in
\tilde {\Omega }$ (see \S 2.7.2).}
\par {\bf 3.4. Remark.} For $m=1$ Formula $3.2.(1)$ is Formula
$2.4.2.(1)$, which are the non-Archimedean analogs of the
Martinelli-Bochner and Cauchy-Green formulas respectively
(see for comparison the classical complex case in \cite{henlei}).
\par {\bf 3.5. Definitions and Notations.} Consider a clopen
compact $\Omega \subset ({\bf K}\oplus \alpha {\bf K})^m$
and a $C^{(q,n+1)}$-function $v: \Omega \times (\partial \Omega )^{\delta }
\to ({\bf K}\oplus \alpha {\bf K})^m$,
$v=(v_1,...,v_m)$, $0<\delta <\infty $, $v=v(z,\zeta )$, $z\in
\Omega $, $\zeta \in (\partial \Omega )^{\delta }$,
$\Psi ^{\epsilon }:= \{ z\in X:$ $d(z,\Psi )< \epsilon \} $ for a
topological space $X$ with a metric $d$ and a subset $\Psi \subset
X$, $d(z, \Psi ):=\inf_{x\in \Psi }d(z,x)$, $0<\epsilon $,
$0\le q\in \bf Z$, $1\le n\in \bf N$. Suppose
\par ${\tilde {\phi }}(s):=-ord_{{\bf K}(\alpha ,\alpha ')}v(z,\zeta )$,
$s=-ord_{{\bf K}(\alpha ,\alpha ')}({\bar {\zeta }}-{\bar z})$ such that
\par $(1)$ $\lim_{s\to \infty }{\tilde {\phi }}(s)=\infty $ \\
for each $z\in \Omega $ and $\zeta \in (\partial \Omega )^{\delta }$.
Put
\par $\eta ^v(z,\zeta ,\lambda ):=(1-\lambda /\beta )\xi 
(v(z,\zeta ))+\lambda \xi ({\bar {\zeta }}-{\bar z})/\beta $, \\
where $\lambda \in B({\bf K},0,1)$. Impose the condition:
\par $(2)$ $\wedge _{k=1}^md_{\zeta }v_k(z,\zeta )\wedge _{j=1}^md\zeta _j
\ne 0$ and $\eta ^v(z,\zeta ,\lambda )\ne 0$ \\
for each $z\in \Omega $ and $\zeta \in (\partial \Omega )^{\delta }$
and $\lambda \in B({\bf K},0,1)$. Let also
\par $(3)$ $\psi (z,\zeta ):=\sum_{j=1}^m(-1)^{j+1}(\zeta _j-z_j)^{-1}
d\zeta _j\wedge _{k\ne j}[(\xi (v(z,\zeta )))^{-1}
d_{\bar {\zeta }}\xi _k(v(z,\zeta )) \wedge
(\xi (\zeta -z))^{-1}d_{\zeta }\xi _k(\zeta -z)]$. \\
If $f\in C^{(0,n-1)}(\partial \Omega ,Y)$, we set:
\par $(4)$ $(L^{v,n}_{\partial \Omega }f)(z):=q_m^{-1}
\mbox{ }_{\zeta \in \partial \Omega }P^n[f(\zeta )\psi (z,\zeta )]$ \\
for each $z\in \Omega $, where $Y$ is a Banach space over $\bf L$
such that ${\bf K}(\alpha )\subset \bf L$. Put also:
\par $(5)$ $\gamma (z,\zeta ,\lambda ):=\sum_{j=1}^m
(-1)^{j+1}(\zeta _j-z_j)^{-1}d\zeta _j\wedge _{k\ne j}
[(\eta ^v(z,\zeta ,\lambda ))^{-1}({\bar {\partial }}_{z,\zeta }+
d_{\lambda })\eta ^v_k(z,\zeta ,\lambda )\wedge
(\xi (\zeta -z))^{-1}d_{\zeta } \xi _k(\zeta -z)]$ \\
for each $z\in \Omega $, $\zeta \in (\partial \Omega )^{\delta }$,
$\lambda \in B({\bf K},0,1)$.
If $f$ is a $C^{(0,n-1)}$-$1$-form on $\partial \Omega $ put:
\par $(6)$ $(R^{v,n}_{\partial \Omega }f)(z):=q_m^{-1}
\mbox{ }_{\zeta \in \partial \Omega , \lambda \in B({\bf K},0,1)}
P^n[f(\zeta )\wedge \gamma (z,\zeta ,\lambda )]$ \\
for each $z\in \Omega $. Suppose that $v$ is such that
\par $(7)$ $d_{z,\zeta }\psi (z,\zeta )=0$ for each $\zeta \ne z$.
In particular, if $v(z,\zeta )={\bar {\zeta }}-{\bar z}$, then
\par $L^{v,n}_{\partial \Omega }f=B^n_{\partial \Omega }f$ and
$R^{v,n}_{\partial \Omega }f=0$, \\
since $\gamma $ is the $(m,m)$-form by $(\zeta ,{\bar {\zeta }})$
for $v(z,\zeta )={\bar {\zeta }}-{\bar z}$.
\par {\bf 3.6. Theorem.} {\it Let $\Omega $ be a clopen compact subset
in $({\bf K}\oplus \alpha {\bf K})^m$ and let $v(z,\zeta )$,
$L^{v,n}_{\partial \Omega }$ and $R^{v,n}_{\partial \Omega }$ be given by
\S 3.5, $f\in C^{(q+1,n-1)}(\Omega ,Y)$. Then
\par $(1)$ $f(z)=(L^{v,n}_{\partial \Omega }f)(z)-(R^{v,n}_{\partial
\Omega }{\bar {\partial }}f)(z)-(B^n_{\Omega }{\bar {\partial }}
f)(z)$ \\
for each $z\in \tilde {\Omega }$ (see \S 2.7.2) such that
\par $(2)$ $(f\gamma )(z_1,...,z_{l-1},z_l+Exp (\eta ),z_{l+1},...,z_m)
=:{\tilde {\psi }}_l(\eta )$ \\
$\in \mbox{ }_SC^{(q+1,n-1)}(\omega _{l,\epsilon },
L(\Lambda {\bf K}(\alpha ),Y))$ \\
for each $l=1,...,m$, $\epsilon =\epsilon _j$, where
$\omega _l:= \omega _l(z):=
\{ \eta \in {\bf K}(\alpha ):$ $(z_1,...,z_{l-1},z_l+ Exp (\eta ),
z_{l+1},...,z_m) \in \Omega \} $, $\omega _{l,\epsilon }
:=\omega _l\setminus Log (B({\bf K}(\alpha ),z_l,\epsilon ))$,
$0<\epsilon _j\in \Gamma _{\bf K}$ for each $j\in \bf N$,
$\lim_{j\to \infty }\epsilon _j=0$, $0\le q \in \bf Z$, $1\le n\in \bf N$
(see \S 2.4.1).}
\par {\bf Proof.} The using of Theorem 3.2 reduces the proof to that
of the formula:
\par $(3)$ $(R^{v,n}_{\partial \Omega }{\bar {\partial }}f)(z)=
(L^{v,n}_{\partial \Omega }f)(z)-(B^n_{\partial \Omega }f)(z)$ \\
for each $z\in \tilde {\Omega }$ and
satisfying Condition $(2)$. In view of $3.2.(2)$ and $3.5.(7)$ we have
$d_{\zeta ,\lambda }\gamma (z,\zeta ,\lambda )=0$, since
$d_{\zeta ,\lambda }[d_{\zeta ,\lambda }(\eta ^v)]=0$.
Therefore, $d_{\zeta ,\lambda }[f(\zeta )\gamma (z,\zeta ,\lambda )]=
({\bar {\partial }}f(\zeta ))\wedge \gamma (z,\zeta ,\lambda )$,
since $(\partial f)\wedge \gamma =0$. From $3.5.(3,5)$ it follows that
\par $\gamma (z,\zeta ,\lambda )|_{\lambda =0}=\psi (z,\zeta )$ and
\par $\gamma (z,\zeta ,\lambda )|_{\lambda =\beta }=
\sum_{j=1}^m(-1)^{j+1}(\zeta _j-z_j)^{-1}d\zeta _j\wedge _{k\ne j}
[(\xi ({\bar {\zeta }}-{\bar z}))^{-1}d_{\bar {\zeta }}
\xi _k({\bar {\zeta }}-{\bar z}) \wedge
(\xi (\zeta - z))^{-1}d_{\zeta } \xi _k(\zeta -z)] $. \\
Mention that $\lambda =P^n1|_0^{\lambda }$, hence $\lambda
\in \mbox{ }_PC^{(q,n)}(B({\bf K},0,1),{\bf K})$. Then by the degree
reasons
\par $\mbox{ }_{\zeta \in \partial \Omega }P^n[f\gamma |_{\lambda =\beta }]
=\mbox{ }_{\zeta \in \partial \Omega }P^n[fw]$, \\
where $w$ is given by $3.1.(1)$, $dim_{\bf K}\Omega =2m$,
$dim_{\bf K}\partial \Omega =2m-1$. Then
\par $q_m^{-1}\mbox{ }_{\zeta \in \partial \Omega ,\lambda \in B}
P^n\{ d_{\zeta ,\lambda }[f(\zeta )\gamma (z,\zeta ,\lambda )] \}$ \\
$=q_m^{-1}\mbox{ }_{\zeta \in \partial \Omega ,\lambda \in B}
P^n\{ ( {\bar {\partial }}f(\zeta )\wedge \gamma (z,\zeta ,\lambda ) \}
=(R^{v,n}_{\partial \Omega }{\bar {\partial }}f)(z)$ \\
for each $z\in \tilde {\Omega }$,
where $B:=B({\bf K},0,1)$. On the other hand,
\par $\partial ((\partial \Omega )\times B)=(-1)^{2m-1}
((\partial \Omega )\times \{ \beta \} -(\partial \Omega )\times
\{ 0 \} )=-(\partial \Omega )\times \{ \beta \} +
(\partial \Omega )\times \{ 0 \} $. In view of Corollary $2.3.2$
\par $\mbox{ }_{(\partial \Omega )\times B}P^n \{ d_{\zeta ,\lambda }
[f\gamma ] \} =-\mbox{ }_{\partial \Omega }P^n[fw]+
\mbox{ }_{\partial \Omega }P^n[f\psi ]$, \\
hence Formula $(3)$ is accomplished.
\par {\bf 3.7. Corollary.} {\it Let $f$ be as in Theorem $3.6$
and ${\bar {\partial }}f=0$ on $\Omega $, then
\par $(i)$ $f(z)=(L^{v,n}_{\partial \Omega }f)(z)$ \\
for each $z\in \tilde {\Omega }$ (see \S 2.7.2)
and satisfying Condition $3.6.(2)$.}
\par {\bf 3.8. Definitions and Remarks.} Let $\Omega $
be a clopen compact subset in $({\bf K}\oplus \alpha {\bf K})^m$,
consider the differential form:
\par $(1)$ ${\tilde w}(z,\zeta ):=\sum_{j=1}^m(-1)^{j+1}
(\zeta _j-z_j)^{-1}d\zeta _j\wedge _{l\ne j} [(\xi ({\bar {\zeta }}-
{\bar z}))^{-1}{\bar {\partial }}_{\zeta ,z}\xi _l({\bar {\zeta }} -
{\bar z})\wedge (\xi (\zeta -z))^{-1}d\xi _l(\zeta -z)]$.
Let $M$ be a compact manifold over $\bf K$ and $\phi : \Omega \to M
\hookrightarrow ({\bf K}(\alpha ))^N$ be a
$\mbox{ }_SC^{(q+1,n-1)}$-diffeomorphism (see \S 2.2.5). Then the
diffeomorphism $\phi _*w$ of the differential form $w$ is the
differential form on $M$. Consider these differential
forms on $M$ also and denote them by the same notation, since
$\{ \phi (\zeta _j):$ $j \} $ are coordinates in $M$.
Therefore, Theorems $3.2, 3.6$ and Corollaries $3.3, 3.7$
are true for $M$ also due to Theorem $2.3.1$ and Corollary
$2.3.2$, where ${\tilde M}:=\phi ^{-1}(\tilde {\Omega })$
(see \S \S 2.2.5 and 2.7.2). If $f$ is a
$C^{(0,n-1)}$-differential form on $M$, then we define:
\par $(2)$ $(B^n_Mf)(z):=q_m^{-1}\mbox{ }_{\zeta \in M}P^n[
f(\zeta )\wedge {\tilde w}(z,\zeta )]$ \\
for each $z\in M$ encompassed by $\partial M$. If $f$ is a
$C^{(0,n-1)}$-differential form on $M$, then we define:
\par $(3)$ $(B^n_{\partial M}f)(z):=q_m^{-1}\mbox{ }_{
\zeta \in \partial  M}P^n[f(\zeta )\wedge {\tilde w}(z,\zeta )]$ \\
for each $z\in M$ encompassed by $\partial M$. Write $\tilde w$ as:
\par $(4)$ ${\tilde w}(z,\zeta )=\sum_{t=0}^{m-1}\Upsilon _t(z,\zeta )$, \\
where $\Upsilon _t$ is of bedegree $(0,t)$ in $z$ and of bedegree
$(m,m-t-1)$ in $\zeta $. Decompose $f$ as:
\par $(5)$ $f=\sum_{l+s=deg (f)}f_{(l,s)}$, \\
where $f_{(l,s)}$ is the $(l,s)$-form on $M$. Then $f_{(l,s)}(\zeta )
\wedge _{j=1}^md\zeta _j=0$ for each $l>0$, hence
\par $B^n_Mf=B^n_Mf_{(0,deg (f))}$. On the other hand,
$f(\zeta )\wedge \Upsilon _t(z,\zeta )=0$, if $deg (f)>q+1$;
\par $\mbox{ }_{\zeta \in M}P^n[f(\zeta )\wedge \Upsilon _t(z,\zeta )]=0$,
when $deg (f)<q+1$ by the definition of the antiderivation. Therefore,
\par $(6)$ $B^n_Mf=\mbox{ }_{\zeta \in M}P^n[f_{(0,deg (f))}(\zeta )
\wedge \Upsilon _{deg (f)-1}(z,\zeta )]$ for $1\le deg (f)\le m$,
\par $(7)$ $B^n_Mf=0$ for $deg (f)=0$ or $deg (f)>m$, similarly
\par $(8)$ $B^n_{\partial M}f=\mbox{ }_{\zeta \in \partial M}P^n[
f_{(0,deg (f))}(\zeta )\wedge \Upsilon _{deg (f)}(z,\zeta )]$
for $0\le deg (f)\le m-1$, 
\par $(9)$ $B^n_{\partial M}f=0$ for $deg (f)\ge m$, \\
hence $B^n_Mf$ is of bidegree $(0,deg (f)-1)$;
$P^n_{\partial M}f$ is of bidegree $(0,deg (f))$. Using the
notation of \S 3.5 define:
\par $(10)$ ${\tilde \psi }(z,\zeta ):=\sum_{j=1}^m(-1)^{j+1}
(\zeta _j-z_j)^{-1}d\zeta _j\wedge _{k\ne j}[(\xi (v(z,\zeta )))^{-1}
{\bar {\partial }}_{z,\zeta }\xi _k(v(z,\zeta ))\wedge (\xi (\zeta -z))^{-1}
d_{\zeta }\xi _k(\zeta -z)]$,
\par $(11)$ ${\tilde {\gamma }}(z,\zeta ,\lambda ):=
\sum_{j=1}^m(-1)^{j+1}(\zeta _j-z_j)^{-1}d\zeta _j\wedge _{k\ne j}
[(\eta ^v(z,\zeta ,\lambda ))^{-1}
({\bar {\partial }}_{z,\zeta }+d_{\lambda })\eta ^v(z,\zeta ,\lambda )
\wedge (\xi (\zeta -z))^{-1}d_{\zeta }\xi _k(\zeta -z)]$. \\
If $f\in C^{(0,n-1)}(\partial M,L(\Lambda {\bf K}(\alpha ),Y))$ put:
\par $(12)$ $(L^{v,n}_{\partial M}f)(z):=
q_m^{-1}\mbox{ }_{\zeta \in \partial M}P^n[f(\zeta )\wedge
{\tilde {\psi }}(z,\zeta )]$  and
\par $(13)$ $(R^{v,n}_{\partial M}f)(z):=
q_m^{-1}\mbox{ }_{\zeta \in \partial M, \lambda \in B}P^n[f(\zeta )\wedge
{\tilde {\gamma }}(z,\zeta ,\lambda )]$  \\
for each $z\in {\tilde M}:=\phi ^{-1}(\tilde {\Omega })$
(see \S \S 2.2.5 and 2.7.2). There exists the
decomposition:
\par $(14)$ ${\tilde {\gamma }}(z,\zeta ,\lambda )=\sum_{t=0}^{m-1}
\Upsilon ^v_t(z,\zeta ,\lambda )$, \\
where $\Upsilon ^v_t(z,\zeta ,\lambda )$ is of bidegree $(0,t)$ in $z$
and $f$ is of bidegree $(m,m-t-1)$ in $(\zeta ,\lambda )$. Let $f$
be a bounded differential form on $\partial M$, $f=\sum f_{(l,s)}$, then
$R^{v,n}_{\partial M}f=R^{v,n}_{\partial M}f_{(0,deg (f))}$ and
$f(\zeta )\wedge \Upsilon ^v_t(z,\zeta ,\lambda )=0$
if $deg (f)>t+1$ by the definition of $P^n$,
$\mbox{ }_{(\zeta ,\lambda )\in {\partial M}\times B}P^n[
f(\zeta )\wedge \Upsilon ^v_t(z,\zeta ,\lambda )]=0$
if $deg (f)<t+1$, hence
\par $(15)$ $R^{v,n}_{\partial M}f=\mbox{ }_{\zeta \in \partial M,
\lambda \in B}P^n[f_{(0,deg (f))}(\zeta )\wedge \Upsilon ^v_{deg (f)-1}
(z,\zeta ,\lambda )]$ if $1\le deg (f)\le m$,
\par $(16)$ $R^{v,n}_{\partial M}f=0$ if $deg (f)=0$ or $deg (f)>m$.
Similarly,
\par $(17)$ ${\tilde {\psi }}(z,\zeta )=\sum_{t=0}^{m-1}
\Upsilon ^v_t(z,\zeta )$, \\
where $\Upsilon ^v_t(z,\zeta )$ is of bidegree $(0,t)$ in $z$ and
of bidegree $(m,m-t-1)$ in $\zeta $, hence
\par $(18)$ $L^{v,n}_{\partial M}f=\mbox{ }_{\zeta \in \partial M}
P^n[f_{(0,deg (f))}\wedge \Upsilon ^v_{deg (f)}(z,\zeta )]$
if $deg (f)\le m-1$,
\par $(19)$ $L^{v,n}_{\partial M}f=0$ if $deg (f)\ge m$.
If $v(z,\zeta )={\bar {\zeta }}-{\bar z}$, then $L^{v,n}_{\partial M}f=
B^n_{\partial M}f$.
\par {\bf 3.9. Theorem.} {\it Let $M$ be a compact manifold
over $\bf K$ and let $B^n_M$ and $B^n_{\partial M}$ be given
by \S 3.8. Suppose that $f$ is the $C^{(q+1,n-1)}$-$(0,t)$-form,
$0\le t\le m$. Then
\par $(1)$ $(-1)^tf(z)=(B^n_{\partial M}f)(z)-(B^n_M{\bar {\partial }}
f)(z)+({\bar {\partial }}B^n_Mf)(z)$ \\
for each $z\in \tilde M$ such that
\par $(2)$ $(f\wedge {\tilde w})\circ \phi (z_1,...,z_{l-1},
z_l+ Exp (\eta ),z_{l+1},...,z_m)=:{\tilde {\psi }}_l
(\eta )$ \\ $
\in \mbox{ }_SC^{(q+1,n-1)}(\omega _{l,\epsilon },
L(\Lambda ({\bf K}(\alpha )),Y))$ for each $l=1,...,m$
and each $\epsilon =\epsilon _j$, where
$\omega _l:= \{ \eta \in {\bf K}(\alpha ):$ $(z_1,...,z_{l-1},
z_l+Exp (\eta ),z_{l+1},...,z_m)\in \Omega \} $,
$\omega _{l,\epsilon }:=\omega _l\setminus Log (B({\bf K}(\alpha ),
z_l,\epsilon ))$, $0<\epsilon _j$ for each $j\in \bf N$,
$\lim_{j\to \infty } \epsilon _j=0$, $0\le q\in \bf Z$,
$1\le n\in \bf N$.}
\par {\bf Proof.} Using the diffeomorphism $\phi $ it is possible
reduce the case to $\Omega \subset ({\bf K}\oplus \alpha {\bf K})^m$.
If $q=0$, then by $3.8.(7)$ $B^n_{\Omega }f=0$ and $f=
B^n_{\partial \Omega }f-B^n_{\Omega }{\bar {\partial }}f$ is Formula
$3.2.(1)$. Since $(2)$ is satisfied, $v$ and $\xi \in C^{(q,n)}$,
then $B^n_{\partial \Omega }f$ and $B^n_{\Omega }{\bar {\partial }}f$
are in $\mbox{ }_SC^{(q,n)}(\omega _{l,\epsilon },
L(\Lambda ({\bf K}(\alpha )),Y))$
for each $l=1,...,m$, $\epsilon =\epsilon _j$.
From the definition of $B^n_M$ it follows, that
$\sup_u \| {\bar {\Phi }}^u(B^n_{\Omega }{\bar {\partial }}f)(z;
h_1^{\otimes u_1},...,h_m^{\otimes u_m}; \zeta _1,...,\zeta _u)-
{\bar {\Phi }}^u(B^n_{\Omega }{\bar {\partial }}f)(y;
h_1^{\otimes u_1},...,h_m^{\otimes u_m}; \zeta _1,...,\zeta _u) \|_{
C^{(q,n)}} \le C_1 \| f \| _{C^{(q,n)}}|1-\pi ^{-2sm}|$,
where $s=s(\zeta -z)$, $u=u_1+...+u_m$,
$0\le u_l\le n$, hence \\
$(B^n_M{\bar {\partial }}f)\in C^{(q,n)}(M,L(\Lambda ({\bf K}(\alpha )),Y))$.
Analogously, $B^n_{\partial M}f$ and $B^n_Mf$ are in $C^{q,n)}(M,
L(\Lambda ({\bf K}(\alpha ),Y))$. It remains to prove, that in the
sence of distributions:
\par $({\bar {\partial }}B^n_{\Omega }f)(z)=(-1)^tf(z)-
(B^n_{\partial M}f)(z)+(B^n_M{\bar {\partial }}f)(z)$ \\
for each $z\in \tilde {\Omega }$ and satisfying Condition $(2)$.
This means, that for each
$\mbox{ }_SC^{(q+1,n-1)}$-form $\nu $, $supp (\nu )\subset
\tilde {\Omega }$, there is satisfied the equality:
\par $(3)$ $(-1)^t\mbox{ }_{\Omega }P^n[B^n_{\Omega }f\wedge
{\bar {\partial }}\nu ]=(-1)^t\mbox{ }_{\Omega }P^n[
f\wedge \nu ]$ \\
$-\mbox{ }_{\Omega }P^n[B^n_{\partial \Omega }f\wedge \nu ]
+\mbox{ }_{\Omega }P^n[B^n_{\Omega }({\bar {\partial f}})\wedge \nu ].$
In view of Formulas $3.8.(6,8)$ $B^n_{\partial \Omega }f$ and
$B^n_{\Omega }{\bar {\partial }}f $ are of bidegree $(0,t)$ and
$B^n_{\Omega }f$ is of bidegree $(0,t-1)$, we can assume that $\nu $ is
of bidegree $(m,m-t)$. Then $(3)$ takes
the form:
\par $(4)$ $(-1)^t\mbox{ }_{(\zeta ,z)\in \Omega ^2 } P^n[
f(\zeta )\wedge {\tilde w}(z,\zeta )\wedge {\bar {\partial }}\nu (z)]=
(-1)^t\mbox{ }_{z\in \Omega }P^n[f(z)\wedge \nu (z)]$ \\
$-\mbox{ }_{(\zeta ,z)\in (\partial \Omega )\times \Omega } P^n[
f(\zeta )\wedge {\tilde w}(z,\zeta )\wedge \nu (z)]+
\mbox{ }_{(\zeta ,z)\in \Omega ^2} P^n[{\bar {\partial }}f(\zeta )
\wedge {\tilde w}(z,\zeta )\wedge \nu (z)]$. Put
\par $(5)$ ${\tilde {\theta }}(z,\zeta ):=\sum_{j=1}^m(-1)^{j+1}
(\zeta _j-z_j)^{-1}(d\zeta _j-dz_j)\wedge _{k\ne j}
[(\xi ({\bar {\zeta }}-{\bar z}))^{-1}{\bar {\partial }}_{\zeta ,z}
\xi _k({\bar {\zeta }}-{\bar z})\wedge (\xi (\zeta -z))^{-1}
\partial _{\zeta ,z}\xi _k(\zeta -z)]$, \\
then from $3.1.(1,3)$ and $d_{\zeta ,z}=d_{\zeta }+d_z$ it follows:
\par $(6)$ $d_{\zeta ,z}{\tilde {\theta }}(z,\zeta )=0$ for
$\zeta \ne z$, since ${\bar {\partial }}(\zeta _j-z_j)=0$,
$d_{\zeta }^2=0$.
\par Then all monomials in ${\tilde {\theta }}(z,\zeta )-
{\tilde w} (z,\zeta )$ contain at least one of the differentials
$dz_1,...,dz_m$. For $\nu (z)$ of bidegree $(m,m-t)$ it contains
the factor $dz_1\wedge ... \wedge dz_m$, hence from $(5,6)$ it
follows:
\par $(7)$ $d_{\zeta ,z}({\tilde w}(z,\zeta )\wedge \nu (z))=
d_{z,\zeta }({\tilde {\theta }}(z,\zeta )\wedge \nu (z))$  \\
$=(-1)^{2m-1}{\tilde {\theta }}(z,\zeta )\wedge d\nu (z)=
-{\tilde w}(z,\zeta ) \wedge {\bar {\partial }}\nu (z)$ \\
for $\zeta \ne z$, since ${\tilde w}(z,\zeta )$ contains
the factor $d\zeta _1\wedge ... \wedge d\zeta _m$.
Hence $(7)$ implies:
\par $(8)$ $d_{z,\zeta }(f(\zeta )\wedge {\tilde w}(z,\zeta )\wedge
\nu (z))=({\bar {\partial }}f(\zeta ))\wedge {\tilde w}(z,\zeta )
\wedge \nu (z)-(-1)^tf(\zeta )\wedge {\tilde w}(z,\zeta )\wedge
{\bar {\partial }}\nu (z)$ for $\zeta \ne z$. Then
\par $(9)$ $\partial (\Omega \times \Omega \setminus U({\epsilon }))
\cap [({\bf K}\oplus \alpha {\bf K})^m\times supp (\nu )]=
((\partial \Omega )\times \Omega )\cup (\Omega \times \partial \Omega )-
\partial U({\epsilon }))\cap [({\bf K}\oplus \alpha {\bf K})^m\times
supp (\nu )]$, where $U({\epsilon }):=\{ (\zeta ,z)\in ({\bf K}\oplus
\alpha {\bf K})^m\times ({\bf K}\oplus \alpha {\bf K})^m:$
$|\zeta -z|<\epsilon \} $, $0<\epsilon <\epsilon _0$,
$0<\epsilon _0<\infty $ is fixed. In view of Corollary $2.3.2$
and Formulas $(8,9)$:
\par $(10)$ $\mbox{ }_{(\partial \Omega )\times \Omega \cup
\Omega \times (\partial \Omega )}P^n[f(\zeta )\wedge {\tilde w}
(z,\zeta )\wedge \nu (z)]-\mbox{ }_{\partial U({\epsilon })}
P^n[f(\zeta )\wedge {\tilde w}(z,\zeta )\wedge \nu (z)]$ \\
$=\mbox{ }_{\Omega ^2\setminus U({\epsilon })}P^n[({\bar {\partial }}
f(\zeta ))\wedge {\tilde w}(z,\zeta )\wedge \nu (z)]-
(-1)^t\mbox{ }_{\Omega ^2\setminus U(\epsilon )}P^n[f(\zeta )
\wedge {\tilde w}(z,\zeta )\wedge {\bar {\partial }}\nu (z)].$ \\
For $B^{-}_{\epsilon }:= \{ \zeta \in ({\bf K}\oplus \alpha {\bf K})^m:$
$|\zeta |<\epsilon \} $ there exists $\partial B^{-}_{\epsilon }$
such that $T(\partial B^{-}_{\epsilon }\times ({\bf K}\oplus \alpha
{\bf K})^m)=\partial U({\epsilon })$, where $T(\zeta ,z):=(z+\zeta ,z)$,
$T: ({\bf K}\oplus \alpha {\bf K})^m\times ({\bf K}\oplus \alpha {\bf K})^m
\to ({\bf K}\oplus \alpha {\bf K})^{2m}$.
The differential form $\nu (z)$ contains the factor $dz_1\wedge
... \wedge dz_m$, hence ${\tilde w}(z,\zeta )\wedge \nu (z)=
w(z,\zeta )\wedge \nu (z)$ and $T^*(f(\zeta )\wedge {\tilde w}
(z,\zeta )\wedge \nu (z))=\sum_{|I|=t}f_I(z+\zeta )
d({\bar z}+{\bar {\zeta }})^{\wedge I}\wedge w(z,\zeta )\wedge \nu (z)$,
where $T^*$ is the pull-back operator on differential forms
(see \S 2.2.5). The degree of $w(z,\zeta )$ is $2m-1$ and
$2m-1=dim_{\bf K}(\partial B^{-}_{\epsilon })$, consequently,
$d({\bar z}+{\bar {\zeta }})^{\wedge I}\wedge w(z,\zeta )|_{(
\partial B^{-}_{\epsilon })\times ({\bf K}\oplus \alpha {\bf K})^m}
=d{\bar z}^{\wedge I}\wedge w(z,\zeta )|_{(\partial B^{-}_{\epsilon })
\times ({\bf K}\oplus \alpha {\bf K})^m}$. \\
Therefore, taking $R>0$ such that $\Omega \subset B(({\bf K}\oplus
\alpha {\bf K})^m,0,R)=:B_R$:
\par $\mbox{ }_{\partial U({\epsilon })}P^n[f(\zeta )\wedge
{\tilde w}(z,\zeta )]=(-1)^t\mbox{ }_{z\in B_R}P^n
[\sum_{|I|=t}\mbox{ }_{\zeta \in T^{-1}(\partial U({\epsilon }))}
P^n[f_I(z+\zeta )\wedge w(z,\zeta )]d{\bar z}^{\wedge I}\wedge \nu (z)$, \\
since $d{\bar z}^{\wedge I}\wedge w(z,\zeta )=(-1)^t{\tilde w}
(z,\zeta )d{\bar z}^{\wedge I}$ for $|I|=t$, $\Omega +\Omega
\subset B_R+B_R=B_R$, where $supp (f)\subset \Omega $
(see Lemma $2.6.1$). In view of Theorem $3.2$
\par $\mbox{ }_{\zeta \in T^{-1}(\partial U({\epsilon }))}
P^n[f_I(z+\zeta )\wedge w(z,\zeta )]=f_I(z)+
\mbox{ }_{\zeta \in T^{-1}(\partial U({\epsilon }))}P^n \{ [
f_I(z+\zeta )-f_I(z)]\wedge w(z,\zeta ) \} $ \\
for $|I|=t$, which tends to $f_I(z)$, when $\epsilon $ tends
to zero, since $supp (f)$ is bounded, where $T^{-1}
(\partial U({\epsilon }))=(\partial B^{-}_{\epsilon })\times
({\bf K}\oplus \alpha {\bf K})^m$, and inevitably \\
$\lim_{\epsilon \to 0} \mbox{ }_{\partial U({\epsilon })}P^n[
f(\zeta )\wedge {\tilde w}(z,\zeta )\wedge \nu (z)]=(-1)^t
\mbox{ }_{\Omega }P^n[f(z)\wedge \nu (z)]$.
\par {\bf 3.10. Theorem.} {\it Let $M$ be a compact manifold
and let $L^{v,n}_{\partial M}$, $R^{v,n}_{\partial M}$, $B^n_M$
be given by \S 3.8. Suppose $f$ is the $C^{(q,n)}$-$(0,t)$-differential
form, $0\le t\le m$. Then
\par $(1)$ $(-1)^tf(z)=(L^{v,n}_{\partial M}f)(z)-$ \\
$(R^{v,n}_{\partial M}{\bar {\partial }}f+B^n_M{\bar {\partial }}f)(z)+
{\bar {\partial }}(R^{v,n}_{\partial M}f+B^n_Mf)(z)$ \\
for each $z\in \tilde M$ such that
\par $(2)$ $(f\wedge {\tilde {\gamma }})\circ \phi (z_1,...,
z_{l-1},z_l+Exp (\eta ),z_{l+1},...,z_m)=:{\tilde {\psi }}_l(\eta )$ \\
$\in \mbox{ }_SC^{(q+1,n-1)}(\omega _{l,\epsilon },
L(\Lambda {\bf K}(\alpha ),Y))$ \\
for each $l=1,...,m$ and each $\epsilon =\epsilon _j$, where
$\omega _l:=\{ \eta \in {\bf K}(\alpha ):$ $(z_1,...,z_{l-1},
z_l+Exp (\eta ),z_{l+1},...,z_m)\in \Omega \} $, $\omega _{l,\epsilon }
:=\omega _l\setminus Log (B({\bf K}(\alpha ),z_l,\epsilon ))$,
$\epsilon _j>0$ for each $j\in \bf N$, $\lim_{j\to \infty }\epsilon _j=0$,
$0\le q\in \bf Z$, $1\le n\in \bf N$.}
\par {\bf Proof.} If $v(z,\zeta )={\bar {\zeta }}-{\bar z}$, then
$L^{v,n}_{\partial M}=B^n_{\partial M}$, $R^{v,n}_{\partial M}=0$
and Formula $3.10.(1)$ reduces to Formula $3.9.(1)$. If $t=0$, then
by $3.8.(7,16)$ $B^n_Mf=0$ and $R^{v,n}_{\partial M}f=0$, hence
$3.10.(1)$ reduces to $3.6.(1)$. Assume $1\le t\le m$. In view of
\S \S 3.8 and 3.9 $L^{v,n}_{\partial M}f$, $R^{v,n}_{\partial M}
{\bar {\partial }}f$, $B^n_M{\bar {\partial }}f$, ${\bar {\partial }}
R^{v,n}_{\partial M}f$ and ${\bar {\partial }}B^n_Mf$ are in
$C^{(q,n)}(M,L(\Lambda {\bf K}(\alpha ),Y))$. Using the diffeomorphism
$\phi $ consider $\Omega $ instead of $M$.
In view  of $3.9.(1)$ it remains to prove:
\par $(3)$ ${\bar {\partial }}(R^{v,n}_{\partial \Omega }f)(z)=
(B^n_{\partial \Omega }f)(z)-(L^{v,n}_{\partial \Omega }
f)(z)+(R^{v,n}_{\partial \Omega }{\bar {\partial }}f)(z)$ \\
for each $z\in \tilde {\Omega }$
and satisfying Condition $(2)$. Consider the differential form:
\par $(4)$ $\kappa :=\sum_{j=1}^m(-1)^{j+1}(\zeta _j-z_j)^{-1}
d\zeta _j\wedge _{k\ne j} [(\eta ^v(z,\zeta ,\lambda ))^{-1}
d_{z,\zeta ,\lambda }\eta ^v_k(z,\zeta ,\lambda )\wedge (\xi (\zeta -
z))^{-1}d_{\zeta }\xi _k(\zeta -z)].$ 
\par In accordance with \S \S 3.1 and 3.5 $\xi $ and $v$ are of class
of smoothness $C^{(q,n)}$, hence $\kappa $ and ${\tilde {\gamma }}$
belong to $C^{(q,n)}(W,L(\Lambda {\bf K}(\alpha ),Y))$
for suitable clopen $W\subset \Omega \times ({\bf K}\oplus \alpha
{\bf K})^m\times B({\bf K},0,1)$ such that $\Omega \times
(\partial \Omega )\times B({\bf K},0,1)\subset W$, $\zeta \ne z$.
Condition $3.5.(7)$ is satisfied for $\xi ({\bar {\zeta }}-{\bar z})=
Exp (\pi ^{-s}({\bar {\zeta }}-{\bar z}))$ and $v(z,\zeta )$ such that
$\xi (v(z,\zeta ))=Exp (\pi ^{-{\tilde {\phi }}(s)}v(z,\zeta ))$,
where ${\tilde {\phi }}(s)$ is given by \S 3.5 and satisfies
Condition $3.5.(1)$. Therefore, the family of such differential
forms $\psi $ and $w$ is nonvoid. In view of Conditions $3.1.(3)$
and $3.5.(7)$ in the sence of distributions:
\par $(5)$ $d_{z,\zeta ,\lambda }\kappa =0$ on $W$, that is,
$\mbox{ }_{z\in \Omega }P^n[(d_{z,\zeta ,\lambda }\kappa )\wedge \nu ]=0$
for each $\nu $ as above.
\par From $\partial _{\zeta }\kappa =0$ and $(5)$ it follows
\par $({\bar {\partial }}_{z,\zeta }+d_{\lambda }+\partial _z)(\kappa )=0$,
together with $\partial _{\zeta }({\tilde {\gamma }})=0$
it implies:
\par $(6)$ $({\bar {\partial }}_{z,\zeta }+d_{\lambda })({\tilde {\gamma }})
+\partial _z({\tilde {\gamma }})+({\bar {\partial }}_{z,\zeta }+d_{\lambda }
+\partial _z)(\kappa -{\tilde {\gamma }})=0$ on $W$, $\zeta \ne z$, \\
since $(\kappa -{\tilde {\gamma }})$ contains a factor $\partial _z\eta ^v_k$
and $\partial _z(\kappa -{\tilde {\gamma }})=0$. The monomials in
$(\kappa -{\tilde {\gamma }})$ with respect to $dz_j$, $d{\bar z}_j$,
$d\zeta _j$, $d{\bar {\zeta }}_j$ and $d\lambda $ and, consequently,
in $({\bar {\partial }}_{z,\zeta }+d_{\lambda }+\partial _z)(
{\tilde {\gamma }} - \kappa )$ contain at least one of the differentaials
$dz_1,...,dz_m$ as a factor. The same is true for $\partial _z({\tilde
{\gamma }})$. The monomials in $({\bar {\partial }}_{z,\zeta }+d_{\lambda })
({\tilde {\gamma }})$ do not contain any of the differentials $dz_j$. Hence
from $(6)$ it follows, that $({\bar {\partial }}_{\zeta }+d_{\lambda })
(\tilde {\gamma })=-{\bar {\partial }}_z{\tilde {\gamma }}.$ Then
\par $(7)$ $d_{\zeta ,\lambda }(f\wedge {\tilde {\gamma }})=
({\bar {\partial }}_{\zeta }+d_{\lambda })(f\wedge {\tilde {\gamma }})=
({\bar {\partial }}f)\wedge {\tilde {\gamma }}+
(-1)^tf\wedge ({\bar {\partial }}_{\zeta }+d_{\lambda }){\tilde
{\gamma }}=({\bar {\partial }}f)\wedge {\tilde {\gamma }}-
{\bar {\partial }}_z(f\wedge {\tilde {\gamma }})$. The applying
of Corollary $2.3.2$ and Formula $(7)$ to the differential form
$f\wedge {\tilde {\gamma }}$ on $(\partial \Omega )\times B$,
where $B:=B({\bf K},0,1)$, gives
\par $(8)$ $\mbox{ }_{(\zeta ,\lambda )\in (\partial \Omega )\times B}
P^n[({\bar {\partial }}f)\wedge {\tilde {\gamma }}]-
{\bar {\partial }}_z\mbox{ }_{(\zeta ,\lambda )\in (\partial \Omega )
\times B}P^n[f\wedge {\tilde {\gamma }}]$ \\
$=\mbox{ }_{(\partial \Omega )\times \{ 0 \} }
P^n[f\wedge {\tilde {\gamma }}]- \mbox{ }_{(\partial \Omega )
\times \{ \beta \} }P^n[f\wedge {\tilde {\gamma }}]$. \\
On the other hand, ${\tilde {\gamma }}|_{\lambda =0}=
\psi $, ${\tilde {\gamma }}|_{\lambda =\beta }=\tilde w$
and Formula $(8)$ is equivalent to Fromula $(3)$ due to
Formulas $3.8.(3,12,13)$.
\par {\bf 3.11. Corollary.} {\it Let $M$ and $f$ be as in
Theorem $3.10$ and $\partial v/\partial {\bar z}=0$ on $M$.
For $t=1,...,m$ put
\par $(1)$ $T^n_t:=(-1)^t(R^{v,n}_{\partial M}+B^n_M).$ Then
\par $(2)$ $f(z)={\bar {\partial }}(T^n_tf)(z)+(T^n_{t+1}
{\bar {\partial }}f)(z)$ \\
for each $z\in \tilde M$ and satisfying
Condition $3.10.(2)$. If ${\bar {\partial }}f=0$, then
$u=T^n_tf$ is a solution of ${\bar {\partial }}u(z)=f(z)$
for each $z\in \tilde M$ and $f$ satisfying $3.10.(2)$.}
\par {\bf Proof.} In view of Formula $3.8.(18)$ $L^{v,n}_{\partial M}
f=\mbox{ }_{\partial M}P^n[f\wedge \Upsilon ^v_t].$ Since
$\partial v(z,\zeta )/\partial {\bar z}=0$ the monomials in
$\Upsilon ^v_t$ of bidegree $(0,t)$ in $z$ vanish if $t\ge 1$.
Therefore, $L^{v,n}_{\partial M}f=0$ and $(2)$ follows from
$3.10.(1)$. Then from $(2)$ it follows ${\bar {\partial }}u(z)=
f(z)$ if ${\bar {\partial }}f(z)=0$ for each $z\in \tilde M$
and satisfying $3.10.(2)$, where $u=T^n_tf$.
\par {\bf 3.12. Definitions.} Let $M$ be a manifold over $\bf K$
satisfying $2.4.2$ with $(q+1,n)$-antiderivationally holomorphic
$\mbox{ }_SC^{(q+1,n-1)}$-transition maps $\phi _i\circ \phi _j^{-1}$
between charts $(U_i,\phi _i)$ and $(U_j,\phi _j)$ for each
$U_i\cap U_j\ne \emptyset $ and let $GL(N,{\bf K}(\alpha ))$
be the group of invertible $N\times N$-matrices with entries
in ${\bf K}(\alpha )$.
\par $(1)$. A $(q+1,n)$-antiderivationally holomorphic vector bundle
over ${\bf K}(\alpha )$ of ${\bf K}(\alpha )$ dimension $N$ over
$M$ is a $\mbox{ }_SC^{(q+1,n-1)}$-vector bundle over $M$ with
the characteristic fibre $({\bf K}(\alpha ))^N$ and with
$(q+1,n)$-antiderivationally holomorphic atlas of local trivializations
of $B$, that is, with a family $ \{ U_j,h_j \} $ such that
$\{ U_j \} $ is a (cl)open covering of $M$, for each $j$, $h_j$ is
a $\mbox{ }_SC^{(q+1,n-1)}$-bundle isomorphism from $B|_{U_j}$
onto $U_j\times ({\bf K}(\alpha ))^N$; the corresponding transition
mappings $g_{i,j}: U_i\cap U_j\to GL(N,{\bf K}(\alpha ))$ defined by
$(z,g_{i,j}(z)v)=h_i\circ h_j^{-1}(z,v)$, $z\in U_i\cap U_j$,
$v\in ({\bf K}(\alpha ))^N$ are $(q+1,n)$-antiderivationally holomorphic
$\mbox{ }_SC^{(q+1,n-1)}$-mappings. Equipped with the atlas
$ \{ B|_{U_j},h_j \} $ the bundle $B$ gets the structure of the
$\mbox{ }_SC^{(q+1,n-1)}-(q+1,n)$-antiderivationally holomorphic manifold.
\par $(2)$. A $\mbox{ }_SC^{(q+1,n-1)}$-bundle homomorphism  between
$\mbox{ }_SC^{(q+1,n-1)}-(q+1,n)$-antiderivationally holomorphic
vector bundles $B_1$ and $B_2$ is called
$\mbox{ }_SC^{(q+1,n-1)}-(q+1,n)$-antiderivationally
holomorphic if it is
$\mbox{ }_SC^{(q+1,n-1)}-(q+1,n)$-antiderivationally holomorphic
as a map between the
$\mbox{ }_SC^{(q+1,n-1)}-(q+1,n)$-antiderivationally holomorphic
manifolds $B_1$ and $B_2$. Similarly is defined a
$\mbox{ }_SC^{(q+1,n-1)}-(q+1,n)$-antiderivationally holomorphic
section of a ${\bf K}(\alpha )$
$\mbox{ }_SC^{(q+1,n-1)}-(q+1,n)$-antiderivationally holomorphic
vector bundle.
\par $(3)$. A $\mbox{ }_SC^{(q+1,n-1)}-(q+1,n)$-antiderivationally
holomorphic vector bundle over $M$ is called
$\mbox{ }_SC^{(q+1,n-1)}-(q+1,n)$-antiderivationally holomorphically
trivial if there exists a
$\mbox{ }_SC^{(q+1,n-1)}-(q+1,n)$-antiderivationally holomorphic
bundle isomorphism from $B$ onto $M\times ({\bf K}(\alpha ))^N$.
$B$ is called $\mbox{ }_SC^{(q+1,n-1)}-(q+1,n)$-antiderivationally
holomorphically trivial over a (cl)open set $U\subset M$ if $B|_U$ is
$\mbox{ }_SC^{(q+1,n-1)}-(q+1,n)$-antiderivationally holomorphically
trivial. A $\mbox{ }_SC^{(q+1,n-1)}-(q+1,n)$-antiderivationally holomorphic
trivialization of $B$ (over $U$) is a
$\mbox{ }_SC^{(q+1,n-1)}-(q+1,n)$-antiderivationally holomorphic
bundle isomorphism from $B$ onto $M\times ({\bf K}(\alpha ))^N$
($B|_U$ onto $U\times ({\bf K}(\alpha ))^N$).
\par $(4)$. A ${\bf K}(\alpha )$-valued differential form
of degree $r$ over $M$ can be defined as a section of the vector
bundle $\Lambda ^rT^*(M)_{{\bf K}(\alpha )}$, where
$T^*(X)_{{\bf K}(\alpha )}$ is the ${\bf K}(\alpha )$ cotangent
bundle of $M$ over scalars $b\in {\bf K}(\alpha )$
(see \cite{lustpnam}). A differential form of degree $r$ with values
in a $\mbox{ }_SC^{(q+1,n-1)}-(q+1,n)$-antiderivationally holomorphic
bundle (or a $B$-valued differential form) over $M$ is a section
of the bundle $\Lambda ^r(T^*(M)_{{\bf K}(\alpha )})
\otimes _{{\bf K}(\alpha )} B$.
\par If $ \{ U_j : j\in J \} $ is a (cl)open covering of $M$ such that
$B$ is $\mbox{ }_SC^{(q+1,n-1)}-$ $(q+1,n)$-antiderivationally holomorphically
trivial over each $U_j$ and $\{ g_{i,j}: i, j \in J \} $ is the
corresponding system of transition functions, then a differential form
with values in $M$ can be identified with a system $ \{ f_j \} $
of $N$-tuplets of differential forms on $U_j$ such that $f_i=
g_{i,j}f_j$ over $U_i\cap U_j$ for each $i, j \in J$.
A differential form $f$ with values in $B$ is called a $(0,t)$-form,
$\mbox{ }_SC^{(q+1,n-1)}-(0,t)$-form, etc. If for each (cl)open
subset $U\subset M$, where $B$ is
$\mbox{ }_SC^{(q+1,n-1)}-(q+1,n)$-antiderivationally holomorphically
trivial, the corresponding $N$-tuple of differential forms on $U$
consists of $(0,t)$-forms, $\mbox{ }_SC^{(q+1,n-1)}-(0,t)$-form,
etc. Each $(s,t)$-form with values in
a $\mbox{ }_SC^{(q+1,n-1)}$ $(q+1,n)$-antiderivationally holomorphic
vector bundle can be identified with some $(0,t)$-forms with values
in some other $n$-antiderivationally holomorphic vector bundle.
\par {\bf 3.13. Definition.} Let $M$ be a 
$\mbox{ }_SC^{(q+1,n-1)}-(q+1,n)$-antiderivationally holomorphic
manifold, let $B$ be a $\mbox{ }_SC^{(q+1,n-1)}-(q+1,n)$-antiderivationally
holomorphic
vector bundle over $M$ and $\{ U_j: j\in J \} $ be a (cl)open
covering of $M$, where $J$ is a set.
A derivationally $(q+1,n-1)$-holomorphic
Cousin data in $M$ means a system $ \{ f_{i,j}: i, j \in J \} $
of derivationally $(q+1,n-1)$-holomorphic
sections $f_{i,j}: U_i\cap U_j\to B$ such that $f_{i,j}+f_{j,k}=
f_{i,k}$ in $U_i\cap U_j\cap U_k$ for each $i, j, k\in J$.
The corresponding Cousin problem consists in finding a system
$\{ f_j: j\in J \} $ of derivationally $(q+1,n-1)$-holomorphic
sections $f_j: U_j\to B$ such that $f_{i,j}=f_i-f_j$ in
$U_i\cap U_j$ for each $i, j \in J$.
\par {\bf 3.14. Theorem.} {\it Let $M$ be a
$\mbox{ }_SC^{(q+1,n-1)}-(q+1,n)$-antiderivationally holomorphic
manifold and let $B$ be a
$\mbox{ }_SC^{(q+1,n-1)}-(q+1,n)$-antiderivationally holomorphic
vector bundle over $M$. Consider two conditions:
\par $(1)$ each derivationally $(q+2,n-1)$-holomorphic
Cousin problem in $B$ has a solution;
\par $(2)$ each $B$-valued
$\mbox{ }_SC^{(q+1,n-1)}-(0,1)$-form on $M$ such that
${\bar {\partial }}f=0$ on $M$ has a section $u: M\to B$
such that ${\bar {\partial }}u=f$ on $M$.
\par Then from $(1)$ it follows $(2)$. From $(2)$ it follows $(1)$
in the class $u\in C^{(q+2,n-1)}$
and ${\bar {\partial }}u\in \mbox{ }_SC^{(q+1,n-1)}$
is $(q+1,n-1)$-antiderivationally holomorphic.}
\par {\bf Proof.} $(1) \Longrightarrow (2)$. At first,
$f$ is a $\mbox{ }_SC^{(q+1,n-1)}$-form means, that there exists
a refinement $\{ {U'}_k: k \} $ of $ \{ U_j \} $ consisting
of clopen ${U'}_k$ such that $g_k({U'}_k)$ is bounded in
$({\bf K}(\alpha ))^N$ and $f|_{{U'}_k}\in
\mbox{ }_SC^{(q+1,n-1)}$, where $At'(M)= \{ ({U'}_k, g_k): k \} $.
Choose $At'(M)$ such that $\bigcup_k{{\tilde U}'}_k=M$.
Denote $ \{ {U'}_k : k \} $ by $ \{ U_j: j
\in J \} $ also such that $\partial {U'}_k$ satisfies
condition of Theorem $2.8$ up to the
$\mbox{ }_SC^{(q+1,n-1)}$-diffeomorphism. Then $2.8$ on each
$U_j$ gives a solution $u_j$ such that $(u_i-u_j)$ are
derivationally $(q+2,n-1)$-holomorphic
on $U_i\cap U_j$ and form derivationally $(q+2,n-1)$-holomorphic
Cousin data in $B$. According to $(1)$ there exists a
derivationally $(q+2,n-1)$-holomorphic
section $h_j: U_j\to B$ such that $u_i-u_j=h_i-h_j$ in
$U_i\cap U_j$. Set $u:=u_i-h_i$ in $U_i$ for each $j\in J$.
\par $(2) \Longrightarrow (1)$ in the class $u\in C^{(q+2,n-1)}$
and ${\bar {\partial }}u\in \mbox{ }_SC^{(q+1,n-1)}$
is $(q+1,n-1)$-antiderivationally holomorphic. Characteristic functions
of clopen compact subsets belong to $C^{\infty }$.
It is possible to take a refinement $At'(M)$ of $At(M)$
such that its charts be satisfying Lemma 2.6.1, that is,
$g_k({U'}_k)$ are balls satisfying $2.6.1$.
Choose $At'(M)$ such that $\bigcup_k{{\tilde U}'}_k=M$. Denote
it also by $At(M)$. Since $M$ is metrizable it has
an atlas consisting of clopen compact charts, hence
$M$ has a $C^{\infty }$-partition of unity,
$\chi _k:=\chi _{U_k}$. For each $i$ and $j$ $f_{k,j}$
is $\mbox{ }_SC^{(q+1,n-1)}-(q+1,n-1)$-antiderivationally holomorphic,
hence $\chi _kf_{k,j}$ is also by Lemma $2.6.1$
$\mbox{ }_SC^{(q+1,n-1)}-(q+1,n-1)$-antiderivationally holomorphic
for suitable refinement $\{ U_j: j\in J \} $, since
$\chi _kf_{k,j}=f_{k,j}|_{(U_k\cap dom (f_{k,j}))}$,
${\bar {\partial }}(\chi _kf_{k,j})=0$. Set $\theta _j:=
-\sum_k \chi _kf_{k,j}$ in $U_j$, hence by Theorem $2.7.2$ there exists
a $C^{(q+2,n-1)}$-solution of the Cousin problem:
$f_{i,j}:=\sum_k \chi _k(f_{i,k}+f_{k,j})=\theta _i-\theta _j$
in $U_i\cap U_j$; ${\bar {\partial }}\theta _i=
{\bar {\partial }}\theta _j$ in $U_i\cap U_j$. Hence by $(2)$ there exists
a section $u: M\to B$ such that ${\bar {\partial }}u=
{\bar {\partial }}\theta _j$ in $U_j$. The setting $h_j=\theta _j-u$
in $U_j$ provides $(1)$.
\par {\bf 3.15. Remark.} Formulas $3.6.(1), 3.9.(1),
3.10.(1)$ are the non-Archimedean analogs of the Leray,
Koppelman and Koppelman-Leray formulas correspondingly.
\par {\bf 3.16. Notes and Definitions.} The local field $\bf K$ is
the disjoint union of balls $B({\bf K},z_j,R)$ for a given $0<R<\infty $,
where $z_j\in \bf K$ for each $j\in \bf N$. Therefore, the antiderivation
operators $\mbox{ }_{B_j}P^n$ on $B_j:=B({\bf K},z_j,R)$ induce
the antiderivation operator $\mbox{ }_{\bf K}P^n$ on $\bf K$ such that
\par $(1)$ $\mbox{ }_{\bf K}(P^n[f])(y):=\sum_{j=1}^{\infty }
(\mbox{ }_{B_j}P^n[f\chi _{B_j}])(y)$ \\
on $C^{(q,n-1)}({\bf K},{\bf L})$, where ${\bf K}\subset
{\bf L}\subset \bf C_p$, $\bf L$ is a field complete relative
to its uniformity. Then
\par $\mbox{ }_PC^{(q,n)}({\bf K^l},Y):=\mbox{ }_{\bf K^l}P^n(
C^{(q,n-1)}({\bf K^l},Y))\oplus Y$  and
\par $\mbox{ }_SC^{(q+1,n-1)}({\bf K^l},Y):= \{ g\in C^{(q+1,n-1)}
({\bf K^l},Y):$ $g(x_1,...,x_l)$ \\
$\in \mbox{ }_{P,x_j}C^{(q+1,n-1)}
({\bf K^l},Y)$ $\mbox{for each}$ $j=1,...,l \} $,
\par $\mbox{ }_{P,x_j}C^{(q+1,n-1)}({\bf K^l},Y):=
\mbox{ }_{\bf K}P^n_{x_j}(C^{(q,n-1)}({\bf K^l},Y))\oplus Y$, \\
where $C^{(q,n-1)}({\bf K^l},Y)$ and $\mbox{ }_PC^{(q,n)}({\bf K^l},Y)$
are supplied with the inductive limit topologies induced by the embeddings \\
$C^{(q,n-1)}(B({\bf K^l},z,R'),Y)\hookrightarrow C^{(q,n-1)}
({\bf K^l},Y)$, $0<R'<\infty $, where $\mbox{ }_{\bf K^l}P^n:=
\mbox{ }_{\bf K}P^n_{x_1}...\mbox{ }_{\bf K}P^n_{x_l}$,
$x_1,...,x_l\in \bf K$, $Y$ is a Banach space over $\bf L$ such that
${\bf K}\subset \bf L$ (see also \cite{luseam3,nari}).
\par Therefore, in the standard way we get
the definition of a locally compact manifold $M$ over $\bf K$ of class
$\mbox{ }_PC^{(q,n)}$ or $\mbox{ }_SC^{(q+1,n-1)}$, that is, transition
mappings of charts $\phi _{i,j}\in \mbox{ }_PC^{(q,n)}$
or  $\phi _{i,j}\in \mbox{ }_SC^{(q+1,n-1)}$, where $V_j$ is clopen
in $M$, $\phi _j(V_j)$ is clopen in $\bf K^l$, $1\le l\in \bf N$,
$l=dim_{\bf K}M$ (see \S 2.2.5). Using charts and
$\mbox{ }_PC^{(q,n)}({\bf K^l},{\bf K^m})$ or
$\mbox{ }_SC^{(q+1,n-1)}({\bf K^l},{\bf K^m})$ we get the uniform space
$\mbox{ }_PC^{(q,n)}(M,N)$ or $\mbox{ }_SC^{(q+1,n-1)}(M,N)$
of all mappings $g: M\to N$ of class $\mbox{ }_PC^{(q,n)}$
or $\mbox{ }_SC^{(q+1,n-1)}$ respectively, where $M$ is the
$\mbox{ }_PC^{(q,n)}$-manifold or $\mbox{ }_SC^{(q+1,n-1)}$-manifold
on $\bf K^l$ and $N$ is the
$\mbox{ }_PC^{(q,n)}$-manifold or $\mbox{ }_SC^{(q+1,n-1)}$-manifold
on $\bf K^m$ correspondingly, that is,
$\psi _i\circ g\circ \phi _j^{-1}$ is of class
$\mbox{ }_PC^{(q,n)}$ or $\mbox{ }_SC^{(q+1,n-1)}$ for each $i$ and $j$
such that its domain is nonempty, where $At(M)= \{ (V_j,\phi _j):$
$j \} $, $At(N)= \{ (W_j,\psi _j):$ $j \} $. The uniformity
in $\mbox{ }_PC^{(q,n)}({\bf K^l},{\bf K^m})$ or
$\mbox{ }_SC^{(q+1,n-1)}({\bf K^l},{\bf K^m})$
induces the uniformity in
$\mbox{ }_PC^{(q,n)}(M,N)$ or $\mbox{ }_SC^{(q+1,n-1)}(M,N)$
respectively (see Remark 2.4 \cite{luseam3}).
\par For a locally compact manifold $M$ over $\bf K$
of class $\mbox{ }_PC^{(q,n)}$ or $\mbox{ }_SC^{(q+1,n-1)}$
let $DifP^{(q,n)}(M)$ or $DifS^{(q+1,n-1)}(M)$ denotes a family
of all diffeomorphisms $f: M\to M$, $f(M)=M$,
$(f-id)\in \mbox{ }_PC^{(q,n)}$ and $(f^{-1}-id)\in \mbox{ }_PC^{(q,n)}$ or
$(f-id)\in \mbox{ }_SC^{(q+1,n-1)}$ and $(f^{-1}-id)
\in \mbox{ }_SC^{(q+1,n-1)}$ respectively, where $id(z)=z$ for each
$z\in M$, $M\hookrightarrow {\bf K}^N$, $\mbox{ }_PC^{(q,n)}(M,M)
\hookrightarrow \mbox{ }_PC^{(q,n)}(M,{\bf K}^N)$,
$\mbox{ }_SC^{(q+1,n-1)}(M,M)\hookrightarrow \mbox{ }_SC^{(q+1,n-1)}
(M,{\bf K}^N)$ such that $(f-id)$ is correctly defined, $N\in \bf N$.
\par {\bf 3.17. Theorem.} {\it $(1).$ The uniform spaces $DifP^{(q,n)}(M)$
and $DifS^{(q+1,n-1)}(M)$ are the topological groups for each
$0\le q\in \bf Z$, $1\le n\in \bf N$.
\par $(2).$ They have embeddings as clopen subsets into
$\mbox{ }_PC^{(q,n)}(M,M)$ and into $\mbox{ }_SC^{(q+1,n-1)}(M,M)$
respectively.
\par $(3)$. The uniform spaces
$\mbox{ }_PC^{(q,n)}(M,N)$, $\mbox{ }_SC^{(q+1,n-1)}(M,N)$,
$DifP^{(q,n)}(M)$ and $DifS^{(q+1,n-1)}(M)$ are complete
and separable.
\par $(4)$. The groups $DifP^{(q,n)}(M)$
and $DifS^{(q+1,n-1)}(M)$ are ultrametrizable, when $M$ is compact.
\par $(5)$. The uniform spaces $DifP^{(q,n)}(M)$ and
$DifS^{(q+1,n-1)}(M)$ have the infinite-dimensional manifolds
structures over $\bf K$.}
\par {\bf Proof.} At first prove, that compositions of diffeomorphisms
preserve classes $\mbox{ }_PC^{(q,n)}(M,M)$ and
$\mbox{ }_SC^{(q+1,n-1)}(M,M)$ respectively.
For this consider two diffeomorphisms $\psi , \phi \in
DifP^{(q,n)}(U^m)$ or $DifS^{(q+1,n-1)}(U^m)$ simultaneously.
A diffeomorphism $\phi $ is called the simplest diffeomorphism,
if it has the coordinate form:
\par $x_j=\phi _j(y_1,...,y_m)=y_j$ for each $j=1,...,k-1, k+1,...,m$,
\par $x_k=\phi _k(y_1,...,y_m)=\phi _k(y_1,...,y_k,...,y_m)$,
where $x_j, y_j \in U$, $x=(x_1,...,x_m)$, $m=dim_{\bf K}M$.
Suppose such marked number $k$ is for $\phi $ and $l$ is for $\psi $.
To prove $\phi \circ \psi \in DifP^{(q,n)}(U^m)$ or
$DifS^{(q+1,n-1)}(U^m)$ it is sufficient to verify that
$\{ \phi _k(y_1,...,y_{k-1},\psi _l(y_1,...,y_m),y_{k+1},...,y_m)
-y_k \} $ is in \\
$\mbox{ }_PC^{(q,n)}(U^m,{\bf K})$ or $\mbox{ }_SC^{(q+1,n-1)}(U^m,{\bf K})$
correspondingly.
\par In $C^0(U^m,{\bf K}^m)$ there exists the polynomial Amice base
$\{ {\bar Q}_n(x):$ $n\in {\bf N_o}^m \} $ and it is also the base in
$C^{(q,n)}(U^m,{\bf K}^m)$, where ${\bf N_o}:= \{ j:$ $0\le j\in
{\bf Z} \} $ (see \cite{ami,luseam3}). The linear ordering
$\triangle $ in $\bf K$ induces the linear ordering $\triangle $
in ${\bf K}^m$ and hence in $U^m$: $x\triangle y$ if and only if
$x_1=y_1$,...,$x_{j-1}=y_{j-1}$, $x_j\triangle y_j$, where
$1\le j\le m$, $y_j\in \bf K$, $y=(y_1,...,y_m)$ (see \S 2.2.1).
Take in particular, $U=B({\bf K},0,1)$. Then $(\beta ,...,\beta )$
is the largest element in $U^m$. Let ${\bf Z_K}:= \{ z\in {\bf K}:$
$z=\sum_{l=0}^tz_l\pi _l,$ $0\le t\in {\bf Z}, z_l\in \{ 0,\theta _1,...,
\theta _{p^n-1} \} \} $, then $\bf Z_K$ is dense in $B({\bf K},0,1)$
and $\bf Z_K$ is countable. There are decompositions
\par $(i)$ $\psi _l(y)=\sum_{n\in {\bf N_o}^m } a(n,\psi _l){\bar Q}_n(y)$
and
\par $(ii)$ $\phi _k(y)=\sum_{n\in {\bf N_o}^m} a(n,\phi _k){\bar Q}_n(y)$,\\
where $a(n,\psi _l)$ and $a(n,\phi _k)\in \bf K$. In view of
the conditions imposed on $\psi _l$ and $\phi _k$ and continuity of
the $\bf K$-linear operators $\mbox{ }_UP^n_{x_j}$: \\
$(iii)$ $\phi _k(y)=\{ \sum_{n\in {\bf N_o}^m} a(n,\partial \phi _k(y)/
\partial y_j)(\mbox{ }_UP^n_{y_j}{\bar Q}_n(y)|^{y_j}_{y_{j,0}}) \} $ \\
$+\phi _k(y_1,...,y_{j-1},y_{j,0},y_{j+1},...,y_m)$ \\
for each $j=1,...,m$ and analogously for $\psi _l$, where $y_{j,0}$
and $y_j\in U$. To show $(\phi _k(y_1,...,y_{k-1},\psi _l(y),y_{k+1},
...,y_m) -y_k)\in \mbox{ }_SC^{(q+1,n-1)}(U^m,{\bf K})$ it is sufficient
to find $h_j: U^m\to \bf K$ such that
\par $(iv)$ $\mbox{ }_UP^n_{y_j}h_j|^{y_j}_{y_{j,0}} = -h_{j,0} +
\phi _k(y_1,...,y_{k-1},\psi _l(y),y_{k+1},...,y_m)-y_k$ \\
for each $j=1,...,m$, where $h_{j,0}\in \bf K$.
\par From $(iii)$ and continuity of the $\bf K$-linear operator
$\mbox{ }_UP^n_{y_j}$ it follows, that to resolve $(iv)$ it is
sufficient to find a solution of the problem:
\par $(v)$ $\mbox{ }_UP^n_{y_j}h|^{y_j}_{y_{j,0}}=
(\mbox{ }_UP^n_{y_j}y^{t^1}|^{y_j}_{y_{j,0}})...
(\mbox{ }_UP^n_{y_j}y^{t^l}|^{y_j}_{y_{j,0}})$ \\
for each $l\in \bf N$ and each
$t^k=(t^k_1,...,t^k_m)\in {\bf N_o}^m$, $k=1,...,l$,
$y^t=y_1^{t_1}...y_m^{t_m}$. On the other hand, \\
$(vi)$ $\mbox{ }_UP^nz^t=\sum_{0\le j\le n-1, k\in \bf N_o}
t(t-1)...(t-j+1)z_k^{t-j}(z_{k+1}-z_k)^{j+1}/(j+1)!$, \\
where $z\in U$, $t\in \bf N$, $j\in \bf Z$.
Moreover, $(\partial /\partial y_j)\mbox{ }_UP^n_{y_j}|_{(
C^{(q,n-1)}(U^m,{\bf K}))}=I$, hence
Equation $(v)$ can be simplified in the considered class
of $\mbox{ }_{P,y_j}C_0^{(q+1,n-1)}(U^m,{\bf K})$-functions
acting on both sides of $(v)$ by $(\partial /\partial y_j)$.
For each $z\in \bf Z_K$
there exists a solution $\mbox{ }_zh(y)$ of $(v)$ for each $y\in U^m$
such that $y_j\triangle z$, since the set $\{ u\in {\bf Z_K}:$ $u
\triangle z \} $ is finite. In view of $(vi)$ and \S 2.1 this
family $\{ \mbox{ }_zh(y):$ $z\in {\bf Z_K} \} $ can be chosen consistent,
that is, $\mbox{ }_zh(y)=\mbox{ }_{\eta }h(y)$ for each $y$ such that
$y_j\triangle \min (z,\eta )$. Therefore, there exists
\par $(vii)$ $h=\lim_{z\to \beta }\mbox{ }_zh$ \\
such that $(v)$ is satisfied for each $y\in U^m$. In particular,
$id\in \mbox{ }_SC^{(q+1,n-1)}(M,M)$.
\par For the class $\mbox{ }_PC^{(q,n)}(U^m,{\bf K})$ it is sufficient
to find solution of the problem
\par $(viii)$ $(\mbox{ }_{U^m}P^nh)(y)=(\mbox{ }_{U^m}P^ny^{t^1})...
(\mbox{ }_{U^m}P^ny^{t^l})$ \\
for each $l\in \bf N$ and each $t^k\in {\bf N_o}^m$, $k=1,...,l$,
$|t|:=t_1+...+t_m\ge 1$. In view of $(vi)$ and \S 2.1 and
$\mbox{ }_{U^m}P^n=\mbox{ }_UP^n_{y_1}...\mbox{ }_UP^n_{y_m}$
there exists a consistent family $\mbox{ }_zh$ satisfying
$(viii)$ for each $z\in {\bf Z_K}^m$ and each $y\triangle z$
such that $\mbox{ }_zh(y)=\mbox{ }_{\eta }h(y)$ for each
$y\triangle \min (z,\eta )$, where $\eta \in {\bf Z_K}^m$, since the set
$\{ u\in {\bf Z_K}^m:$ $u\triangle z \} $ is finite,
$(\partial /\partial y_1)...(\partial /\partial y_m)\mbox{ }_{U^m}P^n|_{
(C^{(q,n-1)}(U^m,{\bf K}))}=I$ and the acting by 
$(\partial /\partial y_1)...(\partial /\partial y_m)$ on both sides
of Equation $(viii)$ simplifies it in the class of
$\mbox{ }_PC_0^{(q,n)}(U^m,{\bf K})$-functions. Then
\par $(ix)$ $h=\lim_{z\to (\beta ,...,\beta )}\mbox{ }_zh$ is the solution
of $(viii)$. Therefore, $(\phi \circ \psi (y)-y)$ and
$(\phi \circ \psi ^{-1}(y)-y)$ belong to $\mbox{ }_SC^{(q+1,n-1)}$
or $\mbox{ }_PC^{(q,n)}$ correspondingly. The proof above also shows, that
if a bijective surjective $\psi $ is in $\mbox{ }_PC^{(q,n)}(M,M)$
or in $\mbox{ }_SC^{(q+1,n-1)}(M,M)$, then $\psi ^{-1}$ is in
$\mbox{ }_PC^{(q,n)}(M,M)$ or in $\mbox{ }_SC^{(q+1,n-1)}(M,M)$ respectively,
by solving the equation of the type $v(id(y)+g(y))=-g(y)$ relative to
the function $v$ for known $g:=\psi -id$.
Hence using charts $({\tilde V}_j,{\tilde {\phi }}_j)$ of
${\tilde A}t(M)$ such that ${\tilde {\phi }}_j({\tilde V}_j)=B
\subset U^m+z_j$ with suitable $z_j\in {\bf K}^m$ for each $j$ and
${\tilde A}t(M)$ is the refinement of $At(M)$ and $B$ satisfies Lemma
$2.6.1$ (or applying the above proof to $B$ instead of $U^m$), we get
that $({\tilde {\phi }}_l\circ \phi \circ \psi ^k\circ {\tilde {\phi }}^{-1}
(y)-y)$ belongs to $\mbox{ }_PC^{(q,n)}$ or $\mbox{ }_SC^{(q+1,n-1)}$
respectively on its domain for each $l$ and $j$, where $k=1$ or
$k=-1$. Together with Lemma $2.6.1$ it provides $\phi \circ \psi ^k\in
DifP^{(q,n)}(M)$ or $\phi \circ \psi ^k\in DifS^{(q+1,n-1)}(M)$
correspondingly for each $k\in \{ -1, 1 \} $.
\par If $M$ is compact, then $\mbox{ }_PC^{(q,n)}(M,Y)$ is normable
for a Banach space $Y$ over $\bf L$, ${\bf K}\subset \bf L$
(see analogously Lemma 3.4 \cite{lustpnam}). Let $V=B(C^{(q,n-1)}(M,Y),
0,1)$, consider $W:= \{ f\in C^{(q+1,n-1)}(M,Y):$ $
f(x_1,...,x_m)\in \mbox{ }_{P,x_j}C^{(q+1,n-1)}(M,Y)\cap (P^n_{x_j}(V)
\oplus Y)$ $\mbox{for each}$ $j=1,...,m \} .$ In view of $\bf K$-convexity
of $V$ the set $W$ is absolutely $\bf K$-convex (disked) and $W$ is
absorbing in $\mbox{ }_SC^{(q+1,n-1)}(M,Y)$, since
$P^n_{x_j}$ are continuous $\bf K$-linear and $V$ is absorbing
in $C^{(q,n-1)}(M,Y)$. Then $W$ is bounded in the weak topology in
$\mbox{ }_SC^{(q+1,n-1)}(M,Y)$. Therefore, the Minkowski
functional on $\mbox{ }_SC^{(q+1,n-1)}(M,Y)$ generated by $W$
induces a norm in $\mbox{ }_SC^{(q+1,n-1)}(M,Y)$ (see
Exer. 6.204 \cite{nari}). Each space $\mbox{ }_{P,x_j}C^{(q+1,n-1)}
(M,Y)$ is complete (see analogously with Lemma $3.4$ \cite{lustpnam}),
since $Y$ is complete.
\par Consider the $\bf K$-linear space $\Psi _j:=
\mbox{ }_{P,x_j}C^{(q+1,n-1)}(M,Y)\cap \mbox{ }_SC^{(q+1,n-1)}
(M,Y)$ and topologies $\tau _{P,j}$ on $\mbox{ }_{P,x_j}C^{(q+1,n-1)}
(M,Y)$ and $\tau _S$ on $\mbox{ }_SC^{(q+1,n-1)}(M,Y)$ induced by
norms in these spaces, then $\tau _S|_{\Psi _j}\subset \tau _{P,j}$
for each $j$ due to continuity of $P^n_{x_j}$ (for $M$ supplied with
coordinates $x_j$ due to $\mbox{ }_PC^{(q,n)}$ or
$\mbox{ }_SC^{(q+1,n-1)}$-diffeomorphism with $\Omega $ as in \S
2.2.5) and definition of $\tau _S$, since $ker (P^n_{x_j})= \{ 0 \} $
and due to the open mapping Theorem $(14.4.1)$ \cite{nari}
there exists the continuous $\bf K$-linear operator
\par $(P^n_{x_j})^{-1}: (\mbox{ }_{P,x_j}C_0^{(q+1,n-1)}(M,Y), \tau _{P,j})
\to  (C^{(q,n-1)}(M,Y), \| * \| _{C^{(q,n-1)}(M,Y)})$, consequently,
\par $(P^n_{x_j})^{-1}: (\Psi _{j,0}, \tau _S|_{\Psi _{j,0}})
\to (C^{(q,n-1)}(M,Y), \| * \| _{C^{(q,n-1)}(M,Y)})$, is continuous,
where $\Psi _{j,0}:=\Psi _j\cap \mbox{ }_{P,x_j}C_0^{(q+1,n-1)}(M,Y)$,
$\Psi _j=\Psi _{j,0}\oplus Y$. Hence $\mbox{ }_SC^{(q+1,n-1)}(M,Y)$
is complete relative to the above norm.
\par For noncompact $M$ using a refinement $At'(M)$ consisting of
compact charts $({V'}_j,{\phi '}_j)$ and the strict inductive limits
of $\mbox{ }_PC^{(q,n)}(\bigcup_{j=1}^l{V'}_j,Y)$
or $\mbox{ }_SC^{(q+1,n-1)}(\bigcup_{j=1}^l{V'}_j,Y)$, $l\in \bf N$,
we get, that $\mbox{ }_PC^{(q,n)}(M,Y)$ and
$\mbox{ }_SC^{(q+1,n-1)}(M,Y)$ are complete relative to their uniformities
(see Theorems $(12.1.6)$ and $(12.1.8)$ \cite{nari}). In view of
Theorem $(12.1.4)$ \cite{nari} these spaces are separable.
\par Let $(f-id)\in \mbox{ }_PC^{(q,n)}(M,M)$ or $(f-id)\in
\mbox{ }_SC^{(q+1,n-1)}(M,M)$ such that $M$ is compact and
$\max_{j,l} \| f_{l,j}-id_{l,j} \| <1$, where $f_{l,j}:=
\phi _l\circ f\circ \phi _j^{-1}$, $dom (f_{l,j})=:U_{l,j}$,
$\| * \| $ is taken of the space $\mbox{ }_PC^{(q,n)}(U_{l,j},{\bf K}^m)$
or $\mbox{ }_SC^{(q+1,n-1)}(U_{l,j},{\bf K}^m)$. In view of the ultrametric
inequality $f_{l,j}$ is the isometry, since
\par $ \| f_{l,j}-id_{l,j} \| =
\sup_n |a(n,f_{l,j}-id_{l,j})| \| {\bar Q}_n \| $,
where $\| * \| $ is the norm in $\mbox{ }_PC^{(q,n)}(U_{l,j},{\bf K}^m)$
or in $\mbox{ }_SC^{(q+1,n-1)}(U_{l,j},{\bf K}^m)$ respectively
induced by the norm in $C^{(q,n-1)}(U_{l,j},{\bf K}^m)$ and the
Minkowski functional as above.
Then $ \| g_{k,l}\circ f_{l,j}-id_{k,j} \| \le \max ( \| g_{k,l}
\circ f_{l,j}-f_{l,j} \| , \| f_{l,j}-id_{l,j} \| )$.
Using partial difference quotients and $P^n$ and expansion
coefficients in the Amice base we get, that
\par $\max_{l,j} \| f_{l,j}^{-1}-id_{l,j} \| \le C \max_{l,j}
\| f_{l,j}-id_{l,j} \| $, $C=const >0$ is independent of $f$
(see the proof of Theorem $2.6$ \cite{luseam3}), consequently,
$DifP^{(q,n)}(M)$ and $DifS^{(q+1,n-1)}(M)$ are topological groups.
For noncompact $M$ having $At (M)$ with compact charts and using
the strict inductive limit topology we can take an entourage of
the diagonal in $\mbox{ }_PC^{(q,n)}(M,M)^2$ or in
$\mbox{ }_SC^{(q+1,n-1)}(M,M)^2$ of the form
$\{ f: \| f_{l,j}-id_{l,j} \| \le |\pi |$ $\mbox{for each}$
$l, j\in \lambda \} $, where $\lambda $ is a finite subset in $\bf N$.
In view of Theorem $A.4$ \cite{lujmsq2} there exists the inverse
mapping $f_{l,j}^{-1}$, which is the local diffeomorphism, when
$dom (f_{l,j})\ne \emptyset $. Then $f|_W=id|_W$ for
$W:=M\setminus \bigcup_{j\in \lambda } {V'}_j$
for each $f\in DifP^{(q,n)}(M)$ and $DifS^{(q+1,n-1)}(M)$
with $W$ dependent on $f$, where
$supp (f):= cl ( \{ x\in M:$ $f(x)\ne x \} )$ is compact,
a finite subset $\lambda $ of $\bf N$ is such that
$supp (f) \subset \bigcup_{j\in \lambda } {V'}_j$. 
This implies that $f(M)=M$ and $f^{-1}(M)=M$, consequently,
$DifP^{(q,n)}(M)$ and $DifS^{(q+1,n-1)}(M)$ are neighborhoods
of $id$ in $\mbox{ }_PC^{(q,n)}(M,M)$ and in
$\mbox{ }_SC^{(q+1,n-1)}(M,M)$ respectively, left shifts in
these groups $L_gf:=g^{-1}f$ imply that these groups are open
in the corresponding to them spaces.
Since $DifP^{(q,n)}(M)$ and $DifS^{(q+1,n-1)}(M)$ are complete,
then they are clopen in $\mbox{ }_PC^{(q,n)}(M,M)$ and
$\mbox{ }_SC^{(q+1,n-1)}(M,M)$ correspondingly (see Theorem $8.3.6$
\cite{eng}).
\par Finally, statements $(4,5)$ follow from the proofs of Theorems
$2.4$ and $3.6$ \cite{luseam3} modified for the considered here
classes of smoothness.
\par {\bf 3.18. Remark and Definition.} Let $M$ and $N$ be two
locally compact $C^{(q,n)}$-manifolds over $\bf K$ and
$f\in C^1(M,N)$, $dim_{\bf K}M=:m_M$, $dim_{\bf K}N=:m_N$.
Denote by ${\cal E}:={\cal E}(f):= \{ z\in M:$ $rang (d_zf)<m_N \} $
and this set is called the set of critical values of $f$. The
nonnegative Haar measure $\nu $ on ${\bf K}^{m_N}$ as the additive group
induces the measure $\mu $ on $N$ with the help of charts,
since $At(N)$ has a disjoint refinement, where $\nu $ is normalized
by the condition $\nu (B({\bf K}^{m_N},0,1))=1$.
\par {\bf 3.19. Theorem.} {\it Let $f: M\to N$ be a $C^l$-mapping
of a $\mbox{ }_SC^{(q+1,n-1)}$-manifold $M$ into a
$\mbox{ }_SC^{(q+1,n-1)}$-manifold $N$, where $l>\max (m_M,m_N)$.
Then $\mu (f({\cal E})=0$ (see \S 3.18).}
\par {\bf Proof.} Using the charts of atlases it is sufficient to prove
the theorem for $f: U\to {\bf K}^{m_N}$, where $U$ is an open
subset in ${\bf K}^{m_N}$. For $m_M=0$ and $m_N=0$ the statement
is evident, therefore, consider $m_M\ge 1$ and  $m_N\ge 1$.
Put ${\cal E}_i:= \{ y\in U:$ $f^{(j)}(y)=0$ $\mbox{for each}$
$j\le i \} $, hence ${\cal E}\supset {\cal E}_1\supset {\cal E}_2
\supset ... $. To finish the proof use the following two lemmas.
\par {\bf 3.19.1. Lemma.} {\it $\mu (f({\cal E}\setminus {\cal E}_1))=0$.}
\par {\bf Proof.} Consider $n\ge 2$, since for $n=1$ there is only
one partial derivative and from $y\in \cal E$ it follows
$y\in {\cal E}_1$. Let $y\in {\cal E}\setminus {\cal E}_1$, then
there exists a nonzero partial derivative, for example,
$\partial f_1(x)/\partial x_1$ at the point $x=y$. There exists
a mapping $h: U\to {\bf K}^{m_N}$ such that $h(x):=(f_1(x),x_2,...,x_{m_N})$
for which $rang (dh(y))=m_N$. In view of Theorem $A.4$ \cite{lujmsq2}
the mapping $h$ is the diffeomorphism of some open $V=V(y)\subset U$
onto a neighborhood $W\ni z:=h(y)$. The set ${\cal E}'$ of critical
points for $g:=f\circ h^{-1}: W\to {\bf K}^{m_N}$ coinsides with
$h(V\cap {\cal  E})$, that is, $g({\cal E}')=f(V\cap {\cal E})$.
Consider the family $g^t: ( \{ t \} \times {\bf K}^{m_M-1})\cap W\to
\{ t \} \times {\bf K}^{m_N-1}$, where $t \in B({\bf K},0,1)$.
The point $b$ is critical for $g^t$ if and only if it is critical for
$g$. In view of the induction hypothesis $\mu [g^t({\cal E}(g^t))]=0$ in
$\{ t \} \times {\bf K}^{m_N-1}$, hence
$\mu (g({\cal E}')\cap ( \{ t \} \times {\bf K}^{m_N-1}))=0$
for each $t\in B({\bf K},0,1)$. From the Fubini theorem
in $L^1({\bf K^m},\mu ,{\bf R})$ it follows, that $\mu (g({\cal E}'))=0$.
\par {\bf 3.19.2. Lemma.} {\it $\mu (f({\cal E}_k))=0$ for each
$k$ such that $1\le k<l$.}
\par {\bf Proof.} Take a covering of ${\cal E}_k$ by a countable number
of balls of radius $\delta >0$, $\delta \le \delta _0$, where
$\delta _0>0$ is sufficiently small. Take one of these balls
$B$. From the definition of ${\cal E}_k$ and the Taylor formula
(see Theorem $29.4$ \cite{sch1} and  Theorem $A.5$ \cite{lujmsq2})
it follows, that $f(x+h)=f(x)+R(x,h)$, where $\| R(x,h) \|
\le b \| h \| ^{k+1}$, $x\in {\cal E}_k$, $x+h\in B$,
$b\le \| f \|_{C^l(U,{\bf K}^{m_N})}<\infty $ for a compact clopen $U$
in ${\bf K}^{m_M}$. Divide $B$ into a disjoint union of $q^{m_M}$ balls
of radius $\delta /q$, $q=p^{-n}$. Let $B_1$ be a ball
of this partition such that $B_1\ni x$. Then each $y\in B_1$
has the form $y=x+h$, where $|h|\le \delta /q$. Then $f(B_1)\subset
B({\bf K}^{m_N},f(x),b/q^{k+1})$, consequently,
$f({\cal E}_k\cap B)$ is contained in the union of $q^{m_N}$ balls
$B_j$ having $\sum_j\mu (B_j)\le q^{m_N}(b/q^{k+1})^{m_N}=
b^{m_N}q^{-m_Nk}$. Then $\lim_{q\to \infty }b^{m_N}q^{-m_Nk}=0$.
\par Therefore, Lemmas $3.19.1, 2$ finish the proof of Theorem $3.19.$
\par {\bf 3.19.3. Corollary.} {\it The set $N\setminus f({\cal E})$
is dense in $N$, where $f\in C^l(M,N)$ and $l>\max (m_M,m_N)$.}
\par {\bf 3.19.4. Corollary.} {\it If $dim_{\bf K}M<dim_{\bf K}N$,
then $\mu (f(M))=0$.}
\par {\bf 3.20. Definitions.} A $C^1$-mapping $f: M\to N$ is called
an immersion, if $rang (df|_x: T_xM\to T_{f(x)}N)=m_M$ for each $x\in M$.
An immersion $f: M\to N$ is called an embedding, if $f$ is bijective.
\par {\bf 3.21. Theorem.} {\it Let $M$ be a compact
$\mbox{ }_SC^{(q+1,n-1)}$ or $\mbox{ }_PC^{(q,n)}$-manifold over a local
field $\bf K$,
$dim_{\bf K}M=m<\infty $. Then there exists a $\mbox{ }_SC^{(q+1,n-1)}$
or $\mbox{ }_PC^{(q,n)}$-embedding $\tau : M\hookrightarrow {\bf K}^{2m+1}$
and a $\mbox{ }_SC^{(q+1,n-1)}$ or $\mbox{ }_PC^{(q,n)}$-immersion
$\theta : M\to {\bf K}^{2m}$ correspondingly. Each continuous mapping
$f: M\to {\bf K}^{2m+1}$ or $f: M\to {\bf K}^{2m}$ can be approximated
by $\tau $ or $\theta $ relative to the norm $\| * \|_{C^0}$.}
\par {\bf Proof.} Let $M\hookrightarrow {\bf K}^N$ be the
$\mbox{ }_SC^{(q+1,n-1)}$ or $\mbox{ }_PC^{(q,n)}$-embedding
of Theorem $2.2.6$. Consider the bundle of all $\bf K$ straight
lines in ${\bf K}^N$. They compose the projective space
${\bf K}P^{N-1}$. Fix the standard orthonormal (in the non-Archimedean
sence) base $ \{ e_1,...,e_N \} $ in ${\bf K}^N$ and projections
on $\bf K$-linear subspaces relative to this base
$P^L(x):=\sum_{e_j\in L}x_je_j$ for the $\bf K$-linear span
$L=span_{\bf K} \{ e_i:$ $i\in \Lambda _L \} $,
$\Lambda _L\subset \{ 1,...,N \} $, where
$x=\sum_{j=1}^Nx_je_j$, $x_j\in \bf K$ for each $j$.
In this base consider the function $(x,y):=\sum_{j=1}^Nx_jy_j$.
Let $l\in {\bf K}P^{N-1}$, take a $\bf K$-hyperplane
denoted by ${\bf K}^{N-1}_l$ and given by the condition:
$(x,[l])=0$ for each $x\in {\bf K}^{N-1}_l$, where $0\ne [l]\in {\bf K}^N$
characterises $l$. Take $\| [l] \| =1$. Then the orthonormal base
$\{ q_1,...,q_{N-1} \} $ in ${\bf K}^{N-1}_l$ and together with $[l]=:q_N$
composes the orthonormal base $\{ q_1,...,q_N \} $ in ${\bf K}^N$
(see also \cite{roo}).
This provides the projection $\pi _l: {\bf K}^N\to {\bf K}^{N-1}_l$
relative to the orthonormal base $ \{ q_1,...,q_N \} $.
The operator $\pi _l$ is $\bf K$-linear, hence $\pi _l\in
\mbox{ }_SC^{(q+1,n-1)}$, since $P^n$
is the $\bf K$-linear operator, $\mbox{ }_UP^n_{x_j}\lambda e_j|^b_a=
\lambda (b-a)e_j$ for each $\lambda \in \bf K$ and $a, b \in U$,
$j=1,...,N$.
\par To construct an immersion it is sufficient, that each
projection $\pi _l: T_xM\to {\bf K}^{N-1}_l$ has $ker [d(\pi _l(x))]=
\{ 0 \} $
for each $x\in M$. The set of all $x\in M$ for which $ker [d(\pi _l(x))]
\ne \{ 0 \} $
is called the set of forbidden directions of the first kind.
Forbidden are those and only those directions $l\in {\bf K}P^{N-1}$
for which there exists $x\in M$ such that $l'\subset T_xM$,
where $l'=[l]+z$, $z\in {\bf K}^N$. The set of all forbidden directions
of the first kind forms the $C^{(q,n-1)}$-manifold $Q$
of dimension $(2m-1)$ with points $(x,l)$, $x\in M$,
$l\in {\bf K}P^{N-1}$, $[l]\in T_xM$, where
$C^{(q,n)}\subset C^{(q+1,n-1)}$ for each $n\ge 1$, $q\ge 0$.
Take $g: Q\to {\bf K}P^{N-1}$
given by $g(x,l):=l$. Then $g$ is of class $C^{(q,n-1)}$.
In view of Theorem $3.19$
$\mu (g(Q))=0$, if $N-1>2m-1$, that is, $2m<N$. In particular,
$g(Q)$ is not contained in ${\bf K}P^{N-1}$ and there exists
$l_0\notin g(Q)$, consequently, there exists $\pi _{l_0}: M\to
{\bf K}^{N-1}_{l_0}$. Since $\mbox{ }_SC^{(q+1,n-1)}$
or $\mbox{ }_PC^{(q,n)}$ respectively is dense in $C^{(q,n-1)}$,
then there exists a mapping $\kappa $ such that
$\kappa \in \mbox{ }_SC^{(q+1,n-1)}$ or
$\kappa \in \mbox{ }_PC^{(q,n)}$ is sufficiently close to
$\pi _{l_0}$ relative to $\| * \|_{C^1}$ correspondingly
such that $\kappa \circ \theta $ is the immersion,
since $M$ is compact. In view of Theorem $3.17$
the composition $\kappa \circ \theta $ is of class
$\mbox{ }_SC^{(q+1,n-1)}$ or $\mbox{ }_PC^{(q,n)}$ correspondingly.
This procedure can be prolonged, when $2m<N-k$, where $k$
is the number of the step of projection. Hence $M$ can be immersed
in ${\bf K}^{2m}$.
\par Consider now the forbidden directions of the second type:
$l\in {\bf K}P^{N-1}$, for which there exist $x\ne y\in M$
simultaneously belonging to $l$ after suitable parrallel translation
$[l]\mapsto [l]+z$, $z\in {\bf K}^N$. The set of the forbidden directions
of the second type forms the manifold $S:=M^2\setminus \Delta $,
where $\Delta := \{ (x,x):$ $x\in M \} $. Consider $\psi :
S\to {\bf K}P^{N-1}$, where $\psi (x,y)$ is the straight $\bf K$-line
with the direction vector $[x,y]$ in the orthonormal base.
Then $\mu (\psi (S))=0$
in ${\bf K}P^{N-1}$, if $2m+1<N$. Then the closure
$cl (\psi (P))$ coinsides with $\psi (P)\cup g(Q)$ in ${\bf K}P^{N-1}$.
Hence there exists $l_0\notin cl (\psi (P))$. Then consider
$\pi _{l_0}: M\to {\bf K}_{l_0}^{N-1}$.
Since $\mbox{ }_SC^{(q+1,n-1)}$
or $\mbox{ }_PC^{(q,n)}$ correspondingly is dense in $C^{(q,n-1)}$,
then there exists a mapping $\kappa $ such that
$\kappa \in \mbox{ }_SC^{(q+1,n-1)}$ or
$\kappa \in \mbox{ }_PC^{(q,n)}$ is sufficiently close to
$\pi _{l_0}$ relative to $\| * \|_{C^1}$
such that $\kappa \circ \tau $ is the embedding,
since $M$ is compact. In view of Theorem $3.17$
the composition $\kappa \circ \tau $ is of class
$\mbox{ }_SC^{(q+1,n-1)}$ or $\mbox{ }_PC^{(q,n)}$ correspondingly.
This procedure can be prolonged, when $2m+1<N-k$, where $k$
is the number of the step of projection. Hence $M$ can be embedded
into ${\bf K}^{2m+1}$.
\par {\bf 3.21.2. Remark.} Theorems $3.19$ and $3.21$
are non-Archimedean analogs of the Sard's and Witney's theorems.
In Theorem $3.21$ classes of smoothness globally on $M$ are
important. Theorem $3.21$ justifies the considered class
of manifolds $M$ in the theorems above about antiderivational
representations of functions.
\par {\bf 3.22. Note and Definition.} The proof of Theorem $3.17$
shows, that the family of all diffeomorphisms of $M$ of the class
$\mbox{ }_PC((t,s))$ as it was defined slightly different in
\cite{lustpnam} also form the topological group.
Moreover, spaces $\mbox{ }_PC((t,s),\Omega \to Y):=
P(l,s)[C((t,s-1),\Omega \to Y)] \oplus Y$ and $\mbox{ }_PC^{(t,s)}
(M,Y)$ are topologically $\bf K$-linearly isomorphic, where
$l=[t]+1$, $[t]$ is the integer part of $t$, $[t]\le t$, $0\le t
\in \bf R$, 
though the antiderivation operators $P(l,s)$ on a clopen subset
$X'=\Omega $ in $B({\bf K}^m,0,1)$ (see \S 2.11 \cite{luambp00}) and
$\mbox{ }_{\Omega }P^s$ above (see \S \S 2.1, 2.2.5) are different.
\par Define by induction spaces $\mbox{ }_S^lC^{\xi +(l,0)}(\Omega ,Y)
:= \{ f\in C^{\xi +(l,0)}(\Omega ,Y):$ $f(x_1,...,x_m)\in
\mbox{ }_UP^{n+l}_{x_j}(\mbox{ }_S^{l-1}C^{\xi +(l-1,0)}(\Omega ,Y))
\oplus Y$ $\mbox{for each}$ $j=1,...,m \} $, where $l\in \bf N$,
$\mbox{ }_S^1C^{\xi +(1,0)}(\Omega ,Y):=\mbox{ }_SC^{\xi +(1,0)}
(\Omega ,Y)$, $\mbox{ }_S^0C^{\xi }(\Omega ,Y):=C^{\xi }(\Omega ,Y)$.
\par {\bf 3.23. Theorem.} {\it Let $M$ be a
$\mbox{ }_S^lC^{(q+l,n-1)}$-manifold over $\bf K$ with $l\ge 2$, then
there exists a clopen neighborhood ${\tilde T}M$ of $M$ in $TM$
and an exponential $\mbox{ }_S^lC^{(q+l,n-1)}$-mapping
$\exp : {\tilde T}M\to M$ of ${\tilde T}M$ on $M$.}
\par {\bf Proof.} As in the proof of Theorem $3.7$ \cite{lustpnam}
it can be shown that the non-Archimedean geodesic equation
$\nabla _{\dot c}{\dot c}=0$ with initial conditions
$c(0)=x_0$, ${\dot c}(0)=y_0$, $x_0\in M$, $y_0\in T_{x_0}M$ has
a unique $\mbox{ }_S^lC^{(q+l,n-1)}$-solution,
$c: B({\bf K},0,1)\to M$. For a chart $(U_j,\phi _j)$
containing $x$, put $\psi _j(b)=\phi _j\circ c(b)$, then
\par $\psi _j(b)=\phi _j(x_0)+\mbox{ }_UP^{q+l+n}(y_0+
\mbox{ }_UP^{q+l+n-1}f)$, \\
$\psi _j(b)=\psi _j(b;x_0,y_0)$, $b\in B({\bf K},0,1)$,
where $f\in \mbox{ }_S^{l-2}C^{(q+l+n-3)}(B({\bf K},0,1),{\bf K}^m)$,
consequently, the mapping $V_1\times B({\bf K}^m,0,\delta )
\ni ({\tilde x}_0,y_0)\mapsto \psi _j(\beta ;x_0,y_0)$
is of class of smoothness $\mbox{ }_S^lC^{(q+l,n-1)}$, where
$0<\delta $, ${\tilde x}_0=\phi _j(x_0)\in V_1\subset V_2\subset
\phi _j(U_j)$, $V_1$ and $V_2$ are clopen, $\delta $ and $V_1$
are sufficiently small, that to satisfy the inclusion
$\psi _j(\beta ;x_0,y_0)\in V_2$ for each $({\tilde x}_0,y_0)\in
V_1\times B({\bf K}^m,0,\delta )$. The rest of the proof see
in \S 3.7 \cite{lustpnam}.
\par {\bf 3.24. Theorem.} {\it Let $\Omega =\Omega _1\times
...\times \Omega _m$ be a polydisk in $({\bf K}\oplus \alpha {\bf K})^m$
and let \\
$\mbox{ }_S{\bar C}^{(q+1,n-1)}(\Omega ,{\bf K}(\alpha )):=
\{ f\in \mbox{ }_SC^{(q+1,n-1)}(\Omega ,{\bf K}(\alpha )):$
${\bar {\partial }}f=0$ $\mbox{on}$ $\Omega \} $, \\
then
$\mbox{ }_S{\bar C}^{(q+1,n-1)}(\Omega ,{\bf K}(\alpha ))
|_{\tilde {\Omega }}$ is the algebra over $\bf K$, where
$\tilde {\Omega }:=\{ z\in \Omega :$ $z_j$ $\mbox{is encompassed by}$
$\partial \Omega _j \} $ for each $j=1,...,m$, $z=(z_1,...,z_m) \} $.}
\par {\bf Proof.} Evidently
$\mbox{ }_S{\bar C}^{(q+1,n-1)}(\Omega ,{\bf K}(\alpha ))$
is the $\bf K$-linear space, since ${\bar {\partial }}(\lambda f)=
\lambda {\bar {\partial }}f$ for each $\lambda \in \bf K$
and ${\bar {\partial }}(f+g)={\bar {\partial }}f+
{\bar {\partial }}g$ for each $f, g\in
\mbox{ }_SC^{(q+1,n-1)}(\Omega ,{\bf K}(\alpha ))$.
It remains to verify that $fg|_{\tilde {\Omega }} \in \mbox{ }_S{\bar C}^{(
q+1,n-1)}(\Omega ,{\bf K}(\alpha ))|_{\tilde {\Omega }}$ for each
$f$ and $g\in \mbox{ }_S{\bar C}^{(q+1,n-1)}(\Omega ,{\bf K}(\alpha ))$,
where as in \S 2.6.1 $\mbox{ }_S{\bar C}^{(q+1,n-1)}(\Omega ,{\bf K}
(\alpha ))|_{\tilde {\Omega }}= \{ h|_{\tilde {\Omega }}:$
$h\in \mbox{ }_S{\bar C}^{(q+1,n-1)}(\Omega ,{\bf K}(\alpha )) \} $.
In view of Theorem $2.7.2.(i)$ if $f$ and $g\in
\mbox{ }_S{\bar C}^{(q+1,n-1)}(\Omega ,{\bf K}(\alpha ))$, then $f$ and
$g$ are locally $z$-analytic on ${\tilde {\Omega }}$, consequently,
$fg$ is locally $z$-analytic on ${\tilde {\Omega }}$. In view of
Formula $2.7.8.(2)$ or by the direct computation:
\par $(i)\quad res_{\xi }(z-\xi )^j=0$ for each $-1\ne j\in \bf Z$ and
each $\xi \in {\tilde {\Omega }}$, \\
since $res_{\xi }h=0$ for each
$h$ having a decomposition $2.7.8.(1)$ with $a_{-1}=a_{-1}(h)=0$,
indeed it is true for the particular $h(\beta )=h(0)$ for a loop $\gamma $
encompassing $0$ and
such that $h(x):=Exp [j Log (\gamma (x))]$ and $j\in \bf Z$,
$x\in B({\bf K},0,1)$, that leads to the general case.
\par On the other hand, $(fg)'(z)$ also is locally $z$-analytic
on ${\tilde {\Omega }}$. Therefore,
\par $\mbox{ }_{\gamma _j}P^n[(fg)'(z)dz_j]=0$ and particularly
\par $(ii)\quad
\mbox{ }_{\gamma _j}P^n[(\partial (f(z)g(z))/\partial z_j)dz_j]=0$ \\
for each loop $\gamma _j$ in $\Omega _j$ encompassed by $\partial
\Omega _j$ (see Theorem $2.5.3$), where $z=(z_1,...,z_m)$,
$\Omega =\Omega _1\times ...\times \Omega _m$, $\Omega _j$ is a
ball in ${\bf K}\oplus \alpha {\bf K}$, $z_j\in {\bf K}\oplus
\alpha {\bf K}$ for each $j=1,...,m$. Then
\par $\partial (\mbox{ }_{\gamma _j}P^n[(\partial (f(z)g(z))/\partial z_j)
dz_j])/ \partial z_j=\partial (fg)(z)/\partial z_j$ and
\par $\mbox{ }_{\gamma _j}P^n[(\partial (fg)(z)/\partial z_j)dz_j]=
(fg)(z_1,...,z_{j-1},\gamma _j(\beta ),z_{j+1},...,z_m)$  \\
$- (fg)(z_1,...,z_{j-1},\gamma _j(0),z_{j+1},...,z_m)$. Moreover,
\par $\mbox{ }_{\gamma _j}P^n[h_j(z)dz_j]=\mbox{ }_BP^n[h_j(z_1,...,
z_{j-1},\gamma _j(\zeta ),z_{j+1},...,z_m)d\gamma _j(\zeta )]$,
\par $\mbox{ }_{\gamma _j}P^n[(\partial (fg)(z)/\partial z_j)dz_j]=
\mbox{ }_BP^n[v_j(z_1,...,z_{j-1},\gamma _j(\zeta ),
z_{j+1},...,z_m)d\gamma _j(\zeta )]$,
\par $\partial (\mbox{ }_BP^n[h_j(z_1,...,z_{j-1},\gamma _j(\zeta ),z_{j+1},
...,z_m)d\gamma _j(\zeta )]/\partial \zeta $ \\
$=h_j(z_1,...,z_{j-1},\gamma _j(\zeta ),z_{j+1},...,z_m)\gamma '(\zeta )$,
\par $\partial (\mbox{ }_BP^n[v_j(z_1,...,z_{j-1},\gamma _j(\zeta ),z_{j+1},
...,z_m)d\gamma _j(\zeta )]/\partial \zeta $ \\
$=v_j(z_1,...,z_{j-1},\gamma _j(\zeta ),z_{j+1},...,z_m)\gamma '(\zeta )$,
$\gamma _j\in \mbox{ }_PC^n(B,{\bf K}(\alpha ))$, where
$v_j(z):=\partial (fg)(z)/\partial z_j$,
$\zeta \in B:=B({\bf K},0,1)$.
Proceeding as in the proof of Theorem $3.17$ with the help
of $(i,ii)$ and Equation $2.2.5.(1)$ (see Equations $3.17.(iii-vii)$)
find $h_j\in
C^{(q,n-1)}(\Omega ,{\bf K}(\alpha ))$ such that $h_j$ is locally
$z$-analytic and $\mbox{ }_{\gamma _j}P^n[h_j(z_1,...,z_{j-1},\zeta _j,
z_{j+1},...,z_m)d\zeta _j]=(fg)(z)-(fg)(z_0)$ for each
$z\in {\tilde {\Omega }}$ and each $j=1,...,m$, where $\gamma _j$
is a path with $\gamma _j(0)=z_{j,0}\in {\tilde {\Omega }}_j$,
$\gamma _j(\beta )=z_j$ for each $j=1,...,m$. This means, that
$(fg)\in \mbox{ }_S{\bar C}^{(q+1,n-1)}(\Omega ,{\bf K}(\alpha ))|_{
\tilde {\Omega }}$, since \\
$\partial \mbox{ }_{\gamma _j}P^n[
h_j(z_1,...,z_{j-1},\zeta _j,z_{j+1},...,z_m)d\zeta _j]/\partial x_j$ \\
$=\partial \mbox{ }_{\gamma _j}P^n[h_j(z_1,...,z_{j-1},
\zeta _j,z_{j+1},...,z_m)d\zeta _j]/\partial z_j=h_j(z)$ and \\
$\partial \mbox{ }_{\gamma _j}P^n[h_j(z_1,...,z_{j-1},
\zeta _j,z_{j+1},...,z_m)d\zeta _j]/\partial y_j=\alpha h_j(z)$
(see Formulas $2.4.1.(i,ii)$) such that $\mbox{ }_UP^n_{x_j}
h_j|^{x_j}_{x_{j,0}}$ and $\mbox{ }_UP^n_{y_j}h_j|^{y_j}_{y_{j,0}}$
as particular cases of $\gamma _j$ along axes $x_j$ and $y_j$
give the desired result.
\par {\bf 3.25. Corollary.} {\it The space $\mbox{ }_S{\bar C}^{(q+1,n-1)}
(\Omega ,{\bf K}(\alpha ))|_{\tilde {\Omega }}$ contains all locally
$z$-analytic functions on $\tilde {\Omega }$.}
\par {\bf Proof.} Mention that $1\in C^{(q,n-1)}(\Omega ,{\bf K}(\alpha ))$
and $\mbox{ }_UP^n_x1|^x_{x_0}=x-x_0$, 
$\mbox{ }_UP^n_y1|^y_{y_0}=y-y_0$,
$\mbox{ }_{\gamma _j}P^n1=\gamma _j(\beta )-\gamma _j(0)=
z_j-z_{j,0}$, where $\gamma _j\subset {\tilde {\Omega }}_j$,
hence $z_j-z_{j,0}\in \mbox{ }_S{\bar C}^{(q+1,n-1)}(\Omega ,
{\bf K}(\alpha ))$ for each $z_j$ and $z_{j,0}\in {\tilde {\Omega }}_j=
\pi _j({\tilde {\Omega }}_j)$. It is possible to take
$\gamma _j$ contained in balls $B$ such that $B\subset \Omega _j$.
Therefore, $\mbox{ }_{\gamma }P^n\sum_{l=1}^k\chi _{B_l}=
\sum_{l=1}^ka_l\mbox{ }_{\gamma _j}P^n\chi _{B_l}\in
\mbox{ }_S{\bar C}^{(q+1,n-1)}(\Omega _j,{\bf K}(\alpha ))$,
where $B_l$ are balls satisfying conditions of Lemma
$2.6.1$, $a_l\in {\bf K}(\alpha )$, $k\in \bf N$.
In view of Theorem $3.24$
each polynomial in $z$ belongs to $\mbox{ }_S{\bar C}^{(q+1,n-1)}
(\Omega ,{\bf K}(\alpha ))|_{\tilde {\Omega }}$. The using of
expansions into series by $z$ of locally $z$-analytic functions
and limits of sequences of polynomials in $z$ and Lemma $2.6.1$
leads to the conclusion that each locally $z$-analytic function
on $\tilde {\Omega }$ belongs to $\mbox{ }_S{\bar C}^{(q+1,n-1)}(\Omega ,
{\bf K}(\alpha ))|_{\tilde {\Omega }}$.
\par {\bf 3.26. Note.} From Corollary $3.25$ it follows, that a
$\mbox{ }_S{\bar C}^{(q+1,n-1)}$-manifold $M$ is locally $z$-analytic
manifold and there exists a refinement ${\tilde A}t(M)= \{ ({\tilde U}_j,
{\tilde {\phi }}_j):$ $j \} $ of $At (M)$ such that transition mappings
${\tilde {\phi }}_l\circ {\tilde {\phi }}_j^{-1}$ are $z$-analytic for
each ${\tilde U}_j\cap {\tilde U}_l\ne \emptyset $. If $f$ is
$z$-analytic on $\tilde {\Omega }$, then $f'$ is $z$-analytic on
$\tilde {\Omega }$. Therefore, there exists a family $\Upsilon $
of the cardinality $card (\Upsilon )={\sf c}:=card({\bf R})$
of all functions $f\in \mbox{ }_S{\bar C}^{(q+1,n-1)}(\Omega ,
{\bf K}(\alpha ))$ and $f$ is not $z$-analytic on $\tilde {\Omega }$,
since a locally $z$-analytic function is not necessarily $z$-analytic.
For example, take $h\in C^{(q,n-1)}(\Omega ,{\bf K}(\alpha ))$
locally $z$-analytic on $\tilde {\Omega }$ and nonanalytic on
$\tilde {\Omega }$ and put $f(z)=\mbox{ }_{\gamma }P^nh$,
where $\gamma (0)=z_0$, $\gamma (\beta )=z$, $\Omega $ is a polydisc
(see \S 3.25). Indeed, each locally plynomial in $z$ nonpolynomial
$h: \Omega \to {\bf K}(\alpha )$ and its iterated antiderivatives along
paths $h_k(z):=\mbox{ }_{\gamma }P^kh_{k-1}$, $k=1,...,n$, $h_0:=h$,
up to order $n$ fit this construction. For nonlocally compact fields
there is the theory of analytic elements \cite{escas}.
\par {\bf 3.27. Corollary.} {\it Let $\bf L$ be a non-Archimedean
field such that ${\bf K}(\alpha )\subset \bf L$ with a valuation
$|*|_{\bf L}$ extending that of ${\bf K}(\alpha )$ and let $\bf L$
be complete relative to $|*|_{\bf L}$. Suppose
$\Omega $ is a clopen compact subset in $({\bf K}\oplus \alpha {\bf K})^m$
and $Y$ is a Banach space over $\bf L$. Then $f\in
\mbox{ }_S{\bar C}^{(q+1,n-1)}
(\Omega ,Y)|_{\tilde {\Omega }}$ if and only if there exists an open subset
$W$ in ${\bf L}^m$ and a locally $z$-analytic function $F$ on $W$,
$z\in W$, such that $W\cap ({\bf K}\oplus \alpha {\bf K})^m
\supset \tilde {\Omega }$ and $F|_{\tilde {\Omega }}=f$.}
\par {\bf Proof.} The valuation group $\Gamma _{{\bf K}(\alpha )}$ is
discrete, hence $Y$ as the ${\bf K}(\alpha )$-linear space has
an orthonormal base $ \{ e_j:$ $j\in \Lambda \} $, where
$\Lambda $ is a set (see Chapter 5 \cite{roo}). Therefore,
$F: W\to Y$ has the decomposition $F(z)=\sum_{j\in \Lambda }
F_j(z)e_j$, where $F_j: W\to {\bf K}(\alpha )$. Since $F$ is
locally $z$-analytic, then $f$ is locally $z$-analytic on
$\tilde {\Omega }$ and in accordance with Corollary $3.25$
$f\in \mbox{ }_S{\bar C}^{(q+1,n-1)}(\Omega ,Y)|_{\tilde {\Omega }}$.
Vice versa, if $f\in \mbox{ }_S{\bar C}^{(q+1,n-1)}
(\Omega ,Y)|_{\tilde {\Omega }}$, then by Theorem $2.7.2$
$f$ is locally $z$-analytic on $\tilde {\Omega }$, consequently,
for each $\zeta \in \tilde {\Omega }$ there exists a ball
$B({\bf K}(\alpha ),\zeta ,R(\zeta ))$ with $0<R(\zeta )<\infty $
on which the power series $2.7.2.(2)$ is uniformly convergent,
that is, $\lim_{|k|\to 0} |a_k|_{\bf L}R^{|k|}=0$, hence
this series is uniformly convergent on $B({\bf L},\zeta ,
R(\zeta ))$ also. Put $W=\bigcup_{\zeta \in \tilde {\Omega }}
B({\bf L},\zeta ,R(\zeta ))$.
\par {\bf 3.28. Definition and Note.} Let $\bf L$, $\Omega =
\Omega (f)$, $W=W(f)$ be satisfying conditions of \S 3.27
with $m=1$.
Let also $T\in L(Y)$ be a bounded $\bf L$-linear operator on a Banach
space $Y$ over $\bf L$ with a nonvoid spectrum
$\sigma (T):= \{ b\in {\bf L}:$ $(bI-T)$ $\mbox{is not invertible in}$
$L(Y) \} $ (see Chapter 6 in \cite{roo}), where $L(X,Y)$ is the Banach
space of all bounded $\bf L$-linear operators $T: X\to Y$
for Banach spaces $X$ and $Y$ over $\bf L$, $ \| T \| :=
\sup_{0\ne x\in X} \| Tx \| / \| x \| $, $L(Y):=L(Y,Y)$.
Suppose in addition, that for each
$z\in W$ with $dist (z,\Omega )<\infty $ there exist
$R\ge dist (z,\Omega )$ and $\zeta \in \Omega $ such that
$B({\bf L},z,R)\subset W$ and $B({\bf L},z,R)\cap \Omega =
B({\bf K}\oplus \alpha {\bf K},\zeta ,R)\subset \Omega $.
Denote by ${\cal F}(T)$ a family of all functions $f$ with
$\psi _f\in \mbox{ }_S{\bar C}^{(q+1,n-1)}(\omega _{\epsilon },
{\bf L})$, where $W$ is a clopen neighborhood of $\sigma (T)$,
$W=W(f)$, $\Omega =W\cap ({\bf K}\oplus \alpha {\bf K})\ne \emptyset $,
$0<dist (\partial \Omega ,\sigma (T)):=\inf_{z\in \partial \Omega }
dist (z, \sigma (T))$, $dist (z,G):=\inf_{y\in G} |z-y|$ for
$G\subset \bf L$ and $z\in \bf L$, $0\le q\in \bf Z$,
$1\le n\in \bf N$, $\psi _f(\eta ):=f(z+Exp (\eta ))$,
$\omega :=\omega (z):= \{ \eta \in {\bf K}(\alpha ):$
$z+Exp (\eta )\in \Omega \} $, $z\in \Omega $,
$\omega _{\epsilon }=\omega \setminus Log (B({\bf K}(\alpha ),z,\epsilon ))$,
$\epsilon =\epsilon _j$, $\epsilon _j>0$ for each $j\in \bf N$,
$\lim_{j\to \infty }\epsilon _j=0$, there exists a locally
$z$-analytic function $\Psi _f$
on $W$ such that $\Psi _f|_{\tilde {\Omega }}=\psi _f$
(see \S 3.27). Put
\par $(1)$ $f(T)=C(\alpha )^{-1}\mbox{ }_{\partial {\Omega }}P^n[
f(\zeta )R(\zeta ;T)d\zeta ]$, where $R(\zeta ;T)=(\zeta I-T)^{-1}$
for $\zeta \in \rho (T):={\bf L}\setminus \sigma (T)$
and the antiderivative is supposed to be convergent in the strong
operator topology sence, that is, $\mbox{ }_{\partial \Omega }
P^n[f(\zeta )R(\zeta ;T)y d\zeta ]$ converges for each $y\in Y$.
There are others definitions of spectral sets
(see Chapter 6 \cite{roo}), but this one is used here.
\par {\bf 3.29. Theorem.} {\it Let $\sigma (T) \ne \emptyset $,
$\sigma (T)\subset \bf L$, $f, g \in {\cal F}(T)$, $a, b\in \bf L$
(see \S 3.28). Then
\par $(i)$ $af+bg\in {\cal F}(T)$ and $aF(T)+bg(T)=(af+bg)(T)$;
\par $(ii)$ $fg\in {\cal F}(T)$ and $f(T)g(T)=(fg)(T)$;
\par $(iii)$ if $f(z)=\sum_{k=0}^{\infty }a_kz_k$ on $W(f)$ such that
$W(f)\supset \sigma (T)$, then $f(T)=\sum_{k=0}^{\infty }a_kT^k$.}
\par {\bf Proof.} Definition 3.28 is correct, since
$xI-T$ is invertible in $L(Y)$ for each
$x\in \rho (T):={\bf L}\setminus \sigma (T)$,
hence $r_{\sigma }(T):=\sup_{x\in \sigma (T)}|x|\le \| T \| $,
where $\rho (T)$ is open in $\bf L$
and $R(x;T)$ is locally $x$-analytic on $\rho (T)$
(see Chapter VII in \cite{dansch} and Chapter 6 in \cite{roo}).
\par $(i)$. It follows from Definition 3.28 and Corollary $3.27$.
\par $(ii)$. In view of Corollary $3.27$ and Theorem $3.24$
$fg\in {\cal F}(T)$, since $W(f)\cap W(g)=:W(fg)\supset \sigma (T)$.
Without loss of generality take $\Omega (f)$ encompassed by $\partial
\Omega (g)$ shrinking $\Omega (f)$ a little if necessary such that
$W(f)\supset \sigma (T)$, $W(f)\subset W(g)$. Then
\par $f(T)g(T)=C(\alpha )^{-2}\mbox{ }_{\zeta \in \partial \Omega (f)}
P^n[f(\zeta )R(\zeta ;T)d\zeta ]\mbox{ }_{\kappa \in \partial \Omega (g)}
P^n[g(\kappa )R(\kappa ;T)d\kappa ]$ \\
$=C(\alpha )^{-2}\mbox{ }_{\kappa \in \partial \Omega (g)}P^n[
\mbox{ }_{\zeta \in \partial \Omega (f)}P^n[f(\zeta )g(\kappa ) \{
R(\zeta ;T)R(\kappa ;T) \} d\zeta \} d\kappa ]$. \\
On the other hand, $R(\zeta ;T)R(\kappa ;T)=(R(\zeta ;T)-R(\kappa ;T))
(\kappa -\zeta )^{-1}.$ Therefore,
\par $(1)$ $f(T)g(T)=C(\alpha )^{-2}\mbox{ }_{\zeta \in \partial \Omega (f)}
P^n[f(\zeta )R(\zeta ;T) \{ \mbox{ }_{\kappa \in \partial \Omega (g)}
P^n[g(\kappa )(\kappa -\zeta )^{-1}d\kappa ] \} d\zeta ]$ \\
$-C(\alpha )^{-2}\mbox{ }_{\kappa \in \partial \Omega (g)}P^n[
g(\kappa )R(\kappa ;T) \{ \mbox{ }_{\zeta \in \partial \Omega (f)}
P^n[f(\zeta )(\kappa -\zeta )^{-1}d\zeta ] \} d\kappa ]$. \\
The second term on the right hand side of $(1)$ is zero,
since $\partial \Omega (f)$ is encompassed by $\partial \Omega (g)$,
$\kappa \in \partial \Omega (g)$, $\zeta \in \partial \Omega (f)$
(see Formulas $2.7.2.(2-4)$). Hence
\par $f(T)g(T)=C(\alpha )^{-1}\mbox{ }_{\zeta \in \partial \Omega (f)}
P^n[f(\zeta )g(\zeta )R(\zeta ;T)d\zeta ]=(fg)(T)$.
\par $(iii)$. It follows from Definition $3.28$ and Formulas $2.7.2.(2-4)$
applied to $f(\zeta )R(\zeta ;T)y$ for each $y\in Y$.
\par {\bf 3.30. Theorem.} {\it Let $\sigma (T)\ne \emptyset $,
$\sigma (T)\subset \bf L$, $f\in {\cal F}(T)$ (see \S 3.28).
Then $f(\sigma (T))=\sigma (f(T))$.}
\par {\bf 3.31. Theorem.} {\it Let $\sigma (T)\ne \emptyset $,
$\sigma (T)\subset \bf L$, $f\in {\cal F}(T)$,
$g\in {\cal F}(f(T))$ (see \S 3.28) and $h(z):=g(f(z))$ for each
$z\in f^{-1}[W(g)\cap f(W(f))].$ Then $h\in {\cal F}(T)$ and
$h(T)=g(f(T))$.}
\par {\bf Proof.} Theorem $3.30$ follows from Theorem $3.29$
analogously to Theorems $VII.3.10$ \cite{dansch} and $3.3.6$
\cite{kadrin}. The function $f$ is locally $z$-analytic on $W(f)$,
$g$ is locally $z$-analytic on $W(g)$, hence $h$ is locally $z$-analytic on
$f^{-1}[W(g)\cap f(W(f))].$ In view of Theorem $3.30$ $\sigma (f(T))
\subset f(W(f))\cap W(g)$, hence $h$ is defined on open $W(h)$ such that
$W(h)\supset \sigma (T)$. Without loss of generality take $W(g)\supset
f(W(f))$. Put
\par $S(\kappa )=C(\alpha )^{-1}\mbox{ }_{\zeta \in
\partial \Omega (f)}P^n[R(\zeta ;T)(\kappa -f(\zeta ))^{-1}d\zeta ]$,
then in accordance with Theorems $3.29$ and $2.4.6$ (applied on
pieces of $\Omega (f)$ affine homotopic to points)
$(\kappa I-T)S(\kappa )=S(\kappa )(\kappa I-T)=I$, consequently,
$S(\kappa )=R(\kappa ;T)$. Therefore,
\par $g(f(T))=C(\alpha )^{-1}\mbox{ }_{\partial \Omega (g)}P^n[
g(\kappa )R(\kappa ;f(T))d\kappa ]$ \\
$=-C(\alpha )^{-2}\mbox{ }_{\partial \Omega (g)}P^n[\mbox{ }_{\partial
\Omega (f)}P^n \{ g(\kappa )R(\zeta ;T)(\kappa -f(\zeta ))^{-1}d\zeta \}
d\kappa ]$ \\
$=C(\alpha )^{-1}\mbox{ }_{\partial \Omega (f)}P^n[R(\zeta ;T)g(f(\zeta ))
d\zeta ]=h(T)$.
\par {\bf 3.32. Proposition.} {\it Let $f_k\in {\cal F}(T)$ for each
$k\in \bf N$ (see \S 3.28) and there exists a clopen subset $W$
in $\bf L$ such that $\sigma (T)\subset W\subset \bigcap_{n=1}^{\infty }
W(f_n).$ If $f_k$ converges to $f$ uniformly on $W$,
then $f_n(T)$ converges to $f(T)$ uniformly on each totally
bounded subset in $Y$.}
\par {\bf Proof.} There exists a sequence $C(\alpha )^{-1}
\mbox{ }_{\partial \Omega }P^n[f_k(\zeta )R(\zeta ;T)d\zeta ]$
in $L(X,Y)$ in the topology of pointwise convergence, where
$L(X,Y)$ denotes the Banach space of continuous $\bf L$-linear
operators $S: X\to Y$ for two Banach spaces $X$ and $Y$
over $\bf L$. In view of Theorem $(11.6.3)$ and Example
$11.202.(g)$ \cite{nari} this sequence converges to a $\bf L$-linear
operator on $Y$ uniformly on each totally bounded subset in $Y$.
\par {\bf 3.33. Definition.} A point $z_0\in \sigma (T)$ is called
an isolated point of a spectrum $\sigma (T)$, if there exists
a neighborhood $U$ of $z_0$ such that $\sigma (T)\cap U= \{ z_0 \} $,
where $U$ satisfies the same conditions of \S 3.28 as $W$.
An isolated point $z_0\in \sigma (T)$ is called a pole of an operator
$T$ or a pole of a spectrum, if a mapping $R(\zeta ;T)$ has a pole
at $z_0$. An order $j(z_0)$ of a pole $z_0$ is an order of $z_0$ as a
pole of $R(\zeta ;T)$.
\par {\bf 3.34. Theorem.} {\it Let $f, g \in {\cal F}(T)$
(see \S 3.28). Then $f(T)=g(T)$ if and only if $f(\zeta )=g(\zeta )$ on a
clopen $W$ such that $\sigma (T)\setminus \bigcup_{l\in \Lambda }
\{ z_l \} \subset W\subset \bf L$, where $z_l\in \Omega \subset
{\bf K}\oplus \alpha {\bf K}$ is a pole for each $l\in \Lambda $,
$\Lambda $ is a finite set and $(f-g)$ at $z_l$ has zero of order
not less than $j(z_l)$ for each $l=1,...,k$. }
\par {\bf Proof.} Without loss of generality take $g=0$ and let $f=0$
on $W\setminus \bigcup_{l\in \Lambda } \{ z_l \} $. Then due to
Theorem $2.4.6$ (applied on each piece affine homotopic to a point)
\par $f(T)=C(\alpha )^{-1}\sum_{l\in \Lambda }\mbox{ }_{\partial
B_l}P^n[f(\zeta )R(\zeta ;T)d\zeta ]$, where $B_l:=B({\bf K}\oplus
\alpha {\bf K},z_l,R_l)$, $0<R_l<\infty $ and $B_l\cap \sigma (T)=
\{ z_l \} $ for each $l\in \Lambda $. Since $f(\zeta )R(\zeta ;T)$ is
regular on $B_l$, then by Theorem $2.7.14$ $f(T)=0$. Vice versa, let
$f(T)=0$, then by Theorem $3.30$ $f(\sigma (T))=0.$ The set $\sigma
(T)\cap ({\bf K}\oplus \alpha {\bf K})$ is compact and it can
be covered by a finite union of balls $B({\bf K}\oplus \alpha {\bf K},
\zeta _j,R_j)$, $0<R_j<\infty $. If $B({\bf K}\oplus \alpha {\bf K},
\zeta _j,R_j)\cap \sigma (T)$ is infinite, then for each limit point
$x$ of the latter set there exists a clopen neighborhood $V_x$ on
which $f|_{V_x}=0$ (see Theorem $2.7.7$). Therefore, $\sigma (T)\cap
(({\bf K}\oplus \alpha {\bf K})\setminus \bigcup_xV_x)$ consists
of a finite number of isolated points $\{ \lambda _l:$ $l=1,...,k \} $,
since $\Omega \supset ({\bf K}\oplus \alpha {\bf K})\cap \sigma (T)$,
$\Omega $ is compact. Let $f$ is not zero on any neighborhood
of $\lambda _1$. Since $\lambda _1\in \sigma (T)$ and $f(\sigma (T))=
\{ 0 \} $, then $f$ has a zero of finite order $j$, hence $g_1(z)
=(\lambda _1-z)^j/f(z)$ is locally $z$-analytic on a neighborhood
of $\lambda _1$. From the proof of Theorem $3.24$ it follows, that
\par $(1)$ $R(\zeta ;T)=\sum_{m=-\infty }^{\infty }
a_m(\lambda _1-\zeta )^m$ \\
on $B(\epsilon ):=B({\bf K}\oplus \alpha {\bf K},\lambda _1,\epsilon )$
for a sufficiently small $0<\epsilon <\infty $, where
\par $(2)$ $a_{-m}=-C(\alpha )^{-1}\mbox{ }_{\partial B(\epsilon )}
P^n[(\lambda _1-\zeta )^{m-1}R(\zeta ;T)d\zeta ]$ 
$=-(\lambda _1I-T)^{m-1}h(T)$, \\
$h(T)$ denotes a function equal to $1$ on $B(\epsilon )$ and zero
on a neighborhood of $\lambda _l$ for each $l\ne 1$ such that
$\psi (\eta )=h(z+Exp(\eta ))$ satisfies $3.28$, that is,
possible due to Lemma $2.6.1$ and Corollary $3.27$, since $Exp
(\eta )$ is locally $\eta $-analytic. Then $a_{-m-1}=
-(\lambda _1I-T)^mh(T)=0$ for each $m\ge j$.
\par {\bf 3.35. Definition and Note.} A subset $V$ of $\sigma (T)$
clopen in $\sigma (T)$ is called a spectral set if it has a clopen
neighborhood $W_V$ satisfying the same conditions of \S 3.28
as $W$ and $W_V\cap (\sigma (T)\setminus V)=\emptyset $.
In accordance with Lemma $2.6.1$ and Theorem $3.24$ consider
$f\in {\cal F}(T)$ such that $f|_V=1$ and $f|_{\sigma (T)\setminus V}=0$,
which is possible due to Corollary $3.27$, since $Exp (\eta )$
is locally $\eta $-analytic. Put $E(V;T):=f(T)$. In view of Theorem
$3.34$ $E(V;T)$ depends on $V$, but not on a concrete choice of $f$
from its definition. If $V\cap \sigma (T)=\emptyset $, put $E(V;T)=0.$
Write also $E(z;T):=E( \{ z \} ;T)$ for a singleton $ \{ z \} $.
An index $j=j(z)$ of $z\in \bf L$ is the smallest integer $j$ such that
$(zI-T)^jy=0$ for each $y\in Y$ with $(zI-T)^{j+1}y=0.$
\par {\bf 3.36. Theorem.} {\it Let $T$, $W$, $\Omega $, ${\bf K}(\alpha )$
be as in \S \S 3.27, 3.28. If $z_0$ is a pole of $T$ of order $j$, then
$z_0\in \Omega $ has the index $j$. An isolated point $z_0\in \sigma (T)$
is a pole of order $j$ if and only if
\par $(i)$ $(z_0I-T)^jE(z_0;T)=0$, $(z_0I-T)^{j-1}E(z_0;T)\ne 0$.}
\par {\bf Proof.} In view of Formulas $3.34.(1,2)$ $z_0$ is a pole of order
$j$ if and only if $(i)$ is satisfied, since $a_{-m-1}=-(z_0I-T)^mE(z_0;T).$
The rest of the proof is analogous to that of Thereom $VII.3.18$
\cite{dansch} due to Corollary $3.27$ and Theorem $3.24$.
\par In view of Theorem $3.29.(ii)$
\par $E(V;T) E(V;T) =E(V;T)$ for each spectral set $V$, \\
that is, $E(V;T)$ is the projection operator on $Y$
(see Chapter 3 \cite{roo}).
\par {\bf 3.37. Theorem.} {\it Let $f\in {\cal F}(T)$ (see \S 3.28)
and let $V$ be a spectral set of $f(T)$. Then $\sigma (T)\cap
f^{-1}(V)$ is the spectral set of $T$ and $E(V;f(T))=E(f^{-1}(V);T)$.}
\par {\bf Proof.} Let $h_V\in {\cal F}(T)$ such that $h_V(z)=1$
on a neighborhood $V_1$ of $V$, $h_V(z)=0$ on a neighborhood $V_2$
of $\sigma (f(T))\setminus V_1$, where $V_1\cap V_2=\emptyset $,
which is possible due to Theorem $3.24$, Corollary $3.27$ and
Lemma $2.6.1$, since $Exp (\eta )$ is locally $\eta $-analytic. Then
$h_V(f(T))=E(V;f(T)).$ In view  of Theorem $3.30$ $\sigma (T)=
f^{-1}(V)\cup f^{-1}(\sigma (f(T))\setminus V)$, where
$f^{-1}(V)\cap f^{-1}(\sigma (f(T))\setminus V)=\emptyset $.
Since $f$ is continuous, then $f^{-1}(V)$ and $f^{-1}(\sigma (T)\setminus
V)$ are clopen in $\sigma (T)$. Therefore, $\sigma (T)\cap f^{-1}(V)=:
\Upsilon $ is the spectral set of $T$. Put $t_{\Upsilon }(z):=
h_V(f(z))$, then $E(\Upsilon ;T)=t_{\Upsilon }(T)$, since
$t_{\Upsilon }\in {\cal F}(T)$ due to Corollary $3.27$.
From Theorem $3.31$ it follows, that $E(V;f(T))=E(\Upsilon ;T)=
E(f^{-1}(V);T)$.
\par {\bf 3.38. Remark.} In the non-Archimedean case the Gelfand-Naimark
Theorem ($IX.3.7$ \cite{dansch}) is not true (see Chapter 6 \cite{roo}).
Therefore, the existence of the projection operator $E(V;T)$ for
each spectral set $V$ does not imply a spectral projection-valued
measure decomposition of $T$ (see also \cite{luddia}).
Here is considered a particular class of operators
satisfying conditions of \S 3.28 for which the operator $E(V;T)$
is defined for each spectral set $V$, $V\subset \sigma (T)$.
Put $Y_V:=E(V;T)Y$. In view of Theorem $3.29.(ii)$ and \S 3.35
$TY_V\subset Y_V$, where $Y_V$ is the $\bf L$-linear subspace
in $Y$, since $E(V;T)$ is $\bf L$-linear, denote $T_V:=T|_{Y_V}$.
\par {\bf 3.39. Theorem.} {\it Let $V$ be a spectral set of
$\sigma (T)\ne \emptyset $ (see \S 3.28). Then $\sigma (T_V)=V$.
If $f\in {\cal F}(T)$, then $f\in {\cal F}(T_V)$ and $f(T_V)=
f(T)_V$. A point $z_0\in V\cap \Omega $ is the pole of $T$ of order
$j$ if and only if $z_0\in \Omega $ is the pole of $T_V$ of order $j$.}
\par {\bf Proof.} Take a marked point $z\in V$ and suppose $z\notin
\sigma (T_V)$. In view of Corollary $3.27$ there exists a function
$g\in {\cal F}(T)$ such that $g|_{V_1}=0$ on a neighborhood $V_1$ of $V$
and $g(\zeta )=(z_0-\zeta )^{-1}$ for each $\zeta \in V_2$,
where $V_2$ is open in $\bf L$,
$V_1\cap V_2=\emptyset $, $V_2\supset \sigma (T)\setminus V$. In view of
Theorem $3.29.(ii)$ $g(T)(zI-T)=(zI-T)g(T)=I-E(V;T)$. Then
$V\subset \sigma (T_V)$ as in Theorem $VII.3.20$ \cite{dansch}.
\par Vice versa, let $z\notin V$. Consider $h\in {\cal F}(T)$
(see \S 3.28) such that $h(\zeta )|_{V_1}=(z-\zeta )^{-1}$
and $h|_{V_2}=0$, where  $V_1$ is chosen such that $z\notin V_1$,
$V_1$ is a neighborhood of $V$, $V_2$ is as above.
Then by Theorem $3.29.(ii)$ $h(T)(zI-T)=(zI-T)h(T)=E(V;T)$.
Therefore, $h(T)_V(zI_V-T_V)=(zI_V-T_V)h(T_V)=I_V$, since
$z\notin \sigma (T_V)$, consequently, $\sigma (T_V)\subset V$ and
$R(z;T_V)=R(z;T)_V$. Take $f\in {\cal F}(T)$ and a neighborhood
$W$ of $\sigma (T)$ as in \S 3.28. Then \\
$f(T)_V=C(\alpha )^{-1}\mbox{ }_{\partial \Omega }P^n[
f(z)R(z;T)dz]_V=C(\alpha )^{-1}\mbox{ }_{\partial \Omega }P^n[
f(z)R(z;T_V)dz]=f(T_V)$ \\
and $E(z;T)E(V;T)=E(z;T)$ for each $z\in V$, hence
$(zI-T)^kE(z;T)=(zI_V-T_V)^kE(z;T)$ for each $k\in \bf N$.
In view of Theorem $3.36$ $z_0\in \Omega \cap V$ is a pole of $T$
of order $j$ if and only if it is a pole of $T_V$ of order $j$.
\par {\bf 3.40. Corollary.} {\it The mapping $E\mapsto E(V;T)$
is the isomorphism of the algebra $\Upsilon $ of all clopen spectral subsets
$V$ of $\sigma (T)$ satisfying conditions of \S 3.28 on
the Boolean algebra $\{ E(V;T):$ $V\in \Upsilon \} $.}
\par {\bf Proof.} In view of Theorem $3.29$ the mapping
$V\mapsto E(V;T)$ is the homomorphism. If $E(V;T)=0$, then $Y_V=0$
and $\sigma (T_V)=\emptyset $, hence $V=\sigma (T_V)=\emptyset $
by Theorem $3.39$. If $V_1, V_2\in \Upsilon $, then evidently
$W_{V_1}\cup W_{V_2}$ and $W_{V_1}\cap W_{V_2}$ (for $V_1\cap V_2\ne
\emptyset $) satisfy conditions of \S 3.28 as $W$.
Consider $\sigma (T)\setminus V$ for $V\in \Upsilon $, then
$W_V\cap (\sigma (T)\setminus V)=\emptyset $ (see \S 3.25),
hence $W\setminus W_V$ satisfies conditions of \S 3.28 as $W$, since each two
balls in $\bf L$ are either disjoint or one of them is contained
in another. Therefore, $\Upsilon $ is the Boolean algebra
and hence $\{ E(V;T):$  $V \in \Upsilon \} $ is the Boolean algebra.
\par {\bf 3.41. Note.} In sections 3.28-3.40 it can be taken
the generalization instead of $\Omega $ for a manifold $M$ which is
$\mbox{ }_SC^{(q+1,n-1)}$-diffeomorphic with $\Omega $.

Address: Theoretical Department,
\par Institute of General Physics,
\par Russian Academy of Sciences,
\par Str. Vavilov 38,
\par Moscow 119991 GSP-1, Russia
\par e-mail: ludkovsk@fpl.gpi.ru

\begin{thebibliography}{99}
\bibitem{ami} Y. Amice. Interpolation p-adique. Bull. Soc. Math. France
   {\bf 92}(1964), 117-180.
\bibitem{berz} M. Berz. Cauchy theory on Levi-Civit\'a fields.
in: Contemporary Mathem. {\bf 319} (2003), 39-52.
Ultrametric functional analysis. Seventh international conference
on $p$-adic functional analysis. June 17-21, 2002. Nijmegen.
Ed. W.H. Schikhof, et. al.
\bibitem{boum} N. Bourbaki. Differentiable and analytical manifolds
   (Moscow: Mir,1975).
\bibitem{dansch} N. Dunford, J.T. Schwartz. Linear operators.
V.V. 1, 2  (New York: Interscience Publishers, 1958, 1963).
\bibitem{eng} R. Engelking. General topology (Moscow: Mir,1986).
\bibitem{escas} A. Escassut. Analytic elements in $p$-adic
analysis (Singapore: World Scientific, 1995).
\bibitem{freput} J. Fresnel, M. van der Put. "G\'eom\'etrie
analytique rigide et applications" (Boston: Birkh\"auser, 1981).
\bibitem{henlei} G. Henkin, J. Leiterer. Theory of functions
on complex manifolds (Basel: Birkh\"auser-Verlag, 1984).
\bibitem{kadrin} R.V. Kadison, J.R. Ringrose. Fundamentlas
of the theory of operator algebras (New York: Academic Press, 1983).
\bibitem{kobl} N. Koblitz. $p$-adic numbers, $p$-adic analysis and
zeta functions (New York: Springer-Verlag, 1977).
\bibitem{koblbr} N. Koblitz. $p$-adic analysis: a short course
on recent work. Lond. Math. Soc., Lecture Notes Series. {\bf 46} (1980).
\bibitem{luambp00} S.V. Ludkovsky. Quasi-invariant measures
on non-Archimedean groups and semigroups of loops and paths,
their representations. Ann. Math. B. Pascal. {\bf 7: 2} (2000),
19-53, 55-80.
\bibitem{lutmf99} S.V. Ludkovsky. Measures on groups of diffeomorphisms
of non-Archimedean manifolds, representations of groups and their 
applications. Theoret. and Math. Phys. {\bf 119: 3} (1999), 698-711.
\bibitem{lustpnam} S.V. Ludkovsky. Non-Archimedean stochastic processes
on non-Archimedean manifolds. Los Alamos Nat. Lab. preprint
40 pages, {\bf math.GM/0212296}, December 2002.
\bibitem{lujmsq2} S.V. Ludkovsky. Quasi-invariant and pseudo-differentiable
measures with values in non-Archimedean fields
on a non-Archimedean Banach space. J. Mathem. Sci. {\bf 112} (2002)
(previous variants: Intern. Centre for Theor.
Phys. Preprint  {\bf IC/96/210},  October 1996;
Los Alamos Nat. Lab. Preprint {\bf math.GM/0106170}).
\bibitem{luseam3} S.V. Ludkovsky. A structure and representations
of diffeomorphism groups of non-Archimedean manifolds.
Southeast Asian Bull. Math. {\bf 26} (2003), 975-1004.
\bibitem{lu8} S.V. Ludkovsky. Embeddings of non-Archimedean
Banach manifolds into non-Archimedean Banach spaces. Uspek. Mat. Nauk.
{\bf 53: 5} (1998), 241-242. 
\bibitem{luddia} S.V. Ludkovsky, B. Diarra. Spectral integration
and spectral theory for non-Archimedean Banach spaces.
Intern. J. of Math. and Mathem. Sciences {\bf 31: 7} (2002), 421-442.
\bibitem{nari} L. Narici, E. Beckenstein. Topological vector
spaces (New York: Marcel Dekker Inc., 1985).
\bibitem{put68} M. van der Put. Alg\'ebres de fonctions continues
$p$-adiques. Indagations Math. {\bf A 30: 4} (1968), 401-420.
\bibitem{roo} A.C.M. van Rooij. Non-Archimedean functional analysis (New
    York: Marcel Dekker Inc., 1978).
\bibitem{sch1} W.H. Schikhof. Ultrametric calculus (Cambr.:
Cambr. Univ. Press, 1984).
\bibitem{sch2} W.H. Schikhof. Non-Archimedean calculus (Nijmegen: Math.
    Inst., Cath. Univ., Report {\bf 7812}, 1978).
\bibitem{sch3} W.H. Schikhof. The set of derivatives in a non-Archimedean
    field. Mathem. Annalen. {\bf 216}(1975), 67-70.
\bibitem{shabat} B.V. Shabat. An introduction to the complex analysis
(Moscow: Nauka, 1985).
\bibitem{span} E. Spanier. Algebraic topology (Moscow: Mir, 1971).
\bibitem{wei} A. Weil. Basic number theory (Berlin: Springer-Verlag, 1973).
\bibitem{weil40} A. Weil. L'integration dans les groupes topologiques
et ses applications. Actual. Scient. et Ind. {\bf 869} (Paris:
Herman, 1940).
\end{thebibliography}
\end{document}